\documentclass[a4paper,12pt,reqno]{amsart}

\usepackage{amsmath,amssymb,amsthm}
\usepackage{latexsym}
\usepackage{color}
\usepackage{graphicx}
\usepackage{mathrsfs}
\usepackage{enumerate}
\usepackage[abbrev]{amsrefs}
\usepackage[T1]{fontenc}
\usepackage{mathtools}
\mathtoolsset{showonlyrefs=true}

\allowdisplaybreaks[4]

\theoremstyle{plain}
\newtheorem{thm}{Theorem}[section]
\newtheorem{lemm}[thm]{Lemma}

\theoremstyle{definition}
\newtheorem{df}[thm]{Definition}
\newtheorem{rem}[thm]{Remark}

\makeatletter

\@addtoreset{equation}{section}
\makeatother

\renewcommand{\div}{\operatorname{div}}
\newcommand{\dB}{\dot{B}}

\newcommand{\supp}{\operatorname{supp}}

\newcommand{\ha}{\widehat{a}}
\newcommand{\hu}{\widehat{u}}
\newcommand{\hV}{\widehat{V}}

\newcommand{\hf}{\widehat{f}}
\newcommand{\hg}{\widehat{g}}

\renewcommand{\leq}{\leqslant}
\renewcommand{\geq}{\geqslant}
\newcommand{\tC}{\widetilde{C}}

\newcommand{\n}[1]{{\left\|#1\right\|}}

\newcommand{\Mp}[1]{\left\{#1\right\}}
\renewcommand{\sp}[1]{\left(#1\right)}

\begin{document}
\title[The compressible Navier--Stokes--Coriolis system]
{Global strong solutions to the compressible Navier--Stokes--Coriolis system for large data}
\author[M.~Fujii]{Mikihiro Fujii}
\address[M.~Fujii]{Graduate School of Science, Nagoya City University, Nagoya, 467-8501, Japan}
\email[Corresponding author]{fujii.mikihiro@nsc.nagoya-cu.ac.jp}
\author[K.~Watanabe]{Keiichi Watanabe}
\address[K.~Watanabe]{{School of General and Management Studies, Suwa University of Science, 5000-1, Toyohira, Chino, Nagano 391-0292, Japan}}
\email{watanabe\_keiichi@rs.sus.ac.jp}
\begin{abstract}
We consider the compressible Navier--Stokes system with the Coriolis force on the $3$D whole space.
In this model, the Coriolis force causes the linearized solution to behave like a $4$th order dissipative semigroup $\{ e^{-t\Delta^2} \}_{t>0}$ with slower time decay rates than the heat kernel, which creates difficulties in nonlinear estimates in the low-frequency part and prevents us from constructing the global strong solutions by following the classical method.
On account of this circumstance, the existence of unique global strong solutions has been open even in the classical Matsumura--Nishida framework.
In this paper, we overcome the aforementioned difficulties and succeed in constructing a unique global strong solution in the framework of the scaling critical Besov space.
Furthermore, our result also shows that the global solution is constructed for arbitrarily {\it large} initial data provided that the speed of the rotation is high and the Mach number is low enough by focusing on the dispersive effect due to the mixture of the Coriolis force and acoustic wave.
\end{abstract}
\maketitle

\section{Introduction}\label{sec:intro}
We are concerned with the  {initial value} problem of the $3$D compressible Navier--Stokes system with the
Coriolis force {in} {the whole space}:
\begin{align}\label{eq:NSC-1}
	\begin{dcases}
		\partial_t \rho + \div ( \rho u ) = 0, 
		& 
		t > 0, x \in \mathbb{R}^3,\\
		\begin{aligned}
			\rho \left( \partial_t u + ( u \cdot \nabla ) u + \Omega ( e_3 \times u ) \right)
			&
			+
			\dfrac{1}{\varepsilon^2} \nabla P ( \rho ) \\
			={}
			&
			\mu \Delta u + ( \mu + \mu' ) \nabla \div u,
		\end{aligned}
		&
		t > 0, x \in \mathbb{R}^3,\\
		\rho(0,x) = \rho_0(x), \quad u(0,x) = u_0(x),
		&
		x \in \mathbb{R}^3.
	\end{dcases}
\end{align}
Here, $\rho=\rho(t,x):(0,\infty)\times \mathbb{R}^3 \to (0,\infty)$ and {$u=u(t,x):(0,\infty)\times \mathbb{R}^3 \to \mathbb{R}^3$} stand for the unknown density and velocity {field} of the fluid, respectively,  {whereas the given initial data are denoted by} $\rho_0=\rho_0(x):\mathbb{R}^3 \to (0,\infty)$ and $u_0=u_0(x):\mathbb{R}^3 \to \mathbb{R}^3$.
For the viscosity term, the constants $\mu$ and $\mu'$ stand for the shear and bulk viscosity coefficients, respectively, which are assumed to satisfy the ellipticity conditions
$\mu > 0$ and $\nu : = 2 \mu + \mu' > 0$.
The constants $\Omega \in \mathbb{R}$ and $\varepsilon > 0$ represent the speed of rotation and Mach numbers, respectively.
The term $\Omega (e_3 \times (\rho u))$ represents the Coriolis force on the flow due to rotation of the fluid around the vertical unit vector $e_3 = (0,0,1)^\top$.
Throughout of this paper, we focus on a solution $(\rho, u)$ around the constant equilibrium state $(\rho_\infty, 0)$  with some positive constant $\rho_\infty$.
We assume that the pressure $P$ is a suitably smooth function of $\rho$ satisfying $c_\infty := P' (\rho_\infty) > 0$.
System~\eqref{eq:NSC-1} arises in studies of geophysical flows and is physically well-justified at mid-latitude regions. Here, we suppose that the centrifugal force balances with the geostrophic force.
We refer to~\cite{Cu-Ro-Be-11} and references therein for further geophysical backgrounds.


The aim of this paper is to show the global existence of a unique large solution to \eqref{eq:NSC-1} in the scaling critical Besov spaces framework\footnote{{This means that we solve the system \eqref{eq:NSC-1} in the homogeneous Besov spaces framework whose norms are invariant under the scaling associated to the system. See \cite{Fuj-Wat-25}*{(1.6)} for the detail on the scaling invariance.}}.
More precisely, for arbitrarily large initial {disturbances} around the constant equilibrium state in the critical Besov spaces,
we prove that there exists a unique global strong solution of~\eqref{eq:NSC-1},
provided that the speed of rotation is high and the Mach number is low enough. 


We survey the literature associated with our study briefly\footnote{ A detailed exposition of the historical background of system \eqref{eq:NSC-1} may be found in the previous paper \cite{Fuj-Wat-25}.} before describing our main result.
Concerning either incompressible or non-rotational case, a lot of contributions have been made.
%
%
{ 
Considering the low Mach number limit $\varepsilon \to +0$ in System~\eqref{eq:NSC-1}, we formally obtain the following incompressible Navier--Stokes--Coriolis system:}
\begin{align}\label{eq:incomp-coriolis}
	\begin{dcases}
		\partial_t u - \mu \Delta u + \Omega( e_3 \times u ) + ( u \cdot \nabla) u + \nabla p = 0, & \qquad t > 0, x \in \mathbb{R}^3,\\
		\div u = 0, & \qquad t \geqslant 0, x \in \mathbb{R}^3,\\
		u(0,x) = u_0(x), & \qquad x \in \mathbb{R}^3.
	\end{dcases}
\end{align}
Mathematical analysis of {the incompressible Navier--Stokes--Coriolis system was started} by Babin, Mahalov, and Nicolaenko \cites{Ba-Ma-Ni-97,Ba-Ma-Ni-00,Ba-Ma-Ni-01}.
{
Later, Chemin, Desjardins, Gallagher, and Grenier \cites{Ch-De-Ga-Gr-02,Ch-De-Ga-Gr-06} and Charve \cite{Ch-06} studied the dispersive effect of the evolution group $\{e^{\pm i\Omega tD_3/|D|}\}_{t \in \mathbb{R}}$ generated by the Coriolis force, and established the Strichartz estimate for the linear solution.
As a consequence of this Strichartz estimate, {some space-time norm of} the linear solution may be small by taking the speed $|\Omega|$ of the rotation sufficiently high.
Iwabuchi and Takada \cites{Iw-Ta-13,Iw-Ta-15} and Koh, Lee, and Takada \cite{Ko-Le-Ta-14-1} improved the Strichartz estimate and it was shown by \cites{Iw-Ta-13,Ko-Le-Ta-14-1} that the global well-posedness of \eqref{eq:incomp-coriolis} for any divergence free initial data $u_0 \in \dot{H}^s(\mathbb{R}^3)$ with $1/2 \leq s < 9/10$, provided that $|\Omega|$ is sufficiently high.}
%
%
    We next recall the results of the viscous compressible fluid.
    For the non-rotational compressible Navier--Stokes system
    \begin{align}\label{eq:comp-NS}
	\begin{dcases}
		\partial_t \rho + \div ( \rho u ) = 0, 
		& 
		t > 0, x \in \mathbb{R}^3,\\
		\rho \left( \partial_t u + ( u \cdot \nabla ) u \right)
		+
		\nabla P ( \rho ) 
		={}
		\mu \Delta u + ( \mu + \mu' ) \nabla \div u,
		&
		t > 0, x \in \mathbb{R}^3,\\
		\rho(0,x) = \rho_0(x), \quad u(0,x) = u_0(x),
		&
		x \in \mathbb{R}^3.
	\end{dcases}
    \end{align}
    Starting with the pioneering work by Matsumura and Nishida \cite{Ma-Ni-79}, there are a number of papers that dealt with the global well-posedness of \eqref{eq:comp-NS}. 
    In the scaling critical framework, Danchin \cite{Da-00} proved 
    the global well-posedness  of \eqref{eq:comp-NS} in $(\dot{B}_{2,1}^{\frac{1}{2}}(\mathbb{R}^3) \cap \dot{B}_{2,1}^{\frac{3}{2}}(\mathbb{R}^3)) \times \dot{B}_{2,1}^{\frac{1}{2}}(\mathbb{R}^3)^3$.
    Later, this result was improved to the case of the general $L^p$-critical framework by \cites{Ch-Da-10,Ha-11,Ch-Mi-Zh-10-1}.
    For the related topics on strong solutions with the critical regularity, see \cites{Dan-Xu-17,Xin-Xu-21} for the large time behavior and \cites{Da-He-16,Da-02-R,Fu-24} for the low Mach number limit.
%
%
    In contrast to either the incompressible or non-rotational case, the well-posedness results on the case of $\Omega \neq 0$ in \eqref{eq:NSC-1} are few. 
    We refer to \cites{Fa-21,Fe-Ga-No-12,Fe-Ga-Ge-No-12,Fe-No-14,Bo-Fa-Pr-22} for the existence and the fast rotation limit for \eqref{eq:NSC-1}  {in an infinite slab $\mathcal D = \mathbb R^2 \times (0, 1)$} in the context of weak solutions.
    The study of the well-posedness for strong solutions to \eqref{eq:NSC-1} on $\mathbb{R}^3$ is first considered by the authors' work \cite{Fuj-Wat-25}.
    In \cite{Fuj-Wat-25}, the authors established the Strichartz estimates for linearized solutions via the analysis of dispersion due to the {mixture} of the Coriolis force and acoustic wave, which led them to prove the {\it long time} existence of the solutions to \eqref{eq:NSC-1} in the scaling critical Besov space framework;
    more precisely, for any initial disturbance {in $(\dot{B}_{2,1}^{\frac{1}{2}}(\mathbb{R}^3) \cap \dot{B}_{2,1}^{\frac{3}{2}}(\mathbb{R}^3)) \times \dot{B}_{2,1}^{\frac{1}{2}}(\mathbb{R}^3)^3$} and arbitrarily large time $0<T<\infty$, there exists a positive constant $\Omega_T=\Omega_T(\rho_0,u_0)$ such that \eqref{eq:NSC-1} possesses a unique solution on the time interval $[0,T)$, provided that $\Omega_T \leq |\Omega| \leq 1/\varepsilon$.
    Unfortunately,
    the proof of the aforementioned result fails if we take $T=\infty$.
    One of the reasons is that the zeroth order Coriolis term $\Omega ( e_3 \times u )$ interacts with the density perturbation {and prevents} us from establishing the energy estimates in the low-frequency region.
    {In addition, if one intends to focus on the eigen frequencies of the linearized system for investigating the properties of the linear part, one needs to solve
    \begin{align}
        \begin{split}
        \lambda^4
        &
        - (4\mu + \mu') |\xi|^2 \lambda^3
        + \left\{\frac{|\xi|^2}{\varepsilon^2} + \mu (5\mu + 2\mu') |\xi|^4 + \Omega^2 \right\} \lambda^2\\
        &
        - \left\{\frac{2\mu}{\varepsilon^2}|\xi|^4 + \mu^2 |\xi|^6 + \Omega^2\mu|\xi|^2 + \Omega^2(\mu + \mu') \xi_3^2 \right\} \lambda
        + \frac{\mu^2}{\varepsilon^2}|\xi|^6 + \frac{\Omega^2}{\varepsilon^2}\xi_3^2 = 0,
        \end{split}
    \end{align}
    see also \cite{Fuj-Wat-25}*{(1.10)}.
    However, it seems to be difficult to analyze the roots $\lambda(\xi)$ since the explicit formula is quite complicated. This implies that we cannot use the asymptotic expansion of $\lambda(\xi)$ as it depends on multiple quantities, that is, $\Omega$, $\varepsilon$, $|\xi|$, $\xi_3$.\footnote{{By the classical perturbation theory \cite{Ka-95}, the asymptotic expansion of the eigenvalues is obtained only for the one parameter case, in general.}}}
    On account of {these circumstances},
    there {seems} to be no results on the global well-posedness of {the initial value problem of} the viscous compressible rotating flow \eqref{eq:NSC-1} on the 3D whole space 
    even in the classical Matsumura--Nishida framework~{\cite{Ma-Ni-79}}.


The aim of the present paper is to overcome the aforementioned difficulty and construct a unique global { solution} to \eqref{eq:NSC-1} in the scaling critical Besov spaces framework.
Moreover, {we take the advantage of} the dispersive nature from the Coriolis force and the acoustic wave into our analysis and show that the global unique solution {may be constructed} for arbitrarily large initial disturbance provided that the speed $|\Omega|$ of the rotation and the {sound speed} $1/\varepsilon$ is sufficiently high in the sense $1 \ll |\Omega| \ll 1/\varepsilon$.


This paper is organized as follows.
In the next section, we provide the precise statements of the main results. 
We focus on the linear analysis in Sections \ref{sec:lin}.
In Section \ref{sec:a-priori}, we establish global a priori estimates.
Finally, we show Theorem \ref{thm:large} in Section \ref{sec:pf-large}.
In Appendix, we put a toolbox on the Chemin--Lerner spaces.


{
    \section{Main results}
    In this section, we give the precise statement of our main theorem.
    To begin with, we prepare the basic notations and definitions of function spaces that are frequently used throughout this paper.
}
    \subsection*{Notations}
    We fix some notations that are to be used throughout the whole article.
    For $\alpha \in \mathbb R$, we set $\langle \alpha \rangle:= \sqrt{1+\alpha^2}$.
    {A generic constant is denoted by $C$,} which may differ in each line.
    In, particular, $C=C(a_1,...,a_n)$ means that $C$ depends only on quantities $a_1,...,a_n$.
    For two real numbers $X$ and $Y$, we use the notation $X \sim Y$ to describe $C^{- 1} X \leq Y \leq C X$ with some positive constant $C$.
    A commutator for two operators $A$ and $B$ is denoted by $[A,B]=AB-BA$.

    Let $\mathscr{S}(\mathbb{R}^3)$ be the set of all Schwartz functions on $\mathbb{R}^3$ and
    $\mathscr{S}'(\mathbb{R}^3)$ be the set of all tempered distributions on $\mathbb{R}^3$.
    For $f \in \mathscr{S}'(\mathbb{R}^3)$, the Fourier transform and the inverse transform of $f$ are denoted by $\mathscr{F}[f]=\widehat{f}$ and $\mathscr{F}^{-1}[f]$, respectively.
    In the following, we recall the definition of the Besov spaces and the Littlewood--Paley theory.
    Let $\phi_0 \in \mathscr{S}(\mathbb{R}^3)$ satisfy
	$0 \leq \phi_0 \leq 1$, 
	$\supp \phi_0 \subset \{\xi \in \mathbb{R}^3 \ ;\ 2^{-1} \leq |\xi| \leq 2 \}$,
	and $\sum_{j \in\mathbb{Z}} \phi_j(\xi) = 1$ for all $\xi \in \mathbb{R}^3\setminus \{0 \}$,
    where $\phi_j(\xi):=\phi(2^{-j}\xi)$.
    Let $\{\Delta_j \}_{j \in \mathbb{Z}}$ be the Littlewood--Paley projection operators defined as 
       $ \Delta_j f := \mathscr{F}^{-1}\phi_j\mathscr{F}f$.
    For $1 \leq p,\sigma \leq \infty$ and $s \in \mathbb{R}$, the Besov space $\dB_{p,\sigma}^s(\mathbb{R}^3)$ is defined by
    \begin{align}
        \dB_{p,\sigma}^s(\mathbb{R}^3)
        :=
        \left\{
        f \in \mathscr{S}'(\mathbb{R}^3)/\mathscr{P}(\mathbb{R}^3)
        \ ;\ 
        \| f \|_{\dB_{p,\sigma}^s} < \infty
        \right\}
    \end{align}
    {with}
    \begin{align}
        \| f \|_{\dB_{p,\sigma}^s}
        :=
        \left\| 
        \{2^{sj} \| \Delta_j f \|_{L^p} \}_{j \in \mathbb{Z}} 
        \right\|_{\ell^{\sigma}(\mathbb{Z})}.
    \end{align}
    Here $\mathscr{P}(\mathbb{R}^3)$ denotes the set of all polynomials on $\mathbb{R}^3$.
    It is well-known that if $(s,\sigma) \in (-\infty, 3 \slash p) \times [1,\infty]$ or $(s,\sigma)=(3/p,1)$, 
    then $\dB_{p,\sigma}^s(\mathbb{R}^3)$ is {identified\footnote{See \cite{Sa-18}*{(2.124) in Theorem 2.31} for the precise identification map.}} by
    \begin{align}
        \left\{
        f \in \mathscr{S}'(\mathbb{R}^3)
        \ ;\ 
        f = \sum_{j \in \mathbb{Z}} \Delta_j f
        {\rm \ in\ }
        \mathscr{S}'(\mathbb{R}^3)
        {\rm \ and\ }
        \| f \|_{\dB_{p,\sigma}^s} < \infty
        \right\}.
    \end{align}
    The truncated Besov semi-norms are defined by
    \begin{align}
        \| f \|_{\dB_{p,\sigma}^s}^{h;\beta}
        :={}&
        \left\| 
        \{2^{sj} \| \Delta_j f \|_{L^p} \}_{j} 
        \right\|_{\ell^{\sigma}(\{j \in \mathbb{Z}\ ; \ \beta < 2^j\})},\\
        \| f \|_{\dB_{p,\sigma}^s}^{m;\alpha,\beta}
        :={}&
        \left\| 
        \{2^{sj} \| \Delta_j f \|_{L^p} \}_{j} 
        \right\|_{\ell^{\sigma}( \{j \in \mathbb{Z}\ ;\ \alpha < 2^j \leq \beta \} )},\\
        \| f \|_{\dB_{p,\sigma}^s}^{\ell;\alpha}
        :={}&
        \left\| 
        \{2^{sj} \| \Delta_j f \|_{L^p} \}_{j} 
        \right\|_{\ell^{\sigma}(\{j \in \mathbb{Z} \ ;\ 2^j \leq \alpha \})}
    \end{align}
    for $0 \leq \alpha < \beta \leq \infty$.
    To control functions {with respect to} the space-time variable, we use the Chemin--Lerner spaces {defined as follows}.
    For $1 \leqslant p,\sigma,r \leqslant \infty$, $s \in \mathbb{R}$ and an interval $I \subset \mathbb R$, define\footnote{To our best knowledge, there seems {to be} no literature that mention on a similar identification for the Chemin--Lerner spaces, but we may easily obtain it by mimicking the proof of [33, Theorem 2.31]. However, we do not need this identification. Since we only use the norm of the Chemin-Lerner spaces for obtaining the global a priori estimates and the $\widetilde{C}$-spaces are not needed for constructing a local unique strong solution.}
    \begin{align}
        \widetilde{L^r}(I ; \dB_{p,\sigma}^s(\mathbb{R}^3))
        :={}&
        \left\{
        F : I \to \mathscr{S}'(\mathbb R^3)\ ;\ 
        \| F \|_{\widetilde{L^r}(I ; \dB_{p,\sigma}^s)} < \infty 
        \right\},\\
        \| F \|_{\widetilde{L^r}(I ; \dB_{p,\sigma}^s)}
        :={}&
        \left\| 
        \{2^{sj} \| \Delta_j F \|_{L^r( I ; L^p)} \}_{j \in \mathbb{Z}} 
        \right\|_{\ell^{\sigma}(\mathbb{Z})}.
    \end{align}
    We also use the notation
	\begin{equation}
		\widetilde C (I; \dB^s_{p,\sigma} (\mathbb{R}^3)) := C (I; \dB^s_{p,\sigma} (\mathbb{R}^3)) \cap \widetilde {L^\infty} (I; \dB^s_{p,\sigma} (\mathbb{R}^3))
	\end{equation}
    {as well as}
    \begin{align}
        \widetilde C (I; \dB^{s_1}_{p_1,\sigma_1} (\mathbb{R}^3) \cap \dB^{s_2}_{p_2,\sigma_2} (\mathbb{R}^3))
        :=
        \widetilde C (I; \dB^{s_1}_{p_1,\sigma_1} (\mathbb{R}^3))
        \cap
        \widetilde C (I; \dB^{s_2}_{p_2,\sigma_2} (\mathbb{R}^3)).
    \end{align}
    Define the truncated type {semi-norms} as 
    \begin{align}
        \| F \|_{\widetilde{L^r}(I ; \dB_{p,\sigma}^s)}^{h;\beta}
        :={}&
        \left\| 
        \{2^{sj} \| \Delta_j F \|_{L^r( I ; L^p)} \}_{j} 
        \right\|_{\ell^{\sigma}(\{j \in \mathbb{Z}\ ; \ \beta < 2^j\})},\\
        \| F \|_{\widetilde{L^r}(I ; \dB_{p,\sigma}^s)}^{m;\alpha,\beta}
        :={}&
        \left\| 
        \{2^{sj} \| \Delta_j F \|_{L^r( I ; L^p)} \}_{j} 
        \right\|_{\ell^{\sigma}( \{j \in \mathbb{Z}\ ;\ \alpha < 2^j \leq \beta \} )},\\
        \| F \|_{\widetilde{L^r}(I ; \dB_{p,\sigma}^s)}^{\ell;\alpha}
        :={}&
        \left\| 
        \{2^{sj} \| \Delta_j F \|_{L^r( I ; L^p)} \}_{j} 
        \right\|_{\ell^{\sigma}(\{j \in \mathbb{Z} \ ;\ 2^j \leq \alpha \})}.
    \end{align}
	We also define the corresponding norms $\| \,\cdot \, \|_{L^r(I ; \dB_{p,\sigma}^s)}^{h;\beta}$, $\| \,\cdot \, \|_{L^r (I ; \dB_{p,\sigma}^s)}^{m;\alpha,\beta}$, and $\| \,\cdot \, \|_{L^r(I ; \dB_{p,\sigma}^s)}^{\ell;\alpha}$, similarly.
    We summarize the basic properties on the Chemin--Lerner spaces in Appendix \ref{sec:a}.
\subsection*{{ Main Theorem}}
    To state our main result precisely, we reformulate the problem.
    For simplicity of notation, we perform a suitable rescaling transform in \eqref{eq:NSC-1} so as to reduce to the case 
    \begin{align}
        \nu = \rho_\infty = c_\infty = 1.
    \end{align} 
    We refer to \cite{Fuj-Wat-25} for the precise explanation of this scaling argument.

    Let $a$ be  {the density perturbation from the constant $\rho_{\infty}=1$} 
    defined by
    \begin{align}
        a(t,x) : = \frac{\rho(t,x) - 1}{\varepsilon} \quad \text{with} \quad
		 a(0,x) = a_0(x) : = \frac{\rho_0(x) - 1}{\varepsilon}.
	\end{align}
    Then,  {the system \eqref{eq:NSC-1} is rewritten as}
    \begin{align}\label{eq:NSC-2}
        \begin{dcases}
            \partial_t a + \dfrac{1}{\varepsilon} \div u = - \div ( a u ), 
            & 
            t > 0, x \in \mathbb{R}^3,\\
            \partial_t u 
            -
            \mathcal{L} u
            + 
            \Omega ( e_3 \times u ) 
            +
            \dfrac{1}{\varepsilon} \nabla a 
            =
            - {\mathcal{N}_{\varepsilon}}[a,u],
            &
            t > 0, x \in \mathbb{R}^3,\\
            a(0,x) = a_0(x), \quad u(0,x) = u_0(x),
            &
            x \in \mathbb{R}^3,
        \end{dcases}
    \end{align}
    where we have {set} $\mathcal{L}u : = \mu \Delta u + ( \mu + \mu') \nabla \div u$ with $\mu>0$ and $2\mu+\mu'=1$, and
    \begin{gather}
        {\mathcal{N}_{\varepsilon}}[a,u]
        :=
        ( u \cdot \nabla ) u
        +
        J( \varepsilon a ) \mathcal{L}u
        +
        \dfrac{1}{\varepsilon}
        K( \varepsilon a )
        \nabla a,\\
        J(a)
        :=
        \frac{a}{1+a},\qquad
        K(a)
        :=
        \frac{P'(1+a)}{1+a}-1.
    \end{gather}
    Now, we state {the} main result of this paper as follows.
    \begin{thm}\label{thm:large}
	Let 
	$a_0 \in \dB_{2,\infty}^{-\frac{3}{2}} ( \mathbb{R}^3 ) \cap \dB_{2,1}^{\frac{3}{2}} ( \mathbb{R}^3 )$
	and
	$u_0 \in ( \dB_{2,\infty}^{-\frac{3}{2}} ( \mathbb{R}^3 ) \cap \dB_{2,1}^{\frac{1}{2}} ( \mathbb{R}^3 ) )^3$.
	Then, there {exist} positive constants $\Omega_0=\Omega_0(\mu,P,a_0,u_0)$ and $c_0=c_0(\mu,P,a_0,u_0)$ such that if {$\Omega \in \mathbb{R}$} and $\varepsilon > 0$ {satisfy}
	\begin{align}\label{O-e}
		\Omega_0 \leq |\Omega| \leq c_0 \frac{1}{\varepsilon},
	\end{align}
	then
	the {system} \eqref{eq:NSC-2} possesses a unique solution $(a,u)$ in the class
	\begin{align}\label{sol-class}
		&
		a \in \tC ( [0,\infty) ; \dB_{2,\infty}^{-\frac{3}{2}} ( \mathbb{R}^3 ) \cap \dB_{2,1}^{\frac{3}{2}} ( \mathbb{R}^3 ) ),\\
		&
		u \in \tC ( [0,\infty) ; \dB_{2,\infty}^{-\frac{3}{2}} ( \mathbb{R}^3 ) \cap \dB_{2,1}^{\frac{1}{2}} ( \mathbb{R}^3 ) )^3
		\cap 
		L^1_{\rm loc}( [0,\infty) ; \dB_{2,1}^{\frac{5}{2}} (\mathbb{R}^3) )^3
	\end{align}
	with {$a(t,x) > -1/\varepsilon$ for all $(t,x) \in [0,\infty) \times \mathbb{R}^3$}.
    \end{thm}
    {
    \begin{rem}
        Let us make some comments on Theorem \ref{thm:large}.
        \begin{enumerate}
            \item 
            Let us mention the assumption \eqref{O-e}.
            Analogously to the incompressible case \cite{Iw-Ta-13}, the condition $\lvert\Omega\rvert \gg 1$ serves to {make} a certain Strichartz norm of the solution {small}.
            We exploit this fact in the \textit{medium}-frequency estimates; see Lemmas \ref{lemm:middle-str-1} and \ref{lemm:middle-str-2}.
            For the condition $|\Omega| \ll 1/\varepsilon$,
            the \textit{low}-frequency estimates established in Lemma \ref{lemm:low-ene} require $|\Omega|\varepsilon \ll 1$ to close the global a priori estimates; see also \eqref{choice}.
            We also use $1/\varepsilon \gg 1$ for {making the \textit{high}-frequency part of the solution small}.
            Summarizing, the single condition \eqref{O-e} suffices to close all global a priori estimates.
            \item             
            As mentioned above, the assumption $|\Omega| \ll 1/\varepsilon$ is technical but reasonable. 
            Indeed, considering the low Mach number limit $\varepsilon \to +0$ with a fixed $\Omega$, we see that the assumption \eqref{O-e} reduces to $|\Omega| \geq \Omega_0$, which is the same condition as for the incompressible viscous rotating fluid; see \cite{Iw-Ta-13}.
            More precisely, it is expected that the large solution constructed in the present paper converges to the large incompressible rotating Navier--Stokes flow given by \cite{Iw-Ta-13} in some space-time norm as $\varepsilon \to +0$, provided that $\Omega \geq \Omega_0$.
            \item 
            In the rotational setting ($\Omega \ne 0$), one is forced to work in the Chemin--Lerner spaces framework for the solution class \eqref{sol-class}, which is a typical difference from the non-rotational setting ($\Omega = 0$). This necessity arises from estimates \eqref{lemm:low-ene-2}, \eqref{lemm:low-ene-3}, and \eqref{lemm:low-ene-4} in Lemma \ref{lemm:low-ene}. This is because the interpolation index on the homogeneous Besov norms appearing on the left-hand sides is infinity. Consequently, one has to investigate in the framework of the Cheimin--Lerner spaces, and terms such as $\n{a}_{\widetilde{L^{\infty}}(0,t;\dB_{q,1}^{\frac{3}{q}-1})}$  appear on the right-hand sides of the estimates.
            \item 
            The lower bound $a>-1/\varepsilon$ implies the positivity of the density, that is, $\rho =1+\varepsilon a > 0$ in the original system \eqref{eq:NSC-1}. 
        \end{enumerate}
        
    \end{rem}
    }
    \subsection*{{ Sketch of the Key Ideas}}
    Let us make comments on the principal difficulties in proving Theorem \ref{thm:large}.
    In light of the critical Besov framework {approach invented} by Danchin \cite{Da-00}, we meet two difficulties in our analysis.
    The first one appears from the linear analysis.
    In fact, as mentioned before, we may not directly apply the well-known method (cf. \cite{Da-00}) for the energy estimates of the linearized system in the low-frequency part since the lowest-order term $\langle \Omega (e_3 \times \widehat{u}), \varepsilon i\xi \widehat{a} \rangle_{\mathbb{C}^3}$ is hard to handle when we estimate the time evolution of $\langle \widehat{u}, \varepsilon i\xi \widehat{a} \rangle_{\mathbb{C}^3}$, which plays a key role in the energy argument for the low-frequency part.
    This forced the authors in the previous paper \cite{Fuj-Wat-25} to restrict the life span of the solutions to finite time intervals. 
    In this paper, we examine this harmful term more precisely and find that the linearized solution behaves as the $4$th order dissipative semigroup $\{e^{-t\Delta^2}\}_{t>0}$ in the low-frequency part; see \eqref{ene-low} below.
    It should be emphasized that this is completely different from the non-rotational case \eqref{eq:comp-NS}, where the linearized solution behaves like the heat kernel $\{ e^{t\Delta} \}_{t>0}$ in the low-frequency region.
    Namely, the situation becomes much worse in our case.
    Indeed, the maximal $L^1$-regularity estimate for $\{e^{-t\Delta^2}\}_{t>0}$ is given as 
    \begin{align}
        \n{e^{-t\Delta^2} f}_{\widetilde{L^1}(0,\infty; \dB_{2,\sigma}^{\frac{5}{2}})}
        \leq
        C
        \n{f}_{\dB_{2,\sigma}^{-\frac{3}{2}}}
    \end{align}
    for $1 \leq \sigma \leq \infty$.
    Then, we are faced with the negative regularity {space} $\dB_{2,\sigma}^{-\frac{3}{2}}(\mathbb{R}^3)$, which {creates} the second difficulty when we estimate the nonlinear terms.
    As was mentioned in \cite{Xin-Xu-21}*{Section 5}, the derivative of order $s=-3/2$ is just barely controllable in Besov spaces based on $L^2(\mathbb{R}^3)$ by choosing $\sigma=\infty$.
    However, if we utilize the dispersive estimates based on $L^q(\mathbb{R}^3)$-norms, even lower negative regularity space appears so that more ingenuity is required.
    To circumvent this difficulty, we rely on the momentum formulation. 
    It is well-known that considering the momentum $m:= \rho u= (1+\varepsilon a)u$, then all nonlinear terms in the equations for $m$ have a divergence form, which has a better regularity structure and rescues the problem of nonlinear estimates from us.
    {Notice that the momentum formulation is applied only in the low-frequency part, since the difficulty in estimating the lower negative regularity space that we are facing occurs in the low-frequency region;
    see also Lemma \ref{lemm:low-ene} and Remark \ref{rem:low-ene} below for further details.
    If we deal with the high-frequency part of $(a, m)$, we meet a derivative loss of the density perturbation $a$ in the higher-order term $\varepsilon \mathcal{L} (au)$ in \eqref{eq:m}. Hence, we use the momentum form only in the low-frequency part.
    }
\section{Linear analysis}\label{sec:lin}
    In this section, we prepare several estimates of the solutions to the following linearized system:
    \begin{align}\label{eq:lin}
        \begin{dcases}
            \partial_t a + \dfrac{1}{\varepsilon} \div u = f, 
            & 
            t > 0, x \in \mathbb{R}^3,\\
            \partial_t u 
            -
            \mathcal{L} u
            + 
            \Omega ( e_3 \times u ) 
            +
            \dfrac{1}{\varepsilon} \nabla a 
            =
            g,
            &
            t > 0, x \in \mathbb{R}^3,\\
            a(0,x) = a_0(x), \quad u(0,x) = u_0(x),
            &
            x \in \mathbb{R}^3,
        \end{dcases}
    \end{align}
    where the initial data $(a_0,u_0)$ and the external forces $(f,g)$ are given functions belonging to some suitable function spaces.
    First, we recall the following elementary energy estimate.
    \begin{lemm}[\cite{Fuj-Wat-25}*{{Lemma 2.1}}]\label{lemm:simple-ene}
        There exists a positive constant $C=C(\mu)$ such that
        for each $\Omega \in \mathbb{R}$ and $\varepsilon>0$,
        the solutions $(a,u)$ to \eqref{eq:lin} satisfies
        \begin{align}\label{simple-ene}
        \begin{split}
        \| \Delta_j(a,u) \|_{L^{\infty}(0,t;L^2)}
        +
        \| \Delta_j \nabla u \|_{L^2(0,t;L^2)}
        \leq
        C
        \| \Delta_j(a_0,u_0) \|_{L^2}
        +
        C
        \| \Delta_j(f,g) \|_{L^1(0,t;L^2)}
        \end{split}
        \end{align}
        for all $t >0$ and $j \in \mathbb{Z}$, provided that the right{-}hand side is finite.
    \end{lemm}
     {
    The {main objective} in this section is to establish the following lemma.}
    \begin{lemm}\label{lemm:ene-low}
	Let $1 \leq r, \sigma \leq \infty$.
	Let $\alpha>0$, $\beta\geqslant 1$, $\Omega \in \mathbb{R}$, and $0 < \varepsilon \leq 1$ satisfy
	$|\Omega|\varepsilon \leq \alpha < \beta/\varepsilon$.
	Then, there exists a positive constant $C=C(\mu,r)$ such that
	\begin{align}
		&
		\begin{aligned}
			\| (a,u) \|_{\widetilde{L^r}(0,t ; \dB_{2,\sigma}^{s + \frac{2}{r}} )}^{m;\alpha,\frac{\beta}{\varepsilon}}
			\leq{}&
			C 
			\beta^{\frac{2}{r}}
			\| (a_0, u_0) \|_{\dB_{2,\sigma}^s}^{m;\alpha,\frac{\beta}{\varepsilon}}
			+
			C 
			\beta^{\frac{2}{r}}
			\| (f, g) \|_{\widetilde{L^1}( 0,t ; \dB_{2,\sigma}^s)}^{m;\alpha,\frac{\beta}{\varepsilon}},
		\end{aligned}\\
		&
		\begin{aligned}\label{est:4max-reg}
			\| (a,u) \|_{\widetilde{L^r}(0,t ; \dB_{2,\sigma}^{s + \frac{4}{r}} )}^{\ell;\alpha}
			\leq{}&
			C 
			\beta^{\frac{2}{r}}
			\left( \alpha^2 + \Omega^2 \varepsilon^2 \right)^{\frac{1}{r}}
			\| (a_0, u_0) \|_{\dB_{2,\sigma}^s}^{\ell;\alpha}\\
			&
			+
			C 
			\beta^{\frac{2}{r}}
			\left( \alpha^2 + \Omega^2 \varepsilon^2 \right)^{\frac{1}{r}}
			\| (f, g) \|_{\widetilde{L^1}( 0,t ; \dB_{2,\sigma}^s)}^{\ell;\alpha}
		\end{aligned}
	\end{align}
	for all $t>0$, provided that the right-hand side is finite.
\end{lemm}
\begin{rem}
    As {will be} seen from \eqref{ene-low} below, the estimate for low{-}frequency part corresponds to the maximal {$L^r$-}regularity of the $4$th order dissipative semigroup $\{e^{-t\Delta^2}\}_{t>0}$.
    Then, in particular, the maximal {$L^1$-}regularity lowers the differential exponent by four orders, which yields difficulties in nonlinear estimates; see Lemma \ref{lemm:low-ene} and Remark \ref{rem:low-ene} below.
\end{rem}
\begin{proof}
	We first intend to prove that there exists a positive constant $C=C(\mu,r)$ such that
	 {
	\begin{align}\label{ene-low}
		\begin{split}
			| ( \ha, \hu )(t,\xi) |
			\leq {} 
			C 
			e^{- \frac{\kappa (|\xi|, \Omega\varepsilon)}{48\beta^2} t}
			| ( \ha_0,\hu_0 )(\xi) | 
			+
			C 
			\int_0^t
			e^{- \frac{\kappa (|\xi|, \Omega\varepsilon)}{48\beta^2} (t-\tau)}
			| ( \hf, \hg )(\tau,\xi)  |
			d\tau
		\end{split}
	\end{align}
	with}
     {
    \begin{align}
         \kappa (|\xi|, \Omega\varepsilon) := \frac{|\xi|^4 }{\Omega^2 \varepsilon^2 + |\xi|^2}
    \end{align}
    }
	for all $\xi \in \mathbb{R}^3$ with $|\xi| \leq  {2}\beta / \varepsilon$ and $t>0$, provided that the right-hand side is finite.
	Applying the Fourier transform to \eqref{eq:lin}, we see that
	\begin{align}\label{eq:lin-F}
		\begin{dcases}
			\partial_t \ha 
			+ 
			\frac{i}{\varepsilon}\xi \cdot \hu 
			= 
			\hf,\\
			\partial_t \hu + \mu |\xi|^2 \hu + (\mu+\mu') \xi (\xi \cdot \hu) + \Omega  ( e_3 \times \hu )
			+
			\frac{i}{\varepsilon}\xi \ha
			=
			\hg.
		\end{dcases}
	\end{align}
	In the {sequel} of this proof, we  {first} assume $|\xi| \leq 2\beta/\varepsilon$.
	We multiply the first equation of \eqref{eq:lin-F} by $\overline{\ha}$ and take the $\mathbb{C}^3$-inner product of the second equation with $\hu$.
	Then, summing them up and taking the real part, we have
	\begin{align}\label{pf:ene-low-1}
		\frac{1}{2}
		\partial_t
		| (\ha, \hu) |^2
		+
		\underline{\mu}
		|\xi|^2 |\hu|^2
		\leq
		| (\hf, \hg) |
		| (\ha, \hu) |,
	\end{align}
	where $\underline{\mu} := \min \{ \mu,1 \}$.
	Multiplying the first equation of \eqref{eq:lin-F} by $i \varepsilon \xi$ and taking the $\mathbb{C}^3$-inner product of it with $\hu$, we obtain 
	\begin{align}\label{pf:ene-low-2}
		\langle i \varepsilon \xi \partial_t \ha,  \hu \rangle_{\mathbb{C}^3}
		-| \xi \cdot  {\hu} |^2
		= 
		\langle i \varepsilon \xi \hf, \hu \rangle_{\mathbb{C}^3}.
	\end{align}
	From the $\mathbb{C}^3$-inner product of $i \varepsilon \xi \ha$ with the second equation of \eqref{eq:lin-F},
	it follows that
	\begin{align}\label{pf:ene-low-3}
		\langle i \varepsilon \xi \ha, \partial_t \hu \rangle_{\mathbb{C}^3}
		+
		|\xi|^2 
		\langle i \varepsilon \xi \ha , \hu \rangle_{\mathbb{C}^3}
		+
		\Omega \varepsilon \langle i \xi  \ha, e_3 \times \hu \rangle
		+
		|\xi|^2 |\ha|^2
		= 
		\langle i \varepsilon \xi \ha, \hg \rangle_{\mathbb{C}^3}.
	\end{align}
	By \eqref{pf:ene-low-2} and \eqref{pf:ene-low-3}, we see that
	\begin{equation}\label{pf:ene-low-4}
		\begin{split}& \partial_t
		\mathrm{Re}\,
		\langle i \varepsilon \xi \ha , \hu         \rangle_{\mathbb{C}^3}
		+
		|\xi|^2 |\ha|^2 \\
		& \leq{}
		\varepsilon 
		|\xi| 
		| (\hf,\hg) |
		| (\ha,\hu) |
		+
		|\Omega|
		\varepsilon
		| \xi |
		| \ha | | \hu | 
		+
		|\xi|^2 |\hu|^ {2}
		+
		\varepsilon |\xi|^3 | \ha | | \hu |\\
		& \leq{}
		2\beta
		| (\hf,\hg) |
		| (\ha,\hu) |
		+
		|\Omega|
		\varepsilon
		| \xi |
		| \ha | | \hu |
		+
		|\xi|^2 |\hu|^2
		+
		2\beta |\xi|^2 | \ha | | \hu |
		\end{split}
	\end{equation}
	 {by $|\xi| \leq 2 \beta \slash \varepsilon$.}
	Let $\delta = \underline{\mu}/(16\beta^2)$
	and define
	\begin{align}
		\hV^2
		=
		\hV^2(t,\xi)
		:=
		( \Omega^2 \varepsilon^2 + |\xi|^2 )
		| (\ha,\hu) |^2
		+
		2\delta |\xi|^2
		\mathrm{Re}\,
		\langle i \varepsilon \xi \ha , \hu         \rangle_{\mathbb{C}^3}.
	\end{align}
	Here, since
	\begin{align}\label{pf:ene-low-5}
		\left|
		\hV^2
		-
		( \Omega^2 \varepsilon^2 + |\xi|^2 )
		| (\ha,\hu) |^2
		\right|
		\leq{}
		2\delta \varepsilon |\xi|^3 |\ha| |\hu|
		\leq{}
		2\delta \beta |\xi|^2 | (\ha,\hu) |^2 
		\leq{}
		\frac{1}{2}|\xi|^2 | (\ha,\hu) |^2,
	\end{align}
	we see that 
	\begin{align}\label{pf:ene-low-6}
		\frac{1}{2}
		( \Omega^2 \varepsilon^2 + |\xi|^2 )
		| (\ha,\hu) |^2
		\leq
		\hV^2
		\leq
		\frac{3}{2}
		( \Omega^2 \varepsilon^2 + |\xi|^2 )
		| (\ha,\hu) |^2.
	\end{align}
	By \eqref{pf:ene-low-1} and \eqref{pf:ene-low-4},
	it holds
	\begin{align}\label{pf:ene-low-7}
		&
		\frac{1}{2} \partial_t \hV^2
		+
		\underline{\mu} \Omega^2 \varepsilon^2
		|\xi|^2 |\hu|^2
		+
		|\xi|^4 |\hu|^2
		+
		\delta |\xi|^4 |\ha|^2\\
		&\quad
		\leq
		( \Omega^2 \varepsilon^2 +  {\langle 2 \delta \beta \rangle} |\xi|^2 )
		| (\hf, \hg) | | (\ha, \hu) |
		+
		\delta
		|\Omega|
		\varepsilon
		|\xi|^3
		| \ha | | \hu |
		+
		\delta
		|\xi|^4 |\hu|^2
		+
		 {2}\delta \beta |\xi|^4 | \ha | | \hu |.
	\end{align}
	Hence, it follows from
	\begin{align}\label{pf:ene-low-8}
		&
		\delta
		|\Omega|
		\varepsilon
		|\xi|^3
		| \ha | | \hu |
		\leq
		 {4} \delta \Omega^2 \varepsilon^2|\xi|^2 | \hu |^2
		+
		\frac{\delta}{4}
		| \xi |^4 | \ha |^2
		\leq
		\underline{\mu} \Omega^2 \varepsilon^2|\xi|^2 | \hu |^2
		+
		\frac{\delta}{4}
		| \xi |^4 | \ha |^2,
		\\
		&\delta
		|\xi|^4 |\hu|^2
		\leq
		\frac{1}{4}\underline{\mu}
		|\xi|^4 |\hu|^2,\\
		&
		2\delta \beta |\xi|^4 | \ha | | \hu |
		\leq
		 {4} \delta \beta^2 | \xi |^4 | \hu |^2
		+
		\frac{\delta}{4}
		| \xi |^4 | \ha |^2
		\leq
		\frac{1}{4}\underline{\mu} | \xi |^4 | \hu |^2
		+
		\frac{\delta}{4}
		| \xi |^4 | \ha |^2
	\end{align}
	that 
	\begin{align}\label{pf:ene-low-9}
		&
		\frac{1}{2} \partial_t \hV^2
		+
		\frac{1}{2}\underline{\mu}|\xi|^4 |\hu|^2
		+
		\frac{\delta}{2} |\xi|^4 |\ha|^2
		\leq
		 {\langle 2 \delta \beta \rangle}
		(\Omega^2 \varepsilon^2 + |\xi|^2)
		| (\hf, \hg) | | (\ha, \hu) |,
	\end{align}
	which implies 
	\begin{align}\label{dff-V}
		\frac{1}{2}
		\partial_t 
		\hV^2
		+
		\frac{\delta}{3}
		 {\kappa (|\xi|, \Omega\varepsilon) }
		\hV^2
		\leq
		4
		\sqrt{\Omega^2 \varepsilon^2 + |\xi|^2}
		| (\hf, \hg) | \hV.
	\end{align}
	Thus, we obtain
	\begin{align}\label{pf:ene-low-10}
		\begin{split}
			\hV(t,\xi)
			\leq{}
			&
			e^{-\frac{\delta}{3}  {\kappa (|\xi|, \Omega\varepsilon) } t}
			\hV(0,\xi)\\
			&
			+
			4
			\sqrt{\Omega^2 \varepsilon^2 + |\xi|^2}
			\int_0^t
			e^{-\frac{\delta}{3}  {\kappa (|\xi|, \Omega\varepsilon)} (t-\tau)}
			| (\hf,\hg)(\tau,\xi) |
			d\tau.
		\end{split}
	\end{align}
	Combining \eqref{pf:ene-low-6} and \eqref{pf:ene-low-10}, 
	we obtain \eqref{ene-low}. \par
	Let $j \in \mathbb{Z}$ satisfy $2^j \leq \beta/\varepsilon$.
	Multiplying \eqref{ene-low} by $\widehat{\phi_j}(\xi)$ and using $\supp \widehat{\phi_j} \subset \{ 2^{j-1} \leq |\xi| \leq 2^{j+1} \}$, we have
	\begin{align}\label{pf:lin-FB-1}
		\begin{split}
			\widehat{\phi_j}(\xi)
			| ( \ha, \hu )(t,\xi) |
			\leq {} &
			C 
			e^{- \frac{c}{\beta^2}  {\kappa (2^j, \Omega\varepsilon)} t}
			\widehat{\phi_j}(\xi)
			| ( \ha_0,\hu_0 ) | \\
			&
			+
			C 
			\int_0^t
			e^{-  \frac{c}{\beta^2}  {\kappa (2^j, \Omega\varepsilon)} (t-\tau)}
			\widehat{\phi_j}(\xi)
			| ( \hf(\tau,\xi), \hg(\tau,\xi) ) |
			d\tau,
		\end{split}
	\end{align}
	which implies
	\begin{align}\label{pf:lin-FB-5}
		\begin{split}
			\| \Delta_j (a,u) \|_{L^2}
			\leq{} &
			C
			e^{-  \frac{c}{\beta^2}  {\kappa (2^j, \Omega\varepsilon)} t}
			\| \Delta_j (a_0,u_0) \|_{L^2}\\
			&
			+
			C
			\int_0^t
			e^{-  \frac{c}{\beta^2}  {\kappa (2^j, \Omega\varepsilon)} (t-\tau)}
			\| \Delta_j(f,g)(\tau) \|_{L^2}
			d\tau.
		\end{split}
	\end{align}
	Taking $L^r(0,t)$-norm of \eqref{pf:lin-FB-5} and using the Hausdorff--Young inequality, we obtain 
	\begin{align}\label{pf:lin-FB-2}
		\begin{split}
			\| \Delta_j (a,u) \|_{L^r( 0,t ; L^2)}
			\leq{}&
			C
			\beta^{\frac{2}{r}}
			\left( 2^{-2j} +  \Omega^2 \varepsilon^2 2^{-4j} \right)^{\frac{1}{r}}
			\| \Delta_j (a_0,u_0) \|_{L^2}\\
			&+
			C
			\beta^{\frac{2}{r}}
			\left( 2^{-2j} +  \Omega^2 \varepsilon^2 2^{-4j} \right)^{\frac{1}{r}}
			\| \Delta_j(f,g) \|_{L^1( 0,t ; L^2)}.
		\end{split}
	\end{align}
	From this, it follows
	\begin{align}\label{pf:lin-FB-3}
		\begin{split}
			2^{(s+\frac{2}{r})j}\| \Delta_j (a,u) \|_{L^r( 0,t ; L^2)}
			\leq{}&
			C
			\beta^{\frac{2}{r}}
			\left( 1 +  \Omega^2 \varepsilon^2 \alpha^{-2} \right)^{\frac{1}{r}}
			2^{sj}\| \Delta_j (a_0,u_0) \|_{L^2}\\
			&+
			C
			\beta^{\frac{2}{r}}
			\left( 1 +  \Omega^2 \varepsilon^2 \alpha^{-2} \right)^{\frac{1}{r}}
			2^{sj}\| \Delta_j(f,g) \|_{L^1( 0,t ; L^2)}\\
			\leq{}&
			C
			\beta^{\frac{2}{r}}
			2^{sj}\| \Delta_j (a_0,u_0) \|_{L^2}
			+
			C
			\beta^{\frac{2}{r}}
			2^{sj}\| \Delta_j(f,g) \|_{L^1( 0,t ; L^2)}
		\end{split}
	\end{align}
	for all $j$ with $ {|\Omega| \varepsilon \leq\;} \alpha < 2^j \leq \beta/\varepsilon$ and
	\begin{align}\label{pf:lin-FB-4}
		\begin{split}
			2^{(s+\frac{4}{r})j}\| \Delta_j (a,u) \|_{L^r( 0,t ; L^2)}
			\leq{}&
			C
			\beta^{\frac{2}{r}}
			\left( \alpha^2 +  \Omega^2 \varepsilon^2  \right)^{\frac{1}{r}}
			2^{sj}\| \Delta_j (a_0,u_0) \|_{L^2}\\
			&+
			C
			\beta^{\frac{2}{r}}
			\left( \alpha^2 +  \Omega^2 \varepsilon^2 \right)^{\frac{1}{r}}
			2^{sj}\| \Delta_j(f,g) \|_{L^1( 0,t ; L^2)}
		\end{split}
	\end{align}
	for all $j$ with $2^j \leq \alpha$.
	Hence, taking $\ell^{\sigma}$-norm of \eqref{pf:lin-FB-3} and \eqref{pf:lin-FB-4} with respect to $j$ that runs in suitable ranges, we complete the proof.
\end{proof}	
    For the high{-}frequency part, we recall the following estimate.
    \begin{lemm}[\cite{Fuj-Wat-25}*{Lemma 2.5}]\label{lemm:high-ene}
        There exist positive constants $\beta_0=\beta_0(\mu)$ and $C=C(\mu)$ such that
        for each $\Omega \in \mathbb{R}$ and $\varepsilon > 0$, the solution to {\eqref{eq:lin}} satisfies
        \begin{align}
            &
            \varepsilon 
            \| a \|_{\widetilde{L^{\infty}}( 0,t ; \dB_{p,1}^{s+1} )}^{h;\frac{\beta}{\varepsilon}}
            +
            \frac{1}{\varepsilon}
            \| a \|_{L^1( 0,t ; \dB_{p,1}^{s+1} )}^{h;\frac{\beta}{\varepsilon}}
            +
            \| u \|_{\widetilde{L^{\infty}}( 0,t ; \dB_{p,1}^{s} ) \cap L^1( 0,t ; \dB_{p,1}^{s+2} ) }^{h;\frac{\beta}{\varepsilon}}\\
            &\quad 
            \leqslant
            C
            \| (\varepsilon a_0, u_0) \|_{\dB_{p,1}^{s+1} \times \dB_{p,1}^{s}}^{h;\frac{\beta}{\varepsilon}}
            +
            \frac{C\varepsilon}{p}
            \left\|
            \| \div u \|_{L^{\infty}}
            \| a \|_{\dB_{p,1}^{s+1}}
            \right\|_{L^1(0,t)}\\
            &\qquad 
            +
            C
            \sum_{j \in \mathbb{Z}}
            2^{(s+1)j}
            \| \Delta_j f + u\cdot \nabla \Delta_j a \|_{{L^1(0,t;L^p)}}
            +
            C
            \| (f,g) \|_{L^1( 0,t ; \dB_{p,1}^{s})}^{h;\frac{\beta}{\varepsilon}}
        \end{align}
        for all $s \in \mathbb{R}$, $\beta \geqslant \beta_0 \sqrt{\langle \Omega \varepsilon^2 \rangle / 2 }$, $1 \leqslant p < \infty$, and $t>0$, provided that the right{-}hand side is finite.
    \end{lemm}
    Finally, we recall the dispersive estimate of the linear solutions which is a key ingredient for showing that we {may} take arbitrarily large initial data for {constructing} the global strong solution to \eqref{eq:NSC-2}.
    \begin{lemm}[\cite{Fuj-Wat-25}*{Proposition 3.1}]
    \label{lemm:lin-str-visc}
        Let $2 \leq q,r \leq \infty$ satisfy
        \begin{align}
            \frac{1}{q} + \frac{1}{r} \leq \frac{1}{2}, \qquad
            ( q,r ) \neq ( \infty,2 ).
        \end{align}
        Then, for $\beta>0$, there exists a positive constant $C=C(\beta,q,r)$ such that
        for any $\varepsilon, t>0$  and $\Omega \in \mathbb{R} \setminus \{ 0\}$ with $|\Omega| \varepsilon < \beta_0/\varepsilon$,     
        the solution $(a,u)$ of \eqref{eq:lin} satisfies 
        \begin{align}
			\| \Delta_j ( a,u ) \|_{L^r( 0,t ; L^q )}
			\leq{}
			C
			|\Omega|^{-\frac{1}{r}}
			2^{3(\frac{1}{2} - \frac{1}{q})j}
			\left(
			\| \Delta_j(a_0, u_0) \|_{L^2}
			+
			\| \Delta_j(f, g) \|_{L^1( 0,t ; L^2)}
			\right)
		\end{align}
		for all $j \in \mathbb{Z}$ with $|\Omega| \varepsilon < 2^j \leq \beta/\varepsilon$, provided that the right-hand side is finite.
        Here, $\beta_0$ denotes the positive constant appearing in Lemma \ref{lemm:high-ene}.
    \end{lemm}
\section{{Global} a priori estimates}\label{sec:a-priori}

In the present section, 
we establish several global a priori estimates of the solutions to \eqref{eq:NSC-2} by separating {them} into three cases: low, middle, and high{-}frequency parts. 
To begin with, we recall the local well-posedness for \eqref{eq:NSC-2}.
\begin{lemm}\label{lemm:LWP}
	Let $\Omega \in \mathbb{R}$, $\varepsilon>0$,
	and
	$(a_0,u_0) \in \dB_{2,1}^{\frac{3}{2}}(\mathbb{R}^3) \times \dB_{2,1}^{\frac{1}{2}}(\mathbb{R}^3)^3$ with $a_0(x) > - 1/ \varepsilon$ for all $x \in \mathbb{R}^3$.
	Then, there exists a positive time $T_{\Omega,\varepsilon}=T(\mu,\Omega,\varepsilon,a_0,u_0)$ such that
	the system \eqref{eq:NSC-2} possesses a unique solution $(a,u)$ in the class
	\begin{align}\label{LWP-class}
		&
        a \in C( [0,T_{\Omega,\varepsilon}] ; \dB_{2,1}^{\frac{3}{2}}(\mathbb{R}^3) ),\\
		&
        u \in C( [0,T_{\Omega,\varepsilon}] ; \dB_{2,1}^{\frac{1}{2}}(\mathbb{R}^3) )^3 \cap L^1 (0,T_{\Omega,\varepsilon} ; \dB_{2,1}^{\frac{5}{2}}(\mathbb{R}^3))^3
	\end{align}
	with $a(t,x) > -1/\varepsilon$ for all $(t,x) \in [0,T_{\Omega,\varepsilon}] \times \mathbb{R}^3$.
	Moreover, the solution admits the following regularity:
    \begin{align}\label{ad-reg-1}
        a \in \widetilde{C}( [0,T_{\Omega,\varepsilon}] ; \dB_{2,1}^{\frac{3}{2}}(\mathbb{R}^3) ),
        \quad
        u \in \widetilde{C}( [0,T_{\Omega,\varepsilon}] ; \dB_{2,1}^{\frac{1}{2}}(\mathbb{R}^3) )^3.
    \end{align}
    Furthermore, if we additionally assume $(a_0,u_0) \in \dB_{2,\infty}^{-\frac{3}{2}}(\mathbb{R}^3)$, then it holds 
    \begin{align}\label{ad-reg-2}
        (a,u) \in C( [0,T_{\Omega,\varepsilon}] ; \dB_{2,\infty}^{-\frac{3}{2}}(\mathbb{R}^3) )^{1+3}.
    \end{align}
\end{lemm}
\begin{proof}
    Although the proof follows the standard argument, we give its sketch for the readers' convenience.
    The unique existence of the local solutions in the class \eqref{LWP-class} was proved in \cite{Fuj-Wat-25}*{Proposition 4.1}. For the regularity \eqref{ad-reg-1}, it follows from Lemma \ref{lemm:low-ene} with $r=\infty$ and $s=1/2$ and \cite{Fuj-Wat-25}*{Lemma 2.5 with $s=1/2$} 
    \begin{align}
        &
        \varepsilon\n{a}_{
        \widetilde{L^{\infty}}(0,T_{\Omega,\varepsilon};\dB_{2,1}^{\frac{3}{2}})}
        +
        \n{(a,u)}_{\widetilde{L^{\infty}}(0,T_{\Omega,\varepsilon};\dB_{2,1}^{\frac{1}{2}})}
        \\
        \quad
        &\leq
        {}
        C\varepsilon
        \n{a_0}_{\dB_{2,1}^{\frac{3}{2}}}
        +
        C\n{u_0}_{\dB_{2,1}^{\frac{1}{2}}}
        +
        C
        \n{\sp{\div(au),N_{\varepsilon}[a,u]}}_{L^1(0,T_{\Omega,\varepsilon};\dB_{2,1}^{\frac{1}{2}})}
        \\
        &\qquad
        +
        C\int_0^{T_{\Omega,\varepsilon}}
        \varepsilon
        \n{a(t)}_{\dB_{2,1}^{\frac{3}{2}}}
        \n{\div u(t)}_{L^{\infty}}
        +
        \sum_{j \in \mathbb{Z}}
        2^{\frac{3}{2}j}
        \n{[\Delta_j u(t) \cdot \nabla ]a(t)}_{L^2}dt.
    \end{align} 
    Then, we may bound the nonlinear terms by the product of the norm of the solutions in the class \eqref{LWP-class}, which is exactly the same argument as in the case of the non-rotational compressible Navier--Stokes equations.

    Finally, we verify \eqref{ad-reg-2}.
    From Lemma \ref{lemm:simple-ene}, it immediately follows that 
    \begin{align}
        \n{(a,u)}_{L^{\infty}(0,T_{\Omega,\varepsilon};\dB_{2,\infty}^{-\frac{3}{2}})}
        \leq
        C
        \n{(a_0,u_0)}_{\dB_{2,\infty}^{-\frac{3}{2}}}
        +
        C
        \n{(\div(au),\mathcal{N}_{\varepsilon}[a,u])}_{\widetilde{L^1}(0,T_{\Omega,\varepsilon};\dB_{2,\infty}^{-\frac{3}{2}})}.
    \end{align}
    Following the completely same strategy as \cite{Xin-Xu-21}*{Section 5}, we deduce that $(a,u) \in L^{\infty}(0,T_{\Omega,\varepsilon};\dB_{2,\infty}^{-\frac{3}{2}}(\mathbb{R}^3))^{1+3}$.
    {It remains to prove $(a,u) \in C([0,T_{\Omega,\varepsilon}];\dB_{2,\infty}^{-\frac{3}{2}}(\mathbb{R}^3))^{1+3}$. 
    This is not obvious since the third index of the Besov space is infinity.
    In the following proof, we only consider the right continuity of $t \mapsto (a(t),u(t)) \in \dB_{2,\infty}^{-\frac{3}{2}}(\mathbb{R}^3)$ at $t=0$ since the continuity at $t>0$ may be shown similarly.
    Integrating \eqref{eq:NSC-2} with respect to the time variable, we see that 
    \begin{align}
        \begin{cases}
            a(t)-a_0 = \displaystyle\int_0^t -\dfrac{1}{\varepsilon}\div u(\tau)-\div (au)(\tau)d\tau, \\
            u(t)-u_0 = \displaystyle\int_0^t \mathcal{L}u(\tau)-\Omega e_3 \times u(\tau) - \dfrac{1}{\varepsilon} \nabla a(\tau) + \mathcal{N}_{\varepsilon}[a,u](\tau)d\tau.
        \end{cases}
    \end{align}
    First, we consider the continuity of $a$.
    It follows from the first equation above that 
    \begin{align}
        &
        \n{a(t)-a_0}_{\dB_{2,\infty}^{-\frac{3}{2}}} 
        \leq 
        C
        \varepsilon^{-1}
        t^{\frac{3}{4}}
        \n{u}_{L^4(0,T_{\Omega,\varepsilon};\dB_{2,\infty}^{-\frac{1}{2}})}
        +
        C
        t^{\frac{3}{4}}
        \n{au}_{L^4(0,T_{\Omega,\varepsilon};\dB_{2,\infty}^{-\frac{1}{2}})}
        \\
        &\quad 
        \leq 
        C
        t^{\frac{3}{4}}
        \sp{\varepsilon^{-1}+\n{a}_{L^{\infty}(0,T_{\Omega,\varepsilon};\dB_{2,1}^{\frac{3}{2}})}}
        \n{u}_{L^4(0,T_{\Omega,\varepsilon};\dB_{2,\infty}^{-\frac{1}{2}})}
        \\
        &\quad 
        \leq 
        C
        t^{\frac{3}{4}}
        \sp{\varepsilon^{-1}+\n{a}_{L^{\infty}(0,T_{\Omega,\varepsilon};\dB_{2,1}^{\frac{3}{2}})}}
        \n{u}_{L^{\infty}(0,T_{\Omega,\varepsilon};\dB_{2,\infty}^{-\frac{3}{2}})}^{\frac{3}{4}}
        \n{u}_{L^1(0,T_{\Omega,\varepsilon};\dB_{2,\infty}^{\frac{5}{2}})}^{\frac{1}{4}}
        \\
        &\quad \to 0
    \end{align}
    as $t \to +0$.
    For the continuity of $u$, we have 
    \begin{align}
        \n{u(t)-u_0}_{\dB_{2,\infty}^{-\frac{3}{2}}} 
        \leq {}&
        Ct\n{u}_{L^{\infty}(0,T_{\Omega,\varepsilon};\dB_{2,1}^{\frac{1}{2}})}
        +
        C|\Omega|t\n{u}_{L^{\infty}(0,T_{\Omega,\varepsilon};\dB_{2,\infty}^{-\frac{3}{2}})}
        \\
        &
        C\varepsilon^{-1}
        t
        \n{a}_{L^{\infty}(0,T_{\Omega,\varepsilon};\dB_{2,\infty}^{-\frac{1}{2}})}
        +
        C\n{\mathcal{N}_{\varepsilon}[a,u]}_{L^1(0,T_{\Omega,\varepsilon};\dB_{2,\infty}^{-\frac{3}{2}})}.
    \end{align}
    Here, we see by the bilinear estimates that 
    \begin{align}
        \n{\mathcal{N}_{\varepsilon}[a,u]}_{L^1(0,t;\dB_{2,\infty}^{-\frac{3}{2}})}
        \leq{}&
        \n{(u \cdot \nabla)u}_{L^1(0,t;\dB_{2,\infty}^{-\frac{3}{2}})}
        +
        \n{J(\varepsilon a) \mathcal{L}u}_{L^1(0,t;\dB_{2,\infty}^{-\frac{3}{2}})}
        \\
        &
        +
        \varepsilon^{-1}
        \n{K(\varepsilon a)\nabla a}_{L^1(0,t;\dB_{2,\infty}^{-\frac{3}{2}})}
        \\
        \leq{}&
        Ct
        \n{u}_{L^{\infty}(0,T_{\Omega,\varepsilon};\dB_{2,\infty}^{\frac{1}{2}})}
        \n{u}_{L^{\infty}(0,t;\dB_{2,1}^{\frac{1}{2}})}\\
        &
        +
        C
        \varepsilon
        t
        \n{a}_{L^{\infty}(0,T_{\Omega,\varepsilon};\dB_{2,1}^{\frac{3}{2}})}
        \n{u}_{L^{\infty}(0,T_{\Omega,\varepsilon};\dB_{2,\infty}^{\frac{1}{2}})}\\
        &
        +
        Ct
        \n{a}_{L^{\infty}(0,T_{\Omega,\varepsilon};\dB_{2,1}^{\frac{3}{2}})}
        \n{a}_{L^{\infty}(0,T_{\Omega,\varepsilon};\dB_{2,\infty}^{-\frac{1}{2}})}.
    \end{align}
    Hence, we have
    \begin{align}
        \n{u(t)-u_0}_{\dB_{2,\infty}^{-\frac{3}{2}}} 
        \leq{}
        &
        Ct\n{u}_{L^{\infty}(0,T_{\Omega,\varepsilon};\dB_{2,1}^{\frac{1}{2}})}
        +
        C|\Omega|t\n{u}_{L^{\infty}(0,T_{\Omega,\varepsilon};\dB_{2,\infty}^{-\frac{3}{2}})}
        \\
        &
        +
        C\varepsilon^{-1}
        t
        \n{a}_{L^{\infty}(0,T_{\Omega,\varepsilon};\dB_{2,\infty}^{-\frac{3}{2}}\cap \dB_{2,1}^{\frac{3}{2}})}
        +
        Ct
        \n{u}_{L^{\infty}(0,T_{\Omega,\varepsilon};\dB_{2,1}^{\frac{1}{2}})}^2\\
        &
        +
        C
        \varepsilon
        t
        \n{a}_{L^{\infty}(0,T_{\Omega,\varepsilon};\dB_{2,1}^{\frac{3}{2}})}
        \n{u}_{L^{\infty}(0,T_{\Omega,\varepsilon};\dB_{2,\infty}^{\frac{1}{2}})}\\
        &
        +
        Ct
        \n{a}_{L^{\infty}(0,T_{\Omega,\varepsilon};\dB_{2,1}^{\frac{3}{2}})}
        \n{a}_{L^{\infty}(0,T_{\Omega,\varepsilon};\dB_{2,\infty}^{-\frac{3}{2}}\cap \dB_{2,1}^{\frac{3}{2}})}
        \\
        \to{}& 0
    \end{align}
    as $t \to +0$.
    Hence, we complete the proof.}
\end{proof}

\subsection{The low{-}frequency estimates}
First of all, we focus on the low{-}frequency analysis, which is the most {harmful} part in this section.
\begin{lemm}\label{lemm:low-ene}
	Let $2 \leq q < 3$ and $2 < r \leq \infty$ satisfy $- 3 \slash q \leq - 3 \slash 2 + 4\slash r$.
	Then, there exists a positive constant $C=C(\mu,P,q,r)$ 
	such that
	if the solution $(a,u)$ to \eqref{eq:NSC-2} satisfies $\varepsilon \| a \|_{L^{\infty}(0,t;L^{\infty})} < 1$,
	then it holds
	\begin{align}
		&
		\begin{aligned}\label{lemm:low-ene-1}
			\| (a,u) \|_{L^{\infty}(0,t;\dB_{2,\infty}^{-\frac{3}{2}})}^{\ell; |\Omega| \varepsilon}
			\leq{}&
			C 
			\mathscr{D}^*_{\varepsilon}[a_0,u_0]\\
			&
			+
			C 
			\left( 
			1 + \varepsilon \| a \|_{\widetilde{L^{\infty}}( 0,t ; \dB_{q,1}^{\frac{3}{q}})} 
			\right)
			\| (a,u) \|_{{\widetilde{L^2}}(0,t;\dB_{q,\infty}^{\frac{3}{q}-1})}^2\\
			&
			+
			C
			\varepsilon \| a \|_{{\widetilde{L^{\infty}}}(0,t;\dB_{2,1}^{\frac{1}{2}})}
			\| u \|_{L^{\infty}(0,t; \dB_{2,\infty}^{ -\frac{1}{2} })},
		\end{aligned}\\
		&
		\begin{aligned}\label{lemm:low-ene-2}
			\| (a,u) \|_{{\widetilde{L^{r}}}(0,t;\dB_{q,\infty}^{\frac{3}{q}-3+\frac{4}{r}})}^{\ell; |\Omega| \varepsilon}
			\leq{}&
			C
			\left(
			|\Omega| 
			\varepsilon 
			\right)^{\frac{2}{r}}
			\mathscr{D}^*_{\varepsilon}[a_0,u_0]\\
			&
			+
			C 
			\left(
			|\Omega| 
			\varepsilon 
			\right)^{\frac{2}{r}}
			\left( 
			1 + \varepsilon \| a \|_{\widetilde{L^{\infty}}( 0,t ; \dB_{q,1}^{\frac{3}{q}})} 
			\right)
			\| (a,u) \|_{{\widetilde{L^2}}(0,t;\dB_{q,\infty}^{\frac{3}{q}-1})}^2\\
			&
			+
			C
			\varepsilon
			\| a \|_{{\widetilde{L^{\infty}}}(0,t;\dB_{2,1}^{\frac{1}{2}})}
			\| u \|_{{\widetilde{L^{r}}}(0,t;\dB_{q,\infty}^{\frac{3}{q}-2+\frac{4}{r}})},
		\end{aligned}\\
		&
			\begin{aligned}\label{lemm:low-ene-3}
				\| (a,u) \|_{{\widetilde{L^{r'}}}(0,t;\dB_{q,\infty}^{\frac{3}{q}-3+\frac{4}{r'}})}^{\ell; |\Omega| \varepsilon}
				\leq{}&
				C
				\left(
				|\Omega| 
				\varepsilon 
				\right)^{\frac{2}{r}}
				\mathscr{D}^*_{\varepsilon}[a_0,u_0]\\
				&
				+
				C 
				\left(
				|\Omega| 
				\varepsilon 
				\right)^{\frac{2}{r}}
				\left( 
				1 + \varepsilon \| a \|_{\widetilde{L^{\infty}}( 0,t ; \dB_{q,1}^{\frac{3}{q}})} 
				\right)
				\| (a,u) \|_{{\widetilde{L^2}}(0,t;\dB_{q,\infty}^{\frac{3}{q}-1})}^2\\
				&
				+
				C
				\varepsilon
				\| a \|_{{\widetilde{L^{\infty}}}(0,t;\dB_{q,\infty}^{\frac{3}{q}-\frac{2}{r}})}
				\| u \|_{{\widetilde{L^{r'}}}(0,t;\dB_{q,\infty}^{\frac{3}{q}-1+\frac{2}{r'}})}\\
				&
				+
				C
				\Omega^2
				\varepsilon^3
				\| a \|_{{\widetilde{L^2}}(0,t;\dB_{q,1}^{\frac{3}{q}})}
				\| u \|_{{\widetilde{L^2}}(0,t;\dB_{q,1}^{\frac{3}{q}}) \cap {\widetilde{L^{r}}}(0,t;\dB_{q,1}^{\frac{3}{q}-1+\frac{2}{r}})},
			\end{aligned} \quad\;
		\\
		&
		\begin{aligned}\label{lemm:low-ene-4}
			\| (a,u) \|_{{\widetilde{L^1}}(0,t;\dB_{2,\infty}^{\frac{5}{2}})}^{\ell; |\Omega| \varepsilon}
			\leq{}&
			C 
			\left(
			|\Omega| 
			\varepsilon 
			\right)^2
			\mathscr{D}^*_{\varepsilon}[a_0,u_0]\\
			&
			+
			C 
			\left(
			|\Omega| 
			\varepsilon 
			\right)^2
			\left( 
			1 + \varepsilon \| a \|_{\widetilde{L^{\infty}}( 0,t ; \dB_{q,1}^{\frac{3}{q}})} 
			\right)
			\| (a,u) \|_{{\widetilde{L^2}}(0,t;\dB_{q,\infty}^{\frac{3}{q}-1})}^2
		\end{aligned}
	\end{align}
	for all $\Omega \in \mathbb{R}$, $\varepsilon>0$ with $|\Omega| \varepsilon < \beta_0/\varepsilon$, and $0 < t< T_{\Omega,\varepsilon}^{\rm max}$,
	where $T_{\Omega,\varepsilon}^{\rm max}$ denotes the maximal existence time 
	and
	we have set
	\begin{align}
		\mathscr{D}^*_{\varepsilon}[a_0,u_0]
		:=
		\| (a_0, u_0) \|_{\dB_{2,\infty}^{-\frac{3}{2}}}
		+
		\varepsilon 
		\| a_0 \|_{\dB_{2,1}^{\frac{3}{2}}}
		\| u_0 \|_{\dB_{2,\infty}^{-\frac{3}{2}}}.
	\end{align}
\end{lemm}
 {
\begin{rem}\label{rem:low-ene}
    If we make use of Lemma \ref{lemm:ene-low} for the solution $(a,u)$ to \eqref{eq:NSC-2} to obtain the aforementioned estimates, we need the nonlinear estimates of the following type: 
    \begin{align}\label{prod:diff}
        \n{(u \cdot \nabla) u}_{L^1(0,t;\dB_{2,\infty}^{-\frac{3}{2}})}^{\ell;\Omega\varepsilon}
        \leq
        C
        \n{u}_{L^2(0,t;\dB_{q,\infty}^{{\frac{3}{q}-1}})}^2,
    \end{align}
    where $q$ should be $q>2$ in order to use the Strichartz type estimates for $u$ in the right{-}hand side. However, Lemma \ref{lemm:prod-1} implies that we need to take $q=2$ in \eqref{prod:diff}. 
    To {circumvent this back-and-forth situation}, we follow the idea mentioned in Section \ref{sec:intro} and {especially} consider the momentum $m:=\rho u =(1+\varepsilon a)u$ \textit{only} for the low{-}frequency part.
    Then, we see that $(a,m)$ solves \eqref{eq:m} below whose nonlinear terms have the divergence form. Then, the {aforementioned problem arising from} \eqref{prod:diff} is refined to the estimate \eqref{pf:low-ene-2} below, and all the other nonlinear terms {become to be} well-controlled.
\end{rem}}
\begin{proof}[Proof of Lemma \ref{lemm:low-ene}]
	Let $m := (1 + \varepsilon a) u$.
	Then, $(a,m)$ satisfies
	\begin{align}\label{eq:m}
		\begin{dcases}
			\partial_t a + \dfrac{1}{\varepsilon} \div m = 0, & t>0,x \in \mathbb{R}^3,\\
			\partial_t m - \mathcal{L} m + \Omega ( e_3 \times m ) + \dfrac{1}{\varepsilon} \nabla a = - \widetilde{\mathcal{N}_{\varepsilon}}[a,u], & t>0,x \in \mathbb{R}^3,\\
			a(0,x) = a_0(x),\quad m(0,x) = m_0(x),  & x \in \mathbb{R}^3,
		\end{dcases}
	\end{align}
	Here, we have put $m_0 := (1 + \varepsilon a_0) u_0$ and 
	\begin{gather}
		\widetilde{\mathcal{N}_{\varepsilon}}[a,u]
		:=
		\div( (1+\varepsilon a) u \otimes u) 
		+ \varepsilon \mathcal{L} (au) 
		+ \nabla(Q(\varepsilon a)a^2),\\
		Q(a) : = 
		\begin{dcases}
			\frac{P(1+a)- P(1) - a}{a^2}, \qquad & (a \neq 0), \\
			\frac{P''(1)}{2}, & (a=0).
		\end{dcases}
	\end{gather}
	Let $\zeta  {\, \in \{} 1,r',r,\infty  {\}}$.
	It follows from the second estimate in Lemma \ref{lemm:ene-low} with $\alpha = |\Omega|\varepsilon$ and $s=-3/2$
	that
	\begin{align}\label{pf:low-ene-1}
		\begin{split}
			\| (a,m) \|_{\widetilde{L^{\zeta}}(0,t;\dB_{2,\infty}^{-\frac{3}{2}+\frac{4}{\zeta}})}^{\ell; |\Omega| \varepsilon}
			\leq{}&
			C 
			\left( |\Omega| \varepsilon \right)^{\frac{2}{\zeta}}
			\| (a_0, m_0) \|_{\dB_{2,\infty}^{-\frac{3}{2}}}^{\ell; |\Omega| \varepsilon}\\
			&
			+
			C 
			\left( |\Omega| \varepsilon \right)^{\frac{2}{\zeta}}
			\| \widetilde{\mathcal{N}_{\varepsilon}}[a,u] \|_{{\widetilde{L^1}}(0,t;\dB_{2,\infty}^{-\frac{3}{2}})}^{\ell; |\Omega| \varepsilon}.
		\end{split}
	\end{align}
	For the estimate of  {$(a_0, m_0)$}, it follows from $m_0 = u_0 + \varepsilon a_0u_0$ and Lemma \ref{lemm:prod-1} that
	\begin{align}\label{pf:low-ene-6}
		\begin{split}
			\| (a_0,m_0) \|_{\dB_{2,\infty}^{-\frac{3}{2}}}^{\ell; |\Omega| \varepsilon}
			\leq
			\| (a_0,u_0) \|_{\dB_{2,\infty}^{-\frac{3}{2}}}^{\ell; |\Omega| \varepsilon}
			+
			C
			\varepsilon 
			\| a_0 \|_{\dB_{2,1}^{\frac{3}{2}}}
			\| u_0 \|_{\dB_{2,\infty}^{-\frac{3}{2}}}
			\leq
			C
			\mathscr{D}^*_{\varepsilon}[a_0,u_0].
		\end{split}
	\end{align}
	\par
	We  {next} focus on the {estimates} for the nonlinear terms.
	By Lemma \ref{lemm:prod-2} with $(p_1,p_2,s_1,s_2,r_1,r_2)=(q,2,3/q,-1/2,\infty,1)$ and Lemma \ref{lemm:compo} with $s_k=3/q-1$ and $r_k = 2$ ($k=1,2,3,4$), we have 
	\begin{align}\label{pf:low-ene-2}
		\begin{split}
			\| \div( (1+\varepsilon a) u \otimes u)  \|_{{\widetilde{L^1}}(0,t;\dB_{2,\infty}^{-\frac{3}{2}})}^{\ell; |\Omega| \varepsilon}
			\leq{}&
			C
			\| (1+\varepsilon a)( u \otimes u) \|_{{\widetilde{L^1}}(0,t;\dB_{2,\infty}^{-\frac{1}{2}})}\\
			\leq{}&
			C
			\left( 1 + \varepsilon \| a \|_{{\widetilde{L^{\infty}}}(0,t;\dB_{q,1}^{\frac{3}{q}})} \right)
			\| u \otimes u \|_{{\widetilde{L^1}}(0,t;\dB_{2,\infty}^{-\frac{1}{2}})} \quad \\
			\leq{}&
			C
			\left( 1 + \varepsilon \| a \|_{\widetilde{L^{\infty}}( 0,t ; \dB_{q,1}^{\frac{3}{q}})} \right)
			\| u \|_{{\widetilde{L^2}}(0,t;\dB_{q,\infty}^{\frac{3}{q}-1})}^2.
		\end{split}
	\end{align}
	By Lemma \ref{lemm:prod-2}, we have 
	\begin{align}\label{pf:low-ene-3}
		\begin{split}
			\| \varepsilon \mathcal{L} (au) \|_{{\widetilde{L^1}}(0,t;\dB_{2,\infty}^{-\frac{3}{2}})}^{\ell; |\Omega| \varepsilon}
			\leq{}&
			C
			|\Omega|
			\varepsilon^2
			\| au \|_{{\widetilde{L^1}}(0,t;\dB_{2,\infty}^{-\frac{1}{2}})}\\
			\leq{}&
			C\beta_0 
			\| a \|_{{\widetilde{L^2}}(0,t;\dB_{q,\infty}^{\frac{3}{q}-1})}
			\| u \|_{{\widetilde{L^2}}(0,t;\dB_{q,\infty}^{\frac{3}{q}-1})}.
		\end{split}
	\end{align}
	Let $\widetilde{Q}(a):=Q(a)-P''(1)/2$,  {which} is a smooth {function}  {with respect to $a$ satisfying} $\widetilde{Q}(0)=0$.
	By Lemmas \ref{lemm:prod-1}, \ref{lemm:prod-2}, and \ref{lemm:compo}, we see that
	\begin{align}\label{pf:low-ene-4}
		\begin{split}
			\| \nabla(Q(\varepsilon a)a^2) \|_{{\widetilde{L^1}}(0,t;\dB_{2,\infty}^{-\frac{3}{2}})}^{\ell; |\Omega| \varepsilon}
			\leq{}&
			C
			\|Q(\varepsilon a) a^2 \|_{{\widetilde{L^1}}(0,t;\dB_{2,\infty}^{-\frac{1}{2}})}^{ {\ell; |\Omega| \varepsilon}}\\
			\leq{}&
			C
			\left( \frac{1}{2}|P''(1)| +  \| \widetilde{Q}(\varepsilon a) \|_{{\widetilde{L^{\infty}}}(0,t;\dB_{q,1}^{\frac{3}{q}})} \right)
			\| a^2 \|_{{\widetilde{L^1}}(0,t;\dB_{2,\infty}^{-\frac{1}{2}})} \quad \quad \\
			\leq{}&
			C
			\left( 1 +  \varepsilon \| a \|_{{\widetilde{L^{\infty}}}(0,t;\dB_{q,1}^{\frac{3}{q}})} \right)
			\| a \|_{{\widetilde{L^2}}(0,t;\dB_{q,\infty}^{\frac{3}{q}-1})}^2.
		\end{split}
	\end{align}
	Combining \eqref{pf:low-ene-6}, \eqref{pf:low-ene-1}, \eqref{pf:low-ene-2}, \eqref{pf:low-ene-3}, and \eqref{pf:low-ene-4}, we obtain 
	\begin{align}\label{pf:low-ene-5}
		\begin{split}
			\| (a,m) \|_{\widetilde{L^{\zeta}}(0,t;\dB_{2,\infty}^{-\frac{3}{2}+\frac{4}{\zeta}})}^{\ell; |\Omega| \varepsilon}
			\leq{}&
			C 
			\left( |\Omega| \varepsilon \right)^{\frac{2}{\zeta}}
			\mathscr{D}^*_{\varepsilon}[a_0,u_0]\\
			&+
			C 
			\left( |\Omega| \varepsilon \right)^{\frac{2}{\zeta}}
			\left( 
			1 + \varepsilon \| a \|_{\widetilde{L^{\infty}}( 0,t ; \dB_{q,1}^{\frac{3}{q}})} 
			\right)
			\| (a,u) \|_{{\widetilde{L^2}}(0,t;\dB_{q,\infty}^{\frac{3}{q}-1})}^2.
		\end{split}
	\end{align}
	In the case of $\zeta=\infty$, we see by $u = m -\varepsilon au$ and Lemma \ref{lemm:prod-2} that
	\begin{align}\label{pf:low-ene-8}
		\begin{split}
			\| u \|_{L^{\infty}(0,t;\dB_{2,\infty}^{-\frac{3}{2}})}^{\ell; |\Omega| \varepsilon}
			\leq
			\| m \|_{L^{\infty}(0,t;\dB_{2,\infty}^{-\frac{3}{2}})}^{\ell; |\Omega| \varepsilon}
			+
			C \varepsilon 
			\| a \|_{{\widetilde{L^{\infty}}}(0,t;\dB_{2,1}^{\frac{1}{2}})}
			\| u \|_{L^{\infty}(0,t;\dB_{2,\infty}^{-\frac{1}{2}})}.
		\end{split}
	\end{align}
	Gathering \eqref{pf:low-ene-5} and \eqref{pf:low-ene-8}, we obtain \eqref{lemm:low-ene-1}.
	In the case of $\zeta = r$, it follows from $u = m -\varepsilon au$, the Bernstein inequality and Lemma \ref{lemm:prod-1} with $(p_1,p_2,s_1,s_2,r_1,r_2)=(2,q,1/2,3/q-2+4/r,\infty,r)$ that
	\begin{align}\label{pf:low-ene-7}
		\begin{split}
			\| u \|_{{\widetilde{L^{r}}}(0,t;\dB_{q,\infty}^{\frac{3}{q}-3+\frac{4}{r}})}^{\ell; |\Omega| \varepsilon}
			\leq{}&
			C
			\| m \|_{{\widetilde{L^{r}}}(0,t;\dB_{2,\infty}^{-\frac{3}{2}+\frac{4}{r}})}^{\ell; |\Omega| \varepsilon}
			+
			C
			\varepsilon
			\| a \|_{{\widetilde{L^{\infty}}}(0,t;\dB_{2,1}^{\frac{1}{2}})}
			\| u \|_{{\widetilde{L^{r}}}(0,t;\dB_{q,\infty}^{\frac{3}{q}-2+\frac{4}{r}})}.\ \ 
		\end{split}
	\end{align}
	By \eqref{pf:low-ene-5} and \eqref{pf:low-ene-7}, we have \eqref{lemm:low-ene-2}.
		For the case of $\zeta = r'$, we see that 
		\begin{align}\label{pf:low-ene-7-1}
			\begin{split}
				\n{u}_{{\widetilde{L^{r'}}}(0,t;\dB_{q,\infty}^{\frac{3}{q}-3+\frac{4}{r'}})}^{\ell; |\Omega| \varepsilon}
				\leq{}&
				C
				\| m \|_{{\widetilde{L^{r'}}}(0,t;\dB_{2,\infty}^{-\frac{3}{2}+\frac{4}{r'}})}^{\ell; |\Omega| \varepsilon}
				+
				\varepsilon
				\n{au}_{{\widetilde{L^{r'}}}(0,t;\dB_{q,\infty}^{\frac{3}{q}-3+\frac{4}{r'}})}^{\ell; |\Omega| \varepsilon},\\
				\leq{}&
				C
				\| m \|_{{\widetilde{L^{r'}}}(0,t;\dB_{2,\infty}^{-\frac{3}{2}+\frac{4}{r'}})}^{\ell; |\Omega| \varepsilon}\\
				&
				+
				C
				\varepsilon
				\| a \|_{{\widetilde{L^{\infty}}}(0,t;\dB_{q,\infty}^{\frac{3}{q}-\frac{2}{r}})}
				\| u \|_{{\widetilde{L^{r'}}}(0,t;\dB_{q,\infty}^{\frac{3}{q}-1+\frac{2}{r'}})}\\
				&
				+
				C
				\varepsilon
				\| a \|_{{\widetilde{L^2}}(0,t;\dB_{q,\infty}^{\frac{3}{q}+2})}^{\ell;4|\Omega|\varepsilon}
				\| u \|_{{\widetilde{L^{r^*}}}(0,t;\dB_{q,\infty}^{\frac{3}{q}-1+\frac{2}{r^*}})}\\
				\leq{}&
				C
				\| m \|_{{\widetilde{L^{r}}}(0,t;\dB_{2,\infty}^{-\frac{3}{2}+\frac{4}{r'}})}^{\ell; |\Omega| \varepsilon}\\
				&
				+
				C
				\varepsilon
				\| a \|_{{\widetilde{L^{\infty}}}(0,t;\dB_{q,\infty}^{\frac{3}{q}-\frac{2}{r}})}
				\| u \|_{{\widetilde{L^{r'}}}(0,t;\dB_{q,\infty}^{\frac{3}{q}-1+\frac{2}{r'}})}\\
				&
				+
				C
				\Omega^2
				\varepsilon^3
				\| a \|_{{\widetilde{L^2}}(0,t;\dB_{q,1}^{\frac{3}{q}})}
				\| u \|_{{\widetilde{L^2}}(0,t;\dB_{q,1}^{\frac{3}{q}}) \cap {\widetilde{L^{r}}}(0,t;\dB_{q,1}^{\frac{3}{q}-1+\frac{2}{r}})},
			\end{split}
		\end{align}
		where $r^*$ is defined by $1/r^*=1/2-1/r$.
		By \eqref{pf:low-ene-5} and \eqref{pf:low-ene-7-1}, {we see that} \eqref{lemm:low-ene-3} holds.
	In the case of $\zeta = 1$, it holds by $u = m -\varepsilon au$ and Lemma \ref{lemm:prod-2} that
	\begin{align}\label{pf:low-ene-10}
		\begin{split}
			\| u \|_{{\widetilde{L^1}}(0,t;\dB_{2,\infty}^{\frac{5}{2}})}^{\ell; |\Omega| \varepsilon}
			\leq{}&
			\| m \|_{{\widetilde{L^1}}(0,t;\dB_{2,\infty}^{\frac{5}{2}})}^{\ell; |\Omega| \varepsilon}
			+
			C
			\varepsilon
			( |\Omega| \varepsilon)^3
			\| au \|_{{\widetilde{L^1}}(0,t;\dB_{2,\infty}^{-\frac{1}{2}})}\\
			\leq{}&
			\| m \|_{{\widetilde{L^1}}(0,t;\dB_{2,\infty}^{\frac{5}{2}})}^{\ell; |\Omega| \varepsilon}
			+
			C
			\beta_0
			( |\Omega| \varepsilon)^2
			\| a \|_{{\widetilde{L^2}}(0,t;\dB_{q,\infty}^{\frac{3}{q}-1})}
			\| u \|_{{\widetilde{L^2}}(0,t;\dB_{q,\infty}^{\frac{3}{q}-1})}.
		\end{split}
	\end{align}
	Combining \eqref{pf:low-ene-5} and \eqref{pf:low-ene-10},
	we have \eqref{lemm:low-ene-4} and complete the proof.
\end{proof}
\subsection{The middle{-}frequency estimates}
 {
In this subsection, we focus on the a priori estimates of the {middle-}frequency part of the solutions.
In Lemmas \ref{lemm:middle-ene-1} and \ref{lemm:middle-ene-2}, we consider the energy estimates, and the Strichartz estimates are treated in Lemmas \ref{lemm:middle-str-1} and \ref{lemm:middle-str-2}.}
\begin{lemm}\label{lemm:middle-ene-1}
	Let $2 \leq q  < 4$ and $2 \leq r \leq \infty$ satisfy
	\begin{align}
		\frac{3}{q} - 1 + \frac{2}{r} > 0,
		\quad
		\frac{2}{r} \leq \frac{3}{q} - \frac{1}{2}.
	\end{align}
	Then, there exists a positive constant $C=C(\mu,P,p,q,r)$ such that
	if the solution $(a,u)$ to \eqref{eq:NSC-2} satisfies 
	$\varepsilon \| a \|_{L^{\infty}(0,t;L^{\infty})} < 1$, 
	then
	\begin{align}\label{lemm:middle-ene-1-pf1}
		\begin{split}
			&
			\| (a,u) \|_{\widetilde{L^{\infty}}(0,t;\dB_{2,1}^{\frac{1}{2}}) \cap L^1(0,t;\dB_{2,1}^{\frac{5}{2}})}^{m;\alpha,\frac{\beta_0}{\varepsilon}}
			\leq{}
			C\| (a_0,u_0) \|_{\dB_{2,1}^{\frac{1}{2}}}^{m;\alpha,\frac{\beta_0}{\varepsilon}}\\
			&\quad
			+
			C
			\| (a,u) \|_{\widetilde{L^2}(0,t ; \dB_{q,1}^\frac{3}{q})}^2
			+
			C 
			\varepsilon
			\| a \|_{\widetilde{L^{\infty}}( 0,t ; \dB_{q,1}^{\frac{3}{q}})}
			\| u \|_{L^1( 0,t ; \dB_{q,1}^{\frac{3}{q}+1})}^{h;\frac{\beta_0}{\varepsilon}}\\
			&
			\quad
			+
			C
			\| (a,u) \|_{\widetilde{L^{r}}(0,t;\dB_{q,1}^{\frac{3}{q}-1+\frac{2}{r}})}
			\left(
			\| a \|_{\widetilde{L^{r'}}(0,t;\dB_{q,1}^{\frac{3}{q}-1+\frac{2}{r'}})}^{\ell;\frac{4\beta_0}{\varepsilon}}
			+
			\| u \|_{\widetilde{L^{r'}}(0,t;\dB_{q,1}^{\frac{3}{q}-1+\frac{2}{r'}})}
			\right)
		\end{split}
	\end{align}
	for all $\alpha,\varepsilon>0$ with $|\Omega| \varepsilon \leq \alpha < \beta_0/\varepsilon$ and $0 < t< T_{\Omega,\varepsilon}^{\rm max}$,
	where $T_{\Omega,\varepsilon}^{\rm max}$ denotes the maximal existence time.
\end{lemm}
\begin{proof}
	If follows from Lemma \ref{lemm:low-ene-4} that 
	\begin{align}
		\| (a,u) \|_{\widetilde{L^{\infty}}(0,t;\dB_{2,1}^{\frac{1}{2}}) \cap L^1(0,t;\dB_{2,1}^{\frac{5}{2}})}^{m;\alpha,\frac{\beta_0}{\varepsilon}}
		\leq{}&
		C\| (a_0,u_0) \|_{\dB_{2,1}^{\frac{1}{2}}}^{m;\alpha,\frac{\beta_0}{\varepsilon}}\\
		&
		+
		C
		\| (\div(au),\mathcal{N}_{\varepsilon}[a,u]) \|_{L^1(0,t ; \dB_{2,1}^{\frac{1}{2}})}^{\ell;\frac{\beta_0}{\varepsilon}}.
	\end{align}
	For the nonlinear estimate, we see {from a} similar argument as in \cite{Fu-24} that 
	\begin{align}
		&
		\begin{aligned}\label{lemm:middle-ene-1-pf2}
			\begin{split}
				\| \div(au) \|_{L^1(0,t ; \dB_{2,1}^{\frac{1}{2}})}^{\ell;\frac{\beta_0}{\varepsilon}}
				\leq{}&
				C
				\| a \|_{\widetilde{L^r}(0,t ; \dB_{2,1}^{\frac{3}{q}-1+\frac{2}{r}})}
				\| u \|_{\widetilde{L^{r'}}(0,t ; \dB_{2,1}^{\frac{3}{q}-1+\frac{2}{r'}})}
				\\
				&+
				C
				\| a \|_{\widetilde{L^{r'}}(0,t ; \dB_{2,1}^{\frac{3}{q}-1+\frac{2}{r'}})}^{\ell;\frac{4\beta_0}{\varepsilon}}
				\| u \|_{\widetilde{L^r}(0,t ; \dB_{2,1}^{\frac{3}{q}-1+\frac{2}{r}})}^{\ell;\frac{\beta_0}{\varepsilon}},
			\end{split}
		\end{aligned}\\
		&
		\begin{aligned}\label{lemm:middle-ene-1-pf3}
			\begin{split}
				\| (u \cdot \nabla)u \|_{L^1(0,t ; \dB_{2,1}^{\frac{1}{2}})}^{\ell;\frac{\beta_0}{\varepsilon}}
				\leq{}
				C
				\| u \|_{\widetilde{L^r}(0,t ; \dB_{2,1}^{\frac{3}{q}-1+\frac{2}{r}})}
				\| u \|_{\widetilde{L^{r'}}(0,t ; \dB_{2,1}^{\frac{3}{q}-1+\frac{2}{r'}})},
			\end{split}
		\end{aligned}\\
		&
		\begin{aligned}\label{lemm:middle-ene-1-pf4}
			\begin{split}
				\| \mathcal{J}(\varepsilon a)\mathcal{L}u \|_{L^1(0,t ; \dB_{2,1}^{\frac{1}{2}})}^{\ell;\frac{\beta_0}{\varepsilon}}
				\leq{}&
				C
				\varepsilon
				\| a \|_{\widetilde{L^{\infty}}(0,t;\dB_{q,1}^{\frac{3}{q}})}
				\| u \|_{L^1(0,t;\dB_{q,1}^{\frac{3}{q}+1})}^{h;\frac{\beta_0}{\varepsilon}}\\
				&
				+
				C
				\| a \|_{\widetilde{L^2}(0,t;\dB_{q,1}^{\frac{3}{q}})}
				\| u \|_{\widetilde{L^2}(0,t;\dB_{q,1}^{\frac{3}{q}})}\\
				&
				+
				C
				\| a \|_{\widetilde{L^r}(0,t ; \dB_{2,1}^{\frac{3}{q}-1+\frac{2}{r}})}
				\| u \|_{\widetilde{L^{r'}}(0,t ; \dB_{2,1}^{\frac{3}{q}-1+\frac{2}{r'}})}
				\\
				&+
				C
				\| a \|_{\widetilde{L^{r'}}(0,t ; \dB_{2,1}^{\frac{3}{q}-1+\frac{2}{r'}})}^{\ell;\frac{4\beta_0}{\varepsilon}}
				\| u \|_{\widetilde{L^r}(0,t ; \dB_{2,1}^{\frac{3}{q}-1+\frac{2}{r}})}^{\ell;\frac{\beta_0}{\varepsilon}},
			\end{split}
		\end{aligned}\\
		&
		\begin{aligned}\label{lemm:middle-ene-1-pf5}
			\begin{split}
				\frac{1}{\varepsilon}\| \mathcal{K}(\varepsilon a)\nabla a \|_{L^1(0,t ; \dB_{2,1}^{\frac{1}{2}})}^{\ell;\frac{\beta_0}{\varepsilon}}
				\leq{}&
				C
				\| a \|_{\widetilde{L^2}(0,t;\dB_{q,1}^{\frac{3}{q}})}^2\\
				&
				+
				C
				\| a \|_{\widetilde{L^r}(0,t ; \dB_{2,1}^{\frac{3}{q}-1+\frac{2}{r}})}
				\| a \|_{\widetilde{L^{r'}}(0,t ; \dB_{2,1}^{\frac{3}{q}-1+\frac{2}{r'}})}^{\ell;\frac{4\beta_0}{\varepsilon}}.
			\end{split}
		\end{aligned}
	\end{align}
	Combining \eqref{lemm:middle-ene-1-pf1}, \eqref{lemm:middle-ene-1-pf2}, \eqref{lemm:middle-ene-1-pf3}, \eqref{lemm:middle-ene-1-pf4}, and \eqref{lemm:middle-ene-1-pf5}, we complete the proof.
\end{proof}
\begin{lemm}\label{lemm:middle-ene-2}
	Let $0<\varepsilon \leq 1$ and $\Omega \in \mathbb{R}$ satisfy $|\Omega|\varepsilon^2 < \beta_0$.
	Let $2 \leq q  < 3$.
	Then, there exists a positive constant $C=C(\mu,P,q)$ such that
	if the solution $(a,u)$ to \eqref{eq:NSC-2} satisfies $\varepsilon \| a \|_{L^{\infty}(0,t;L^{\infty})} < 1$, then
	\begin{align}
		&
		\| (a,u) \|_{\widetilde{L^{\infty}}(0,t; \dB_{2,1}^{-\frac{1}{2}})}^{\ell;\frac{\beta_0}{\varepsilon}}
		+
		\| a \|_{\widetilde{L^2}(0,t;\dB_{2,1}^{\frac{1}{2}})}^{m;|\Omega|\varepsilon,\frac{\beta_0}{\varepsilon}}
		+
		\| u \|_{\widetilde{L^2}(0,t;\dB_{2,1}^{\frac{1}{2}})}^{\ell;\frac{\beta_0}{\varepsilon}}\\
		&\quad
		\leq
		C
		\| (a_0,u_0) \|_{\dB_{2,1}^{-\frac{1}{2}}}^{\ell;\frac{\beta_0}{\varepsilon}}
		+
		C
		\| (a,u) \|_{\widetilde{L^2}(0,t;\dB_{q,\infty}^{\frac{3}{q}-1})}
		\| (a,u) \|_{\widetilde{L^2}(0,t;\dB_{q,1}^{\frac{3}{q}})}
	\end{align}
	for all $0 < t< T_{\Omega,\varepsilon}^{\rm max}$,
	where $T_{\Omega,\varepsilon}^{\rm max}$ denotes the maximal existence time.
\end{lemm}
\begin{proof}
	It follows from Lemma \ref{lemm:simple-ene} that
	\begin{align}\label{lemm:middle-ene-2-1}
		\begin{split}
			&\| (a,u) \|_{\widetilde{L^{\infty}}(0,t; \dB_{2,1}^{-\frac{1}{2}})}^{\ell;\frac{\beta_0}{\varepsilon}}
			+
			\| u \|_{\widetilde{L^2}(0,t;\dB_{2,1}^{\frac{1}{2}})}^{\ell;\frac{\beta_0}{\varepsilon}}\\
			&\quad
			\leq{}
			C
			\| (a_0,u_0) \|_{\dB_{2,1}^{-\frac{1}{2}}}^{\ell;\frac{\beta_0}{\varepsilon}}
			+
			C
			\| (\div(au), \mathcal{N}_{\varepsilon}[a,u]) \|_{L^1(0,t;\dB_{2,1}^{-\frac{1}{2}})}^{\ell;\frac{\beta_0}{\varepsilon}}.
		\end{split}
	\end{align}
	By Lemma \ref{lemm:ene-low}, it holds  
	\begin{align}\label{lemm:middle-ene-2-2}
		\begin{split}
			\| a \|_{\widetilde{L^2}(0,t;\dB_{2,1}^{\frac{1}{2}})}^{m;|\Omega|\varepsilon,\frac{\beta_0}{\varepsilon}}
			\leq{}&
			C
			\| (a_0,u_0) \|_{\dB_{2,1}^{-\frac{1}{2}}}^{m;|\Omega|\varepsilon,\frac{\beta_0}{\varepsilon}}
			+
			C
			\| (\div(au), \mathcal{N}_{\varepsilon}[a,u]) \|_{L^1(0,t;\dB_{2,1}^{-\frac{1}{2}})}^{m;|\Omega|\varepsilon,\frac{\beta_0}{\varepsilon}}\\
			\leq{}&
			C
			\| (a_0,u_0) \|_{\dB_{2,1}^{-\frac{1}{2}}}^{\ell;\frac{\beta_0}{\varepsilon}}
			+
			C
			\| (\div(au), \mathcal{N}_{\varepsilon}[a,u]) \|_{L^1(0,t;\dB_{2,1}^{-\frac{1}{2}})}^{\ell;\frac{\beta_0}{\varepsilon}}.
		\end{split}
	\end{align}
	We focus on the nonlinear estimates.
	From Lemma \ref{lemm:prod-2} with $s_1 =s_4 = 3/q - 1$ and $s_2 = s_3 = 3/q$, we see that
	\begin{align}\label{lemm:middle-ene-2-3}
		\begin{split}
			\| \div(au) \|_{L^1(0,t;\dB_{2,1}^{-\frac{1}{2}})}^{\ell;\frac{\beta_0}{\varepsilon}}
			\leq{}&
			C
			\| au \|_{L^1(0,t;\dB_{2,1}^{\frac{1}{2}})}\\
			\leq{}&
			C
			\| a \|_{\widetilde{L^2}(0,t;\dB_{q,\infty}^{\frac{3}{q}-1})}
			\| u \|_{\widetilde{L^2}(0,t;\dB_{q,1}^{\frac{3}{q}})}\\
			&
			+
			C
			\| u \|_{\widetilde{L^2}(0,t;\dB_{q,\infty}^{\frac{3}{q}-1})}
			\| a \|_{\widetilde{L^2}(0,t;\dB_{q,1}^{\frac{3}{q}})}\\
			\leq{}&
			C
			\| (a,u) \|_{\widetilde{L^2}(0,t;\dB_{q,\infty}^{\frac{3}{q}-1})}
			\| (a,u) \|_{\widetilde{L^2}(0,t;\dB_{q,1}^{\frac{3}{q}})}.
		\end{split}
	\end{align}
	By Lemma \ref{lemm:prod-2} with $s_1 = s_2 = s_3 = s_4 = 3/q - 1$, we have
	\begin{align}\label{lemm:middle-ene-2-4}
		\begin{split}
			\| (u\cdot \nabla)u \|_{L^1(0,t;\dB_{2,1}^{-\frac{1}{2}})}^{\ell;\frac{\beta_0}{\varepsilon}}
			\leq{}&
			C 
			\| u \|_{\widetilde{L^2} ( 0,t ; \dB_{q,\infty}^{\frac{3}{q}-1})}
			\| \nabla u \|_{\widetilde{L^2} ( 0,t ; \dB_{q,1}^{\frac{3}{q}-1})}\\
			\leq{}&
			C 
			\| u \|_{\widetilde{L^2} ( 0,t ; \dB_{q,\infty}^{\frac{3}{q}-1})}
			\| u \|_{\widetilde{L^2} ( 0,t ; \dB_{q,1}^{\frac{3}{q}})}.
		\end{split}
	\end{align}
	We also see by Lemmas \ref{lemm:prod-2} with $s_1 = s_2 = 3/q - 1$, $s_3 = 3/q$, and $s_4 = 3/q - 2$ and \ref{lemm:compo} that
	\begin{align}\label{lemm:middle-ene-2-5}
		\begin{split}
			\| J(\varepsilon a) \mathcal{L}u \|_{L^1(0,t;\dB_{2,1}^{-\frac{1}{2}})}^{\ell;\frac{\beta_0}{\varepsilon}}
			\leq{}&
			C
			\| J(\varepsilon a) \|_{\widetilde{L^2}( 0,t ; \dB_{q,\infty}^{\frac{3}{q}-1})}
			\| \mathcal{L} u \|_{\dB_{q,1}^{\frac{3}{q}-1}}^{\frac{4\beta_0}{\varepsilon}}\\
			&
			+
			C
			\| J(\varepsilon a)\|_{\widetilde{L^2}( 0,t ; \dB_{q,1}^{\frac{3}{q}})}
			\|\mathcal{L}u\|_{\widetilde{L^2}( 0,t ; \dB_{q,1}^{\frac{3}{q}-2})}\\
			\leq{}&
			C
			\| a \|_{\widetilde{L^2}( 0,t ; \dB_{q,\infty}^{\frac{3}{q}-1})}
			\| u \|_{\widetilde{L^2}( 0,t ; \dB_{q,1}^{\frac{3}{q}})}\\
			&
			+
			C\varepsilon 
			\| a \|_{\widetilde{L^2}( 0,t ; \dB_{q,1}^{\frac{3}{q}})}
			\| u \|_{\widetilde{L^2}( 0,t ; \dB_{q,1}^{\frac{3}{q}})}.
		\end{split}
	\end{align}
	Similarly, we have 
	\begin{align}\label{lemm:middle-ene-2-6}
		\begin{split}
			\frac{1}{\varepsilon}\| K(\varepsilon a) \nabla a \|_{L^1(0,t;\dB_{2,1}^{-\frac{1}{2}})}^{\ell;\frac{\beta_0}{\varepsilon}}
			\leq{}&
			\frac{C}{\varepsilon}
			\| K(\varepsilon a) \|_{\widetilde{L^2}( 0,t ; \dB_{q,\infty}^{\frac{3}{q}-1})}
			\| \nabla a \|_{\widetilde{L^2}( 0,t ; \dB_{q,1}^{\frac{3}{q}-1})}^{\frac{4\beta_0}{\varepsilon}}\\
			&
			+
			\frac{C}{\varepsilon}
			\|\nabla a\|_{\widetilde{L^2}( 0,t ; \dB_{q,\infty}^{\frac{3}{q}-2})}
			\| K(\varepsilon a) \|_{\widetilde{L^2}( 0,t ; \dB_{q,1}^{\frac{3}{q}})}\\
			\leq{}&
			C
			\| a \|_{\widetilde{L^2}( 0,t ; \dB_{q,\infty}^{\frac{3}{q}-1})}
			\| a \|_{\widetilde{L^2}( 0,t ; \dB_{q,1}^{\frac{3}{q}})}.
		\end{split}
	\end{align}
	 {Collecting} 
	\eqref{lemm:middle-ene-2-1}, 
	\eqref{lemm:middle-ene-2-2}, 
	\eqref{lemm:middle-ene-2-3}, 
	\eqref{lemm:middle-ene-2-4}, 
	\eqref{lemm:middle-ene-2-5}, and 
	\eqref{lemm:middle-ene-2-6}, 
	 {the proof is complete}.
\end{proof}
\begin{lemm}\label{lemm:middle-str-1}
	Let 
	$\Omega \in \mathbb{R} \setminus \{ 0 \}$ and 
	$\alpha, \varepsilon>0$ 
	satisfy 
	$|\Omega|\varepsilon < \alpha < \beta_0/\varepsilon$.
	Let $2 \leq q \leq 4$ and $2 \leq r \leq \infty$ satisfy
	\begin{align}
		\label{cond-lemm7.4}
		\frac{1}{q} + \frac{1}{r} \leq \frac{1}{2},\quad
		\frac{2}{r} \leq \frac{3}{q} - \frac{1}{2}.
	\end{align}
	Then, there exists a positive constant $C=C(\mu,P,q,r)$ such that
	if the solution $(a,u)$ to \eqref{eq:NSC-2} satisfies $\varepsilon \| a \|_{L^{\infty}(0,t;L^{\infty})} < 1$, then
	\begin{align}
		&
		\| (a,u) \|_{\widetilde{L^r}( 0,t ; \dB_{q,1}^{\frac{3}{q}-1+\frac{2}{r}})}^{m;|\Omega|\varepsilon,\alpha}
		\leq{}
		C
		\alpha^{\frac{2}{r}}
		|\Omega|^{-\frac{1}{r}}
		\| (a_0,u_0) \|_{\dB_{2,1}^{\frac{1}{2}}}\\
		&\quad
		+
		C
		\alpha^{\frac{2}{r}}
		|\Omega|^{-\frac{1}{r}}
		C\| (a,u) \|_{\widetilde{L^2}(0,t ; \dB_{q,1}^\frac{3}{q})}^2
		+
		C 
		\alpha^{\frac{2}{r}}
		|\Omega|^{-\frac{1}{r}}
		\varepsilon
		\| a \|_{\widetilde{L^{\infty}}( 0,t ; \dB_{q,1}^{\frac{3}{q}})}
		\| u \|_{L^1( 0,t ; \dB_{q,1}^{\frac{3}{q}+1})}^{h;\frac{\beta_0}{\varepsilon}}\\
		&
		\quad
		+
		C
		\alpha^{\frac{2}{r}}
		|\Omega|^{-\frac{1}{r}}
		\| (a,u) \|_{\widetilde{L^{r}}(0,t;\dB_{q,1}^{\frac{3}{q}-1+\frac{2}{r}})}
		\left(
		\| a \|_{\widetilde{L^{r'}}(0,t;\dB_{q,1}^{\frac{3}{q}-1+\frac{2}{r'}})}^{\ell;\frac{4\beta_0}{\varepsilon}}
		+
		\| u \|_{\widetilde{L^{r'}}(0,t;\dB_{q,1}^{\frac{3}{q}-1+\frac{2}{r'}})}
		\right)
	\end{align}
	for all $0 < t< T_{\Omega,\varepsilon}^{\rm max}$,
	where $T_{\Omega,\varepsilon}^{\rm max}$ denotes the maximal existence time.
\end{lemm}
\begin{proof}
	It follows from Lemma \ref{lemm:lin-str-visc} that 
	\begin{align}
		&
		\| (a,u) \|_{\widetilde{L^r}( 0,t ; \dB_{q,1}^{\frac{3}{q}-1+\frac{2}{r}})}^{\ell;\frac{\beta_0}{\varepsilon}}\\
		&\quad
		\leq{}
		C|\Omega|^{-\frac{1}{r}}
		\| (a_0,u_0) \|_{\dB_{2,1}^{\frac{1}{2}+\frac{2}{r}}}^{m;|\Omega|\varepsilon,\alpha}
		+
		C
		|\Omega|^{-\frac{1}{r}}
		\| (\div(au),\mathcal{N}_{\varepsilon}[a,u]) \|_{L^1(0,t;\dB_{2,1}^{\frac{1}{2}+\frac{2}{r}})}^{m;|\Omega|\varepsilon,\alpha}
		\\
		&\quad\leq{}
		C
		\alpha^{\frac{2}{r}}
		|\Omega|^{-\frac{1}{r}}
		\| (a_0,u_0) \|_{\dB_{2,1}^{\frac{1}{2}}}
		+
		C
		\alpha^{\frac{2}{r}}
		|\Omega|^{-\frac{1}{r}}
		\| (\div(au),\mathcal{N}_{\varepsilon}[a,u]) \|_{L^1(0,t;\dB_{2,1}^{\frac{1}{2}})}^{\ell;\frac{\beta_0}{\varepsilon}}.
	\end{align}
	Here, the estimates of the nonlinear terms are obtained by the same argument as {in} Lemma \ref{lemm:middle-ene-1}; see \eqref{lemm:middle-ene-1-pf2}, \eqref{lemm:middle-ene-1-pf3}, \eqref{lemm:middle-ene-1-pf4}, and \eqref{lemm:middle-ene-1-pf5}.
	Thus, we complete the proof.
\end{proof}
\begin{lemm}\label{lemm:middle-str-2}
	Let
	$\Omega \in \mathbb{R} \setminus \{ 0 \}$ and $0 < \varepsilon \leq 1$ 
	satisfy 
	$|\Omega|\varepsilon^2 < \beta_0$. 
	Let $2 \leq q < 3$ and $2 \leq r \leq \infty$ satisfy \eqref{cond-lemm7.4}.
	Then, there exists a positive constant $C=C(\mu,P,q,r)$ such that
	if the solution $(a,u)$ to \eqref{eq:NSC-2} satisfies $\varepsilon \| a \|_{L^{\infty}(0,t;L^{\infty})} < 1$, then
	\begin{align}
		&
		\| (a,u) \|_{\widetilde{L^{r}}( 0,t ; \dB_{q,1}^{\frac{3}{q}-3+\frac{4}{r}})}^{m;|\Omega|\varepsilon,\frac{\beta_0}{\varepsilon}}
		\leq{}
		C
		|\Omega|^{-\frac{1}{r}}
		\| (a_0,u_0) \|_{\dB_{2,\infty}^{-\frac{3}{2}} \cap \dB_{2,1}^{\frac{1}{2}}}\\
		&\quad
		+
		C
		|\Omega|^{-\frac{1}{r}}
		\left( 1 + \varepsilon \| a \|_{\widetilde{L^{\infty}}( 0,t ; \dB_{q,1}^{\frac{3}{q}})} \right)
		\left(
		\| a \|_{\widetilde{L^2}( 0,t ; \dB_{2,\infty}^{\frac{1}{2}})}
		+
		\| u \|_{\widetilde{L^2}( 0,t ; \dB_{2,1}^{\frac{1}{2}})}
		\right)^2\\
		&\quad
		+
		C
		|\Omega|^{-\frac{1}{r}}
		\| (a,u) \|_{\widetilde{L^2}(0,t ; \dB_{2,1}^\frac{3}{2})}^2
		+
		C 
		|\Omega|^{-\frac{1}{r}}
		\varepsilon
		\| a \|_{L^{\infty}( 0,t ; \dB_{q,1}^{\frac{3}{q}})}
		\| u \|_{L^1( 0,t ; \dB_{q,1}^{\frac{3}{q}+1})}^{h;\frac{\beta_0}{\varepsilon}}\\
		&
		\quad
		+
		C
		|\Omega|^{-\frac{1}{r}}
		\| (a,u) \|_{\widetilde{L^{r}}(0,t;\dB_{q,1}^{\frac{3}{q}-1+\frac{2}{r}})}
		\left(
		\| a \|_{\widetilde{L^{r'}}(0,t;\dB_{q,1}^{\frac{3}{q}-1+\frac{2}{r'}})}^{\ell;\frac{4\beta_0}{\varepsilon}}
		+
		\| u \|_{\widetilde{L^{r'}}(0,t;\dB_{q,1}^{\frac{3}{q}-1+\frac{2}{r'}})}
		\right)
	\end{align}
	for all $0 < t< T_{\Omega,\varepsilon}^{\rm max}$,
	where $T_{\Omega,\varepsilon}^{\rm max}$ denotes the maximal existence time.
\end{lemm}
\begin{proof}
	By Lemma \ref{lemm:lin-str-visc}, we have
	\begin{align}\label{pf::middle-str-2-1}
		\begin{split}
			\| (a,u) \|_{\widetilde{L^r}( 0,t ; \dB_{q,\infty}^{\frac{3}{q}-3+\frac{4}{r}})}^{m;|\Omega|\varepsilon,\frac{\beta_0}{\varepsilon}}
			\leq{}&
			C|\Omega|^{-\frac{1}{r}}
			\| (a_0,u_0) \|_{\dB_{2,\infty}^{-\frac{3}{2}+\frac{4}{r}}}^{m;|\Omega|\varepsilon,\frac{\beta_0}{\varepsilon}}\\
			&+
			C
			|\Omega|^{-\frac{1}{r}}
			\| (\div(au),\mathcal{N}_{\varepsilon}[a,u]) \|_{\widetilde{L^1}(0,t;\dB_{2,\infty}^{-\frac{3}{2}+\frac{4}{r}})}^{m;|\Omega|\varepsilon,\frac{\beta_0}{\varepsilon}}.
		\end{split}
	\end{align}
	Since $-3/2 \leq -3/2+4/r \leq 1/2$, it holds
	\begin{align}\label{pf::middle-str-2-2-1}
		\| (a_0,u_0) \|_{\dB_{2,\infty}^{-\frac{3}{2}+\frac{4}{r}}}^{m;|\Omega|\varepsilon,\frac{\beta_0}{\varepsilon}}
		\leq
		C
		\| (a_0,u_0) \|_{\dB_{2,\infty}^{-\frac{3}{2}} \cap \dB_{2,1}^{\frac{1}{2}}},
	\end{align}
	and
	\begin{align}\label{pf::middle-str-2-2-2}
		\begin{split}
			\| (\div(au),\mathcal{N}_{\varepsilon}[a,u]) \|_{\widetilde{L^1}(0,t;\dB_{2,\infty}^{-\frac{3}{2}+\frac{4}{r}})}^{m;|\Omega|\varepsilon,\frac{\beta_0}{\varepsilon}}
			\leq{}&
			C
			\| (\div(au),\mathcal{N}_{\varepsilon}[a,u]) \|_{\widetilde{L^1}(0,t;\dB_{2,\infty}^{-\frac{3}{2}})}^{\ell;\frac{\beta_0}{\varepsilon}}\\
			&
			+
			C
			\| (\div(au),\mathcal{N}_{\varepsilon}[a,u]) \|_{L^1(0,t;\dB_{2,1}^{\frac{1}{2}})}^{\ell;\frac{\beta_0}{\varepsilon}}.
		\end{split}
	\end{align}
	As the {estimate} for the second term of the right{-}hand side in \eqref{pf::middle-str-2-2-2} is obtained by \eqref{lemm:middle-ene-1-pf1}, \eqref{lemm:middle-ene-1-pf2}, \eqref{lemm:middle-ene-1-pf3}, \eqref{lemm:middle-ene-1-pf4}, and \eqref{lemm:middle-ene-1-pf5} in the proof of Lemma \ref{lemm:middle-ene-1},  {we have}
	\begin{align}\label{pf::middle-str-2-8}
		\begin{split}
			&
			\| (\div(au),\mathcal{N}_{\varepsilon}[a,u]) \|_{\widetilde{L^1}(0,t;\dB_{2,1}^{ {\frac12}} )}^{\ell;\frac{\beta_0}{\varepsilon}}\\
			&\quad
			\leq
			C
			\| (a,u) \|_{\widetilde{L^2}(0,t ; \dB_{q,1}^\frac{3}{q})}^2
			+
			C 
			\varepsilon
			\| a \|_{L^{\infty}( 0,t ; \dB_{q,1}^{\frac{3}{q}})}
			\| u \|_{L^1( 0,t ; \dB_{q,1}^{\frac{3}{q}+1})}^{h;\frac{\beta_0}{\varepsilon}}\\
			&
			\qquad
			+
			C
			\| (a,u) \|_{\widetilde{L^{r}}(0,t;\dB_{q,1}^{\frac{3}{q}-1+\frac{2}{r}})}
			\left(
			\| a \|_{\widetilde{L^{r'}}(0,t;\dB_{q,1}^{\frac{3}{q}-1+\frac{2}{r'}})}^{\ell;\frac{4\beta_0}{\varepsilon}}
			+
			\| u \|_{\widetilde{L^{r'}}(0,t;\dB_{q,1}^{\frac{3}{q}-1+\frac{2}{r'}})}
			\right).
		\end{split}
	\end{align}
	Thus, {it} {remains to} consider the estimates of the first term of the right{-}hand side  {of \eqref{pf::middle-str-2-2-2}}.
	From Lemma \ref{lemm:prod-2}, it follows that
	\begin{align}\label{pf::middle-str-2-3}
		\| \div(au) \|_{\widetilde{L^1}(0,t;\dB_{2,\infty}^{-\frac{3}{2}})}^{\ell;\frac{\beta_0}{\varepsilon}}
		\leq
		C
		\| au \|_{\widetilde{L^1}(0,t;\dB_{2,\infty}^{-\frac{1}{2}})}
		\leq
		C 
		\| a \|_{\widetilde{L^2}(0,t ; \dB_{2,\infty}^{\frac{1}{2}})}
		\| u \|_{\widetilde{L^2}(0,t ; \dB_{2,\infty}^{\frac{1}{2}})}.
	\end{align}
	It holds by Lemma \ref{lemm:prod-1} that
	\begin{align}\label{pf::middle-str-2-4}
		\| (u \cdot \nabla) u \|_{\widetilde{L^1}(0,t;\dB_{2,\infty}^{-\frac{3}{2}})}^{\ell;\frac{\beta_0}{\varepsilon}}
		\leq{}
		C\| u \|_{\widetilde{L^2}( 0,t ;\dB_{2,1}^{\frac{1}{2}})}
		\| \nabla u \|_{\widetilde{L^2}(0,t ;\dB_{2,\infty}^{-\frac{1}{2}})}
		\leq{}
		C\| u \|_{\widetilde{L^2}( 0,t ;\dB_{2,1}^{\frac{1}{2}})}^2.
	\end{align}
	Using Lemmas \ref{lemm:prod-1} and \ref{lemm:compo}, there holds
	\begin{align}\label{pf::middle-str-2-5}
		\begin{split}
			\| J(\varepsilon a) \mathcal{L}u \|_{\widetilde{L^1}(0,t;\dB_{2,\infty}^{-\frac{3}{2}})}^{\ell;\frac{\beta_0}{\varepsilon}}
			\leq{}&
			C
			\| J(\varepsilon a)\|_{\widetilde{L^2}(0,t;\dB_{2,\infty}^{\frac{1}{2}})}
			\|\mathcal{L}u \|_{\widetilde{L^2}(0,t;\dB_{2,1}^{-\frac{1}{2}})}\\
			\leq{}&
			C\varepsilon
			\| a \|_{\widetilde{L^2}(0,t;\dB_{2,\infty}^{\frac{1}{2}})}
			\| u \|_{\widetilde{L^2}(0,t;\dB_{2,1}^{\frac{3}{2}})}\\
			\leq{}&
			C
			\| a \|_{\widetilde{L^2}(0,t;\dB_{2,\infty}^{\frac{1}{2}})}^2
			+
			C
			\| u \|_{\widetilde{L^2}(0,t;\dB_{2,1}^{\frac{3}{2}})}^2.
		\end{split}
	\end{align}
	Here, we see that
	\begin{align}
		\frac{1}{\varepsilon}K(\varepsilon a) \nabla a
		=
		\frac{1}{\varepsilon^2}
		\nabla(G(\varepsilon a)),\qquad 
		G(a)
		= 
		\int_0^{a} K(s) ds.
	\end{align}
	Since $G$ is a smooth function satisfying $G(0) = G'(0) = 0$, there exists a smooth function $H$ satisfying
	$G(a) = (G''(0)\slash2 + H(a)) a^2$ and $H(0) = 0$.
	Thus, we infer from Lemmas \ref{lemm:prod-2} and \ref{lemm:compo} that
	\begin{align}\label{pf::middle-str-2-6}
		\begin{split}
			\frac{1}{\varepsilon}
			\| K(\varepsilon a) \nabla a \|_{\widetilde{L^1}(0,t;\dB_{2,\infty}^{-\frac{3}{2}})}^{\ell;\frac{\beta_0}{\varepsilon}}
			\leq {}&
			\frac{C}{\varepsilon^2}
			\| G(\varepsilon a) \|_{\widetilde{L^1}(0,t;\dB_{2,\infty}^{-\frac{1}{2}})}\\
			\leq{}&
			C
			\left( \frac{1}{2}|G''(0)| + \| H(\varepsilon a) \|_{\widetilde{L^{\infty}}( 0,t ; \dB_{q,1}^{\frac{3}{q}})} \right)\| a ^2 \|_{\widetilde{L^1}(0,t;\dB_{2,\infty}^{-\frac{1}{2}})} \quad \enskip \\
			\leq {}&
			C 
			\left( 1 + \varepsilon \| a \|_{\widetilde{L^{\infty}}( 0,t ; \dB_{q,1}^{\frac{3}{q}})} \right)
			\| a^2 \|_{\widetilde{L^1}(0,t;\dB_{2,\infty}^{-\frac{1}{2}})}\\
			\leq {}&
			C
			\left( 1 + \varepsilon \| a \|_{\widetilde{L^{\infty}}( 0,t ; \dB_{q,1}^{\frac{3}{q}})} \right)
			\| a \|_{\widetilde{L^2}(0,t ; \dB_{2,\infty}^{\frac{1}{2}})}^2.
		\end{split}
	\end{align}
	It follows from \eqref{pf::middle-str-2-3}, \eqref{pf::middle-str-2-4}, \eqref{pf::middle-str-2-5}, and \eqref{pf::middle-str-2-6} that
	\begin{align}\label{pf::middle-str-2-7}
		\begin{split}
			&
			\| (\div(au),\mathcal{N}_{\varepsilon}[a,u]) \|_{L^1(0,t;\dB_{2,\infty}^{ {- \frac32}})}^{\ell;\frac{\beta_0}{\varepsilon}}
			\leq
			C
			\| (a,u) \|_{\widetilde{L^2}(0,t ; \dB_{2,1}^\frac{3}{2})}^2\\
			&\quad
			+
			C
			\left( 1 + \varepsilon \| a \|_{\widetilde{L^{\infty}}( 0,t ; \dB_{q,1}^{\frac{3}{q}})} \right)
			\left(
			\| a \|_{\widetilde{L^2}( 0,t ; \dB_{2,\infty}^{\frac{1}{2}})}
			+
			\| u \|_{\widetilde{L^2}( 0,t ; \dB_{2,1}^{\frac{1}{2}})}
			\right)^2.
		\end{split}
	\end{align}
	Combining \eqref{pf::middle-str-2-1}, \eqref{pf::middle-str-2-2-1}, \eqref{pf::middle-str-2-2-2}, \eqref{pf::middle-str-2-8}, and \eqref{pf::middle-str-2-7}, we complete the proof.
\end{proof}
\subsection{The high{-}frequency estimates}
The following lemma provides the high{-}frequency estimates.
\begin{lemm}\label{lemm:a-propri-high}
	Let $\Omega \in \mathbb{R}$ and $\varepsilon>0$ satisfy $|\Omega| \varepsilon^2 \leq 1$.
	Let $1 \leq p,q < \infty$ satisfy $1\slash3 < 1\slash p + 1 \slash q \leq 1$.
	Then, there exist positive constants $\beta_0 = \beta_0(\mu)$ and $C=C(\mu,P,p,q)$ such that
	if the solution $(a,u)$ to \eqref{eq:NSC-2} satisfies  $\varepsilon \| a \|_{L^{\infty}(0,t;L^{\infty})} < 1$, then
	\begin{align}
		&
		\varepsilon \| a \|_{\widetilde{L^{\infty}}( 0,t ; \dB_{p,1}^{\frac{3}{p}})}^{h;\frac{\beta_0}{\varepsilon}}
		+
		\frac{1}{\varepsilon} \| a \|_{L^1( 0,t ; \dB_{p,1}^{\frac{3}{p}})}^{h;\frac{\beta_0}{\varepsilon}}
		+
		\| u \|_{\widetilde{L^{\infty}}( 0,t ; \dB_{p,1}^{\frac{3}{p}-1}) \cap L^1( 0,t ; \dB_{p,1}^{\frac{3}{p}+1})}^{h;\frac{\beta_0}{\varepsilon}}\\
		&\quad
		\leq
		C
		\left( \varepsilon \| a_0 \|_{\dB_{p,1}^{\frac{3}{p}}}^{h;\frac{\beta_0}{\varepsilon}} + \| u_0 \|_{\dB_{p,1}^{\frac{3}{p}-1}}^{h;\frac{\beta_0}{\varepsilon}} \right)
		+
		C
		\| (a,u) \|_{\widetilde{L^2}(0,t ; \dB_{p,1}^{\frac{3}{p}})}
		\| (a,u) \|_{\widetilde{L^2}(0,t ; \dB_{q,1}^{\frac{3}{q}})}\\
		&\qquad
		+
		C
		\varepsilon
		\| a \|_{L^{\infty}(0,t ; \dB_{q,1}^{\frac{3}{q}})}
		\| u \|_{L^1(0,t ; \dB_{p,1}^{\frac{3}{p}+1})}^{h;\frac{\beta_0}{\varepsilon}}
		+
		C
		\varepsilon
		\| a \|_{L^{\infty}(0,t ; \dB_{p,1}^{\frac{3}{p}})}
		\| u \|_{L^1(0,t ; \dB_{q,1}^{\frac{3}{q}+1})}^{h;\frac{\beta_0}{\varepsilon}}
	\end{align}
	for all $0 < t< T_{\Omega,\varepsilon}^{\rm max}$,
	where $T_{\Omega,\varepsilon}^{\rm max}$ denotes the maximal existence time.
\end{lemm}
\begin{proof}
	It follows from Lemma \ref{lemm:high-ene} with $f=-\div(au)=-u \cdot \nabla a - (\div u) a$ and $g=-\mathcal{N}_{\varepsilon}[a,u]$ that
	\begin{align}
		&
		\varepsilon \| a \|_{\widetilde{L^{\infty}}( 0,t ; \dB_{p,1}^{\frac{3}{p}})}^{h;\frac{\beta_0}{\varepsilon}}
		+
		\frac{1}{\varepsilon} \| a \|_{L^1( 0,t ; \dB_{p,1}^{\frac{3}{p}})}^{h;\frac{\beta_0}{\varepsilon}}
		+
		\| u \|_{\widetilde{L^{\infty}}( 0,t ; \dB_{p,1}^{\frac{3}{p}-1}) \cap L^1( 0,t ; \dB_{p,1}^{\frac{3}{p}+1})}^{h;\frac{\beta_0}{\varepsilon}}\\
		&\quad
		\leq
		C
		\left( \varepsilon \| a_0 \|_{\dB_{p,1}^{\frac{3}{p}}}^{h;\frac{\beta_0}{\varepsilon}} + \| u_0 \|_{\dB_{p,1}^{\frac{3}{p}-1}}^{h;\frac{\beta_0}{\varepsilon}} \right)
		+
		\frac{C\varepsilon}{p}\left\| \| \div u \|_{L^{\infty}} \| a \|_{\dB_{p,1}^{\frac{d}{p}}} \right\|_{L^1(0,t)} 
		\\
		&\qquad
		+
		C
		\varepsilon
		\sum_{j \in \mathbb{Z}}
		2^{\frac{3}{p}j}
		\| [u\cdot\nabla, \Delta_j]a \|_{L^1(0,t;L^p)}
		+
		C
		\varepsilon \| (\div u) a \|_{L^1(0,t;\dB_{p,1}^{\frac{3}{p}})}\\
		&
		\qquad +
		\left\| \left(\div(au), \mathcal{N}_{\varepsilon}[a,u] \right) \right\|_{L^1(0,t;\dB_{p,1}^{\frac{3}{p}-1})}.
	\end{align}
	From  {\cite{Fu-24}*{Lemma 2.8}}, it follows that
	\begin{align}
		&
		\varepsilon
		\sum_{j \in \mathbb{Z}}
		2^{\frac{3}{p}j}
		\| [u\cdot\nabla, \Delta_j]a \|_{L^1(0,t;L^p)}
		+
		\varepsilon \| (\div u) a \|_{L^1(0,t;\dB_{p,1}^{\frac{3}{p}})}
		\\
		&\quad
		+
		\frac{\varepsilon}{p}\left\| \| \div u \|_{L^{\infty}} \| a \|_{\dB_{p,1}^{\frac{3}{p}}} \right\|_{L^1(0,t)} 
		+
		\left\| \left(\div(au), \mathcal{N}_{\varepsilon}[a,u] \right) \right\|_{L^1(0,t;\dB_{p,1}^{\frac{3}{p}-1})}\\
		&\quad
		\leq
		C
		\varepsilon
		\| a \|_{L^{\infty}(0,t ; \dB_{q,1}^{\frac{3}{q}})}
		\| u \|_{L^1(0,t ; \dB_{p,1}^{\frac{3}{p}+1})}^{h;\frac{\beta_0}{\varepsilon}}
		+
		C
		\varepsilon
		\| a \|_{L^{\infty}(0,t ; \dB_{p,1}^{\frac{3}{p}})}
		\| u \|_{L^1(0,t ; \dB_{q,1}^{\frac{3}{q}+1})}^{h;\frac{\beta_0}{\varepsilon}}\\
		&
		\qquad
		+
		C
		\| (a,u) \|_{\widetilde{L^2}(0,t ; \dB_{p,1}^{\frac{3}{p}})}
		\| (a,u) \|_{\widetilde{L^2}(0,t ; \dB_{q,1}^{\frac{3}{q}})},
	\end{align}
	which {completes} the proof.
\end{proof}

\section{Proof of Theorem \ref{thm:large}}\label{sec:pf-large}
 {Now, we are in a position to present the proof of Theorem \ref{thm:large}.
The proof follows the standard continuation arguments via the global a priori estimates.
To this end, we first define the two norms for the solutions to \eqref{eq:comp-NS}.}
\begin{df}
	Let $\Omega \in \mathbb{R}$, $\alpha,\varepsilon>0$ with $|\Omega|\varepsilon < \alpha < \beta_0/\varepsilon$, and $1 \leq q,r \leq \infty$.
	We define the energy norm $\| (a,u) \|_{ {\mathscr{E}_{\Omega,\varepsilon}}(t)}$ and the auxiliary norm $\| (a,u) \|_{ {\mathscr{A}_{\Omega,\varepsilon;\alpha}^{q,r}}(t)}$ of the solution $(a,u)$ by
	\begin{align}
		&
		\begin{aligned}
			\| (a,u) \|_{ {\mathscr{E}_{\Omega,\varepsilon}}(t)}
			:={}
			&
			\| (a,u) \|_{\widetilde{L^{\infty}}(0,t ; \dB_{2,1}^{-\frac{1}{2}})}^{\ell;\frac{\beta_0}{\varepsilon}}
			+
			\| u \|_{\widetilde{L^2}(0,t ; \dB_{2,1}^{\frac{1}{2}})}^{\ell;\frac{\beta_0}{\varepsilon}}\\
			&
			+
			\| (a,u) \|_{L^{\infty}( 0,t ; \dB_{2,\infty}^{-\frac{3}{2}}) \cap \widetilde{L^1}( 0,t ; \dB_{2,\infty}^{\frac{5}{2}} )}^{\ell;|\Omega|\varepsilon}\\
			&
			+
			\| (a,u) \|_{\widetilde{L^{\infty}}( 0,t ; \dB_{2,1}^{\frac{1}{2}} ) \cap L^1( 0,t ; \dB_{2,1}^{\frac{5}{2}} )}^{m;|\Omega|\varepsilon,\frac{\beta_0}{\varepsilon}}
			+
			\| a \|_{\widetilde{L^2}(0,t ; \dB_{2,1}^{\frac{1}{2}})}^{m;|\Omega|\varepsilon,\frac{\beta_0}{\varepsilon}}\\
			&
			+
			\varepsilon 
			\| a \|_{\widetilde{L^{\infty}}( 0,t ; \dB_{2,1}^{\frac{3}{2}} )}^{h;\frac{\beta_0}{\varepsilon}}
			+
			\frac{1}{\varepsilon} 
			\| a \|_{{L^1}( 0,t ; \dB_{2,1}^{\frac{3}{2}} )}^{h;\frac{\beta_0}{\varepsilon}}
			+
			\| u \|_{\widetilde{L^{\infty}}( 0,t ; \dB_{2,1}^{\frac{1}{2}} ) \cap L^1( 0,t ; \dB_{2,1}^{\frac{5}{2}} )}^{h;\frac{\beta_0}{\varepsilon}},
		\end{aligned}\\
		&
		\begin{aligned}
			\| (a,u) \|_{ {\mathscr{A}_{\Omega,\varepsilon;\alpha}^{q,r}}(t)}
			:={}
			&
            \| (a,u) \|_{
                \widetilde{L^{r'}}( 0,t ; \dB_{q,\infty}^{\frac{3}{q}-3+\frac{4}{r'}} )
            }^{\ell;|\Omega|\varepsilon}
			+
			\| (a,u) \|_{
				\widetilde{L^{r}}( 0,t ; \dB_{q,\infty}^{\frac{3}{q}-3+\frac{4}{r}} )
			}^{\ell;\frac{\beta_0}{\varepsilon}}
			\\
			&
			+
			\| (a,u) \|_{\widetilde{L^r}( 0,t ; \dB_{q,1}^{\frac{3}{q}-1+\frac{2}{r}} )}^{m;|\Omega|\varepsilon,\alpha}
			+
			\| (a,u) \|_{\widetilde{L^{\infty}}( 0,t ; \dB_{q,1}^{\frac{3}{q}-1} ) \cap L^1( 0,t ; \dB_{q,1}^{\frac{3}{q}+1} )}^{m;\alpha,\frac{\beta_0}{\varepsilon}}
			\\
			&
			+
			\varepsilon 
			\| a \|_{\widetilde{L^{\infty}}( 0,t ; \dB_{q,1}^{\frac{3}{q}} )}^{h;\frac{\beta_0}{\varepsilon}}
			+
			\frac{1}{\varepsilon} 
			\| a \|_{{L^1}( 0,t ; \dB_{q,1}^{\frac{3}{q}} )}^{h;\frac{\beta_0}{\varepsilon}}
			+
			\| u \|_{\widetilde{L^{\infty}}( 0,t ; \dB_{q,1}^{\frac{3}{q}-1} ) \cap L^1( 0,t ; \dB_{q,1}^{\frac{3}{q}+1} )}^{h;\frac{\beta_0}{\varepsilon}}
		\end{aligned}
	\end{align}
	for all $t>0$.
\end{df}
{
\begin{rem}
    Let us explain the solution norms defined above.
    The $\mathscr{E}_{\Omega,\varepsilon}(t)$-norm corresponds the energy quantity based on $L^2(\mathbb{R}^3)$.
    The terms 
    \begin{gather}
        \| (a,u) \|_{\widetilde{L^{\infty}}( 0,t ; \dB_{2,1}^{\frac{1}{2}} ) \cap L^1( 0,t ; \dB_{2,1}^{\frac{5}{2}} )}^{m;|\Omega|\varepsilon,\frac{\beta_0}{\varepsilon}}
        ,\quad
        \varepsilon 
        \| a \|_{\widetilde{L^{\infty}}( 0,t ; \dB_{2,1}^{\frac{3}{2}} )}^{h;\frac{\beta_0}{\varepsilon}}
        ,\\
        \frac{1}{\varepsilon} 
        \| a \|_{L^1( 0,t ; \dB_{2,1}^{\frac{3}{2}} )}^{h;\frac{\beta_0}{\varepsilon}}
        ,
        \quad
        \| u \|_{\widetilde{L^{\infty}}( 0,t ; \dB_{2,1}^{\frac{1}{2}} ) \cap L^1( 0,t ; \dB_{2,1}^{\frac{5}{2}} )}^{h;\frac{\beta_0}{\varepsilon}}
    \end{gather}
    are components of well-known solution {norms} for the compressible Navier--Stokes system in critical Besov spaces, while the terms
    \begin{gather}
        \| (a,u) \|_{\widetilde{L^{\infty}}(0,t ; \dB_{2,1}^{-\frac{1}{2}})}^{\ell;\frac{\beta_0}{\varepsilon}},
        \quad
        \| u \|_{\widetilde{L^2}(0,t ; \dB_{2,1}^{\frac{1}{2}})}^{\ell;\frac{\beta_0}{\varepsilon}},
        \\
        \| (a,u) \|_{L^{\infty}( 0,t ; \dB_{2,\infty}^{-\frac{3}{2}}) \cap \widetilde{L^1}( 0,t ; \dB_{2,\infty}^{\frac{5}{2}} )}^{\ell;|\Omega|\varepsilon},
        \quad
        \| a \|_{\widetilde{L^2}(0,t ; \dB_{2,1}^{\frac{1}{2}})}^{m;|\Omega|\varepsilon,\frac{\beta_0}{\varepsilon}}
    \end{gather}
    have lower regularities than the critical spaces, which are to be used for controlling the $4$th order dissipative semigroup $\{ e^{-t\Delta^2} \}_{t>0}$ {appearing} in the linear structure in the low{-}frequency region; see \eqref{ene-low}.
    As the $\mathscr{E}_{\Omega,\varepsilon}(t)$-norm may not be small for large initial data, 
    we follow the idea of \cites{Fu-24,Fuj-Wat-25} to introduce the auxiliary $\mathscr{A}_{\Omega,\varepsilon;\alpha}^{q,r}(t)$-norm based on $L^q(\mathbb{R}^3)$ which corresponds to the Strichartz norm and will {be small} for $|\Omega| \gg 1$ and $0 <\varepsilon \ll 1$.
    Note that $\mathscr{A}_{\Omega,\varepsilon;\alpha}^{q,r}(t)$-norm contains the non-critical regularities
    \begin{align}
        \| (a,u) \|_{
            \widetilde{L^{r'}}( 0,t ; \dB_{q,\infty}^{\frac{3}{q}-3+\frac{4}{r'}} )
        }^{\ell;|\Omega|\varepsilon},
        \quad
        \| (a,u) \|_{
            \widetilde{L^{r}}( 0,t ; \dB_{q,\infty}^{\frac{3}{q}-3+\frac{4}{r}} )
        }^{\ell;\frac{\beta_0}{\varepsilon}},
    \end{align}
    due to the similar circumstance in $\mathscr{E}_{\Omega,\varepsilon}(t)$-norm.
\end{rem}
In the following two lemmas, we provide relations between the above two solution norms and several Chemin--Lerner norms, which are used in {establishing} the global a priori estimates of the solutions.}
\begin{lemm}\label{lemm:E}
	Let $\Omega \in \mathbb{R}$ and $0<\varepsilon \leq 1$ satisfy
	$|\Omega| \varepsilon \leq 1$ and $|\Omega|\varepsilon  < \beta_0/\varepsilon$.
	Then, there exists an absolute positive constant $C$ such that
	\begin{align}
		\varepsilon \| a \|_{\widetilde{L^{\infty}}( 0,t ; \dB_{2,1}^{\frac{3}{2}})}
		&
		+
		\| u \|_{L^1( 0,t ; \dB_{2,1}^{\frac{5}{2}} )}^{h;\frac{\beta_0}{\varepsilon}}
		+
		\| (a,u) \|_{\widetilde{L^2}( 0,t ; \dB_{2,1}^{\frac{3}{2}})}
		+
		\| (a,u) \|_{\widetilde{L^{\infty}}(0,t ; \dB_{2,1}^{\frac{1}{2}})}\\
		&
		+
		\| (a,u) \|_{\widetilde{L^2}( 0,t ; \dB_{2,1}^{\frac{1}{2}})}
		+
		\| (a,u) \|_{\widetilde{L^{\infty}}(0,t ; \dB_{2,1}^{-\frac{1}{2}})}
		\leq
		C
		\| (a,u) \|_{ {\mathscr{E}_{\Omega,\varepsilon}}(t)}.
	\end{align}
\end{lemm}
 {
We omit the proof of the above lemma since it is easily obtained.
}
\begin{lemm}\label{lemm:AE} 
	Let $\Omega \in \mathbb{R}$, $0<\varepsilon \leq 1$ and $\alpha>0$ satisfy
	$|\Omega| \varepsilon \leq 1$ and $|\Omega|\varepsilon < \alpha < \beta_0/\varepsilon$.
	Let $2 \leq q \leq \infty$ and $2 < r < \infty$.
	Then, there exists a positive constant $C=C(q,r)$ such that
	\begin{align}
		&
		\begin{aligned}\label{AE-1}
			\varepsilon 
			\| a \|_{\widetilde{L^{\infty}}( 0,t ; \dB_{q,1}^{\frac{3}{q}})}
			\leq
			C
			\varepsilon \alpha
			\| (a,u) \|_{ {\mathscr{E}_{\Omega,\varepsilon}}(t)}
			+
			C
			\| (a,u) \|_{ {\mathscr{A}_{\Omega,\varepsilon;\alpha}^{q,r}}(t)},
		\end{aligned}\\
		&
		\begin{aligned}\label{AE-2}
			\| (a,u) \|_{\widetilde{L^{r}}( 0,t ; \dB_{q,1}^{\frac{3}{q}-1+\frac{2}{r}})}
			\leq
			C
			\| (a,u) \|_{ {\mathscr{A}_{\Omega,\varepsilon;\alpha}^{q,r}}(t)},
		\end{aligned}\\
		&
			\begin{aligned}\label{AE-3}
				&
				\| (a,u) \|_{\widetilde{L^2}( 0,t ; \dB_{q,1}^{\frac{3}{q}})}
				+
				\| (a,u) \|_{\widetilde{L^2}( 0,t ; \dB_{q, \infty}^{\frac{3}{q}-1})}
				\\
				&\quad
				\leq{}
				C
				\| (a,u) \|_{ {\mathscr{A}_{\Omega,\varepsilon;\alpha}^{q,r}}(t)}
				+
				C
				\| (a,u) \|_{ {\mathscr{A}_{\Omega,\varepsilon;\alpha}^{q,r}}(t)}^{\frac{r}{2(r-1)}}
				\| (a,u) \|_{ {\mathscr{E}_{\Omega,\varepsilon}}(t)}^{\frac{r-2}{2(r-1)}},
			\end{aligned}
		\\
		&
		\begin{aligned}\label{AE-4}
			\begin{split}
				&
				\| a \|_{\widetilde{L^{r'}}(0,t;\dB_{q,1}^{\frac{3}{q}-1+\frac{2}{r'}})}^{\ell;\frac{4\beta_0}{\varepsilon}}
				+
				\| u \|_{\widetilde{L^{r'}}(0,t;\dB_{q,1}^{\frac{3}{q}-1+\frac{2}{r'}})}
				\leq{}
					C
					\| (a,u) \|_{ {\mathscr{A}_{\Omega,\varepsilon;\alpha}^{q,r}}(t)}.
			\end{split}
		\end{aligned}
	\end{align}
\end{lemm}
\begin{proof}
	We first show \eqref{AE-1}.
	For the high{-}frequency part of $\varepsilon \| a \|_{\widetilde{L^{\infty}}( 0,t ; \dB_{q,1}^{\frac{3}{q}})}$, we  {have} 
	\begin{align}
		\varepsilon 
		\| a \|_{\widetilde{L^{\infty}}( 0,t ; \dB_{q,1}^{\frac{3}{q}})}^{h;\frac{\beta_0}{\varepsilon}}
		&
		\leq{}
		\| (a,u) \|_{ {\mathscr{A}_{\Omega,\varepsilon;\alpha}^{q,r}}(t)}.
	\end{align}
	For the low{-}frequency part of $\varepsilon \| a \|_{\widetilde{L^{\infty}}( 0,t ; \dB_{q,1}^{\frac{3}{q}})}$, the Bernstein inequality  {implies} that
	\begin{align}
		\varepsilon
		\| a \|_{\widetilde{L^{\infty}}( 0,t ; \dB_{q,1}^{\frac{3}{q}})}^{\ell;\frac{\beta_0}{\varepsilon}}
		\leq{}&
		C\varepsilon
		\| a \|_{\widetilde{L^{\infty}}( 0,t ; \dB_{2,1}^{\frac{3}{2}})}^{\ell;|\Omega|\varepsilon}
		+
		C\varepsilon 
		\| a \|_{\widetilde{L^{\infty}}( 0,t ; \dB_{2,1}^{\frac{3}{2}})}^{m;|\Omega|\varepsilon,\alpha}
		+
		C
		\| a \|_{\widetilde{L^{\infty}}( 0,t ; \dB_{q,1}^{\frac{3}{q}-1})}^{m;\alpha,\frac{\beta_0}{\varepsilon}}
		\\
		\leq{}&
		C\varepsilon
		(|\Omega|\varepsilon)^{3}
		\| a \|_{\widetilde{L^{\infty}}( 0,t ; \dB_{2,\infty}^{-\frac{3}{2}})}^{\ell;|\Omega|\varepsilon}
		+
		C\varepsilon \alpha
		\| a \|_{\widetilde{L^{\infty}}( 0,t ; \dB_{2,1}^{\frac{1}{2}})}^{m;|\Omega|\varepsilon,\alpha}
		+
		C
		\| a \|_{\widetilde{L^{\infty}}( 0,t ; \dB_{q,1}^{\frac{3}{q}-1})}^{m;\alpha,\frac{\beta_0}{\varepsilon}}
		\\
		\leq{}&
		C\varepsilon \alpha \| (a,u) \|_{ {\mathscr{E}_{\Omega,\varepsilon}}(t)}
		+
		C\| (a,u) \|_{ {\mathscr{A}_{\Omega,\varepsilon;\alpha}^{q,r}}(t)}.
	\end{align}
	Combining the above two estimates, we obtain \eqref{AE-1}.
	For the proof of \eqref{AE-2}, it follows from the Bernstein inequality and the interpolation inequality that 
	\begin{align}
		&
		\| (a,u) \|_{\widetilde{L^{r}}( 0,t ; \dB_{q,1}^{\frac{3}{q}-1+\frac{2}{r}})}\\
		&\quad\leq{}
		\| (a,u) \|_{\widetilde{L^{r}}( 0,t ; \dB_{q,1}^{\frac{3}{q}-1+\frac{2}{r}})}^{\ell;|\Omega|\varepsilon}
		+
		\| (a,u) \|_{\widetilde{L^{r}}( 0,t ; \dB_{q,1}^{\frac{3}{q}-1+\frac{2}{r}})}^{m;|\Omega|\varepsilon,\frac{\beta_0}{\varepsilon}}
		\\
		&\qquad
		+
		\| a \|_{\widetilde{L^{r}}( 0,t ; \dB_{q,1}^{\frac{3}{q}-1+\frac{2}{r}})}^{h;\frac{\beta_0}{\varepsilon}}
		+
		\| u \|_{\widetilde{L^{r}}( 0,t ; \dB_{q,1}^{\frac{3}{q}-1+\frac{2}{r}})}^{h;\frac{\beta_0}{\varepsilon}}\\
		&\quad 
		\leq{}
		(|\Omega|\varepsilon)^{2-\frac{2}{r}}
		\| (a,u) \|_{\widetilde{L^{r}}( 0,t ; \dB_{q,1}^{\frac{3}{q}-3+\frac{4}{r}})}^{\ell;|\Omega|\varepsilon}
		+
		\| (a,u) \|_{\widetilde{L^{r}}( 0,t ; \dB_{q,1}^{\frac{3}{q}-1+\frac{2}{r}})}^{m;|\Omega|\varepsilon,\frac{\beta_0}{\varepsilon}}
		\\
		&\qquad
		+
		\left( \frac{\beta_0}{\varepsilon} \right)^{-1+\frac{2}{r}}
		\| a \|_{\widetilde{L^{r}}( 0,t ; \dB_{q,1}^{\frac{3}{q}})}^{h;\frac{\beta_0}{\varepsilon}}
		+
		\| u \|_{\widetilde{L^{r}}( 0,t ; \dB_{q,1}^{\frac{3}{q}-1+\frac{2}{r}})}^{h;\frac{\beta_0}{\varepsilon}}\\
		&\quad
		\leq{}
		\| (a,u) \|_{\widetilde{L^{r}}( 0,t ; \dB_{q,1}^{\frac{3}{q}-3+\frac{4}{r}})}^{\ell;|\Omega|\varepsilon}
		+
		\| (a,u) \|_{\widetilde{L^{r}}( 0,t ; \dB_{q,1}^{\frac{3}{q}-1+\frac{2}{r}})}^{m;|\Omega|\varepsilon,\frac{\beta_0}{\varepsilon}}
		\\
		&\qquad 
		+
		C\varepsilon
		\| a \|_{\widetilde{L^{\infty}}( 0,t ; \dB_{q,1}^{\frac{3}{q}})}^{h;\frac{\beta_0}{\varepsilon}}
		+
		\frac{C}{\varepsilon}
		\| a \|_{L^1( 0,t ; \dB_{q,1}^{\frac{3}{q}})}^{h;\frac{\beta_0}{\varepsilon}}
		+
		\| u \|_{\widetilde{L^{\infty}}( 0,t ; \dB_{q,1}^{\frac{3}{q}-1}) \cap L^1( 0,t ; \dB_{q,1}^{\frac{3}{q}+1})}^{h;\frac{\beta_0}{\varepsilon}}\\
		&
		\quad 
		\leq{}
		C\| (a,u) \|_{ {\mathscr{A}_{\Omega,\varepsilon;\alpha}^{q,r}}(t)}.
	\end{align}
	Next, we prove \eqref{AE-3}.
	By the interpolation inequality, the Bernstein inequality, and \eqref{AE-2}, it holds
		\begin{align}
			\label{pf:AE-3-1}
			\notag
			\| (a,u) \|_{\widetilde{L^2}( 0,t ;\dB_{q,1}^{\frac{3}{q}} )}
			\leq{}&
			\| (a,u) \|_{\widetilde{L^2}( 0,t ;\dB_{q,1}^{\frac{3}{q}} )}^{\ell;|\Omega|\varepsilon}
			+
			\| (a,u) \|_{\widetilde{L^2}( 0,t ;\dB_{q,1}^{\frac{3}{q}} )}^{m;|\Omega|\varepsilon,\frac{\beta_0}{\varepsilon}}
			+
			\| (a,u) \|_{\widetilde{L^2}( 0,t ;\dB_{q,1}^{\frac{3}{q}} )}^{h;\frac{\beta_0}{\varepsilon}}\\
			\notag
			\leq{}&
			C|\Omega|\varepsilon
			\| (a,u) \|_{\widetilde{L^2}( 0,t ;\dB_{q,\infty}^{\frac{3}{q}-1} )}^{\ell;|\Omega|\varepsilon}
			\\
			\notag
			&
			+
			\sp{\| (a,u) \|_{\widetilde{L^{r}}( 0,t ;\dB_{q,1}^{\frac{3}{q}-1+\frac{2}{r}} )}^{h;|\Omega|\varepsilon,\frac{\beta_0}{\varepsilon}}}^{\frac{r}{2(r-1)}}
			\sp{\| (a,u) \|_{\widetilde{L^1}( 0,t ;\dB_{q,1}^{\frac{3}{q}+1} )}^{h;|\Omega|\varepsilon,\frac{\beta_0}{\varepsilon}}}^{\frac{r-2}{2(r-1)}}\\
			\leq{}&
			C\| (a,u) \|_{ {\mathscr{A}_{\Omega,\varepsilon;\alpha}^{q,r}}(t)}
			+
			C
			\| (a,u) \|_{ {\mathscr{A}_{\Omega,\varepsilon;\alpha}^{q,r}}(t)}^{\frac{r}{2(r-1)}}
			\| (a,u) \|_{ {\mathscr{E}_{\Omega,\varepsilon}}(t)}^{\frac{r-2}{2(r-1)}}.
		\end{align}
	Using \eqref{pf:AE-3-1} and a similar argument as above, we have
	\begin{align}
		\| (a,u) \|_{\widetilde{L^2}( 0,t ;\dB_{q, \infty}^{\frac{3}{q}-1} )}
		\leq{}&
		\| (a,u) \|_{\widetilde{L^2}( 0,t ;\dB_{q, \infty}^{\frac{3}{q}-1} )}^{\ell;\frac{\beta_0}{\varepsilon}}
		+
		\| (a,u) \|_{\widetilde{L^2}( 0,t ;\dB_{q, \infty}^{\frac{3}{q}-1} )}^{h;\frac{\beta_0}{\varepsilon}}\\
		\leq{}&
		C
		\left(
		\| (a,u) \|_{\widetilde{L^{r}}( 0,t ;\dB_{q,\infty}^{\frac{3}{q}-3+\frac{4}{r}} )}^{\ell;\frac{\beta_0}{\varepsilon}}
		\right)^{\frac{r}{2(r-1)}}
		\| (a,u) \|_{\widetilde{L^1}( 0,t ;\dB_{q,\infty}^{\frac{3}{q}+1} )}^{\frac{r-2}{2(r-1)}}\\
		&
		+
		\frac{\varepsilon}{\beta_0}
		\| (a,u) \|_{\widetilde{L^2}( 0,t ;\dB_{q,1}^{\frac{3}{q}} )}\\
		\leq{}&
		C\| (a,u) \|_{ {\mathscr{A}_{\Omega,\varepsilon;\alpha}^{q,r}}(t)}
			+
			C
			\| (a,u) \|_{ {\mathscr{A}_{\Omega,\varepsilon;\alpha}^{q,r}}(t)}^{\frac{r}{2(r-1)}}
			\| (a,u) \|_{ {\mathscr{E}_{\Omega,\varepsilon}}(t)}^{\frac{r-2}{2(r-1)}}.
	\end{align}
	Thus, we obtain \eqref{AE-3}.
	Finally, we show \eqref{AE-4}.
	From \eqref{pf:AE-3-1}, the Bernstein inequality, and the interpolation inequality, it follows that
		\begin{align}
			\| a \|_{\widetilde{L^{r'}}(0,t;\dB_{q,1}^{\frac{3}{q}-1+\frac{2}{r'}})}^{\ell;\frac{4\beta_0}{\varepsilon}}
			\leq{}&
			\| a \|_{\widetilde{L^{r'}}(0,t;\dB_{q,1}^{\frac{3}{q}-1+\frac{2}{r'}})}^{\ell;|\Omega|\varepsilon}
			+
			\| a \|_{\widetilde{L^{r'}}(0,t;\dB_{q,1}^{\frac{3}{q}-1+\frac{2}{r'}})}^{m;|\Omega|\varepsilon,\frac{\beta_0}{\varepsilon}}
			+
			\| a \|_{\widetilde{L^{r'}}(0,t;\dB_{q,1}^{\frac{3}{q}-1+\frac{2}{r'}})}^{m;\frac{\beta_0}{\varepsilon},\frac{4\beta_0}{\varepsilon}}\\
			\leq{}&
			(|\Omega|\varepsilon)^{\frac{2}{r}}
			\| a \|_{\widetilde{L^{r'}}(0,t;\dB_{q,\infty}^{\frac{3}{q}-3+\frac{4}{r'}})}^{\ell;|\Omega|\varepsilon}
			\\
			&
			+
			\| a \|_{L^{\infty}(0,t;\dB_{q,1}^{\frac{3}{q}-1}) \cap L^1(0,t;\dB_{q,1}^{\frac{3}{q}+1})}^{m;|\Omega|\varepsilon,\frac{\beta_0}{\varepsilon}}\\
			&
			+
			C\varepsilon
			\| a \|_{\widetilde{L^{\infty}}( 0,t ; \dB_{q,1}^{\frac{3}{q}})}^{h;\frac{\beta_0}{\varepsilon}}
			+
			\frac{C}{\varepsilon}
			\| a \|_{L^1( 0,t ; \dB_{q,1}^{\frac{3}{q}})}^{h;\frac{\beta_0}{\varepsilon}}\\
			\leq{}&
			C\| (a,u) \|_{ {\mathscr{A}_{\Omega,\varepsilon;\alpha}^{q,r}}(t)}.
		\end{align}
	Similar observation yields
		\begin{align}
			\| u \|_{\widetilde{L^{r'}}(0,t;\dB_{q,1}^{\frac{3}{q}-1+\frac{2}{r'}})}^{\ell;\frac{4\beta_0}{\varepsilon}}
			\leq{}&
			\| u \|_{\widetilde{L^{r'}}(0,t;\dB_{q,1}^{\frac{3}{q}-1+\frac{2}{r'}})}^{\ell;|\Omega|\varepsilon}
			+
			\| u \|_{\widetilde{L^{r'}}(0,t;\dB_{q,1}^{\frac{3}{q}-1+\frac{2}{r'}})}^{h;|\Omega|\varepsilon}
			\\
			\leq{}&
			(|\Omega|\varepsilon)^{\frac{2}{r}}
			\| u \|_{\widetilde{L^{r'}}(0,t;\dB_{q,\infty}^{\frac{3}{q}-3+\frac{4}{r'}})}^{\ell;|\Omega|\varepsilon}
			+
			\| u \|_{L^{\infty}(0,t;\dB_{q,1}^{\frac{3}{q}-1}) \cap L^1(0,t;\dB_{q,1}^{\frac{3}{q}+1})}^{h;|\Omega|\varepsilon}
			\\
			\leq{}&
			C\| (a,u) \|_{ {\mathscr{A}_{\Omega,\varepsilon;\alpha}^{q,r}}(t)}.
		\end{align}
	Hence, we obtain \eqref{AE-4} and complete the proof.
\end{proof}
Now, we are in a position to present the proof of Theorem \ref{thm:large}.
\begin{proof}[Proof of Theorem \ref{thm:large}]
	Let $q$ and $r$ satisfy
	\begin{align}
		2 < q < 3, \quad 
		2 < r < \infty ,\quad 
		\frac{1}{q} + \frac{1}{r} \leq \frac{1}{2}, \quad
		\frac{2}{r} \leq \frac{3}{q} - \frac{1}{2}, \quad 
		-\frac{3}{q} \leq -\frac{3}{2} + \frac{4}{r} < \frac{3}{q}.
	\end{align}
	Note that these $q$ and $r$ satisfy all assumptions of lemmas in Section \ref{sec:a-priori} and Lemma \ref{lemm:AE}.
	Let $\delta$ be a positive constant to be determined later.
	Since 
	\begin{align}
		(a_0,u_0) 
		&
		\in (\dB_{2,\infty}^{-\frac{3}{2}}(\mathbb{R}^3) \cap \dB_{2,1}^{\frac{3}{2}}(\mathbb{R}^3))  \times (\dB_{2,\infty}^{-\frac{3}{2}}(\mathbb{R}^3) \cap \dB_{2,1}^{\frac{1}{2}}(\mathbb{R}^3))^3\\
		&
		\subset (\dB_{2,1}^{\frac{1}{2}}(\mathbb{R}^3) \cap \dB_{2,1}^{\frac{3}{2}}(\mathbb{R}^3))  \times \dB_{2,1}^{\frac{1}{2}}(\mathbb{R}^3)^3,
	\end{align}
	there exists a constant $\alpha_{\delta} = \alpha_{\delta}(a_0,u_0) \geqslant 1 $ such that
	\begin{align}
		\| (a_0,u_0) \|_{\dB_{2,1}^{\frac{1}{2}}}^{h;\alpha_{\delta}}
		+
		\| a_0 \|_{\dB_{2,1}^{\frac{3}{2}}}^{h;\alpha_{\delta}}
		\leq \delta.
	\end{align}
	For this $\alpha_{\delta}$,
	there exist $\Omega_{\delta}=\Omega_{\delta}(a_0,u_0) \geqslant 1$ and $c_{\delta} = c_{\delta}(a_0,u_0)>0$ such that 
	\begin{gather}\label{choice}
		|\Omega| \varepsilon \leq 1, \quad 
		\frac{\beta_0}{\varepsilon} > \alpha_{\delta}, \quad 
		\alpha_{\delta}^{\frac{2}{r}}|\Omega|^{-\frac{1}{r}} \leq 1,\\
		(|\Omega| \varepsilon)^{\frac{2}{r}}\mathscr{D}_{\varepsilon}^*[a_0,u_0] \leq \delta, \quad
		\alpha_{\delta}^{\frac{2}{r}}|\Omega|^{-\frac{1}{r}}
		\| (a_0,u_0) \|_{\dB_{2,\infty}^{-\frac{3}{2}} \cap \dB_{2,1}^{\frac{1}{2}}}
		\leq \delta,\quad
		\varepsilon
		\alpha_{\delta}
		\mathscr{D}_{\varepsilon}[a_0,u_0]
		\leq 
		\delta
	\end{gather}
	for all $\Omega \in \mathbb{R}\setminus \{ 0\}$ and $0 < \varepsilon \leq 1$ with $\Omega_{\delta} \leq |\Omega| \leq c_{\delta}/\varepsilon$.
	Here, we have set 
	\begin{align}
		\mathscr{D}_{\varepsilon}[a_0,u_0]
		:=
		\mathscr{D}^*_{\varepsilon}[a_0,u_0]
		+
		\| (a_0,u_0) \|_{\dB_{2,1}^{\frac{1}{2}}}
		+
		\varepsilon \| a_0 \|_{\dB_{2,1}^{\frac{3}{2}}}
		.
	\end{align}
	Let $(a,u)$ be the local solution on $[0,T_{\Omega,\varepsilon}^{\rm max})$ constructed in  {Lemma \ref{lemm:LWP}}, where $T_{\Omega,\varepsilon}^{\rm max}$ denotes the maximal existence time.

	\vskip1pc
	\noindent
	{\it Step.1{:} A priori estimates for the energy norm.}
	It follows from Lemmas \ref{lemm:low-ene}, \ref{lemm:E}, and \ref{lemm:AE} that
	\begin{align}\label{pf:main-E-1}
		\begin{split}
			&\| (a,u) \|_{L^{\infty}( 0,t ;\dB_{2,\infty}^{-\frac{3}{2}}) \cap \widetilde{L^1}( 0,t ; \dB_{2,\infty}^{\frac{5}{2}})}^{\ell;|\Omega|\varepsilon}
			\leq{}
			C
			\mathscr{D}_{\varepsilon}^*[a_0,u_0]
			+
			\\
			&\quad 
			+
			C\left( 
			1 
			+
			\varepsilon \alpha_{\delta}
			\| (a,u) \|_{ {\mathscr{E}_{\Omega,\varepsilon}}(t)}
			+
			\| (a,u) \|_{ {\mathscr{A}_{\Omega,\varepsilon;\alpha_\delta}^{q,r}}(t)}
			\right)\\
			&\quad \qquad 
			\times 
			\left(
			\| (a,u) \|_{ {\mathscr{A}_{\Omega,\varepsilon;\alpha_\delta}^{q,r}}(t)}^2
			+
			\| (a,u) \|_{ {\mathscr{A}_{\Omega,\varepsilon;\alpha_\delta}^{q,r}}(t)}^{\frac{r}{r-1}}
			\| (a,u) \|_{ {\mathscr{E}_{\Omega,\varepsilon}}(t)}^{\frac{r-2}{r-1}}
			\right)\\
			&\quad
			+
			C
			\varepsilon
			\| (a,u) \|_{ {\mathscr{E}_{\Omega,\varepsilon}}(t)}^2.
		\end{split}
	\end{align}
	By Lemmas \ref{lemm:middle-ene-2} and \ref{lemm:AE}, we have
	\begin{align}\label{pf:main-E-2}
		\begin{split}
			&
			\| (a,u) \|_{\widetilde{L^{\infty}}(0,t; \dB_{2,1}^{-\frac{1}{2}})}^{\ell;\frac{\beta_0}{\varepsilon}}
			+
			\| a \|_{\widetilde{L^2}(0,t;\dB_{2,1}^{\frac{1}{2}})}^{m;|\Omega|\varepsilon,\frac{\beta_0}{\varepsilon}}
			+
			\| u \|_{\widetilde{L^2}(0,t;\dB_{2,1}^{\frac{1}{2}})}^{\ell;\frac{\beta_0}{\varepsilon}}\\
			&\quad
			\leq
			C
			\mathscr{D}_{\varepsilon}[a_0,u_0]
			+
			C
			\| (a,u) \|_{ {\mathscr{A}_{\Omega,\varepsilon;\alpha_\delta}^{q,r}}(t)}^2
			+
			C
			\| (a,u) \|_{ {\mathscr{A}_{\Omega,\varepsilon;\alpha_\delta}^{q,r}}(t)}^{\frac{r}{r-1}}
			\| (a,u) \|_{ {\mathscr{E}_{\Omega,\varepsilon}}(t)}^{\frac{r-2}{r-1}}.
		\end{split}
	\end{align}
	We see by Lemmas \ref{lemm:middle-ene-1} and \ref{lemm:AE} that
	\begin{align}\label{pf:main-E-3}
		\begin{split}
			\| (a,u) \|_{\widetilde{L^{\infty}}( 0,t ; \dB_{2,1}^{\frac{1}{2}} ) \cap L^1( 0,t ; \dB_{2,1}^{\frac{5}{2}} )}^{m;|\Omega|\varepsilon,\frac{\beta_0}{\varepsilon}} 
			\leq{}&
			C
			\mathscr{D}_{\varepsilon}[a_0,u_0]
			+
			C
			\| (a,u) \|_{ {\mathscr{A}_{\Omega,\varepsilon;\alpha_\delta}^{q,r}}(t)}^{\frac{r}{r-1}}
			\| (a,u) \|_{ {\mathscr{E}_{\Omega,\varepsilon}}(t)}^{\frac{r-2}{r-1}} \quad \\
			&
			+
			C
			\left(
			\varepsilon \alpha_{\delta}
			\| (a,u) \|_{ {\mathscr{E}_{\Omega,\varepsilon}}(t)}
			+
			\| (a,u) \|_{ {\mathscr{A}_{\Omega,\varepsilon;\alpha_\delta}^{q,r}}(t)}
			\right)^2.
		\end{split}
	\end{align}
	From Lemmas \ref{lemm:a-propri-high} {(}with $p=2${)}, \ref{lemm:E}, and \ref{lemm:AE}, there holds
	\begin{align}\label{pf:main-E-4}
		\begin{split}
			&
			\varepsilon \| a \|_{\widetilde{L^{\infty}}( 0,t ; \dB_{2,1}^{\frac{3}{2}})}^{h;\frac{\beta_0}{\varepsilon}}
			+
			\frac{1}{\varepsilon} \| a \|_{L^1( 0,t ; \dB_{2,1}^{\frac{3}{2}})}^{h;\frac{\beta_0}{\varepsilon}}
			+
			\| u \|_{\widetilde{L^{\infty}}( 0,t ; \dB_{2,1}^{\frac{1}{2}}) \cap L^1( 0,t ; \dB_{p,1}^{\frac{5}{2}})}^{h;\frac{\beta_0}{\varepsilon}}\\
			&\quad
			\leq
			C
			\mathscr{D}_{\varepsilon}[a_0,u_0]
			+
			C
			\| (a,u) \|_{ {\mathscr{A}_{\Omega,\varepsilon;\alpha_\delta}^{q,r}}(t)}^{\frac{r}{2(r-1)}}
			\| (a,u) \|_{ {\mathscr{E}_{\Omega,\varepsilon}}(t)}^{\frac{r-2}{2(r-1)}+1}\\
			&\qquad
			+
			C
			\left(
			\varepsilon \alpha_{\delta}
			\| (a,u) \|_{ {\mathscr{E}_{\Omega,\varepsilon}}(t)}
			+
			\| (a,u) \|_{ {\mathscr{A}_{\Omega,\varepsilon;\alpha_\delta}^{q,r}}(t)}
			\right)
			\| (a,u) \|_{ {\mathscr{E}_{\Omega,\varepsilon}}(t)}.
		\end{split}
	\end{align}
	Hence, collecting \eqref{pf:main-E-1}, \eqref{pf:main-E-2}, \eqref{pf:main-E-3}, and \eqref{pf:main-E-4}, we obtain 
	\begin{align}\label{E}
		\begin{split}
			&\| (a,u) \|_{ {\mathscr{E}_{\Omega,\varepsilon}}(t)}
			\leq 
			C_1
			\mathscr{D}_{\varepsilon}[a_0,u_0]+
			\\
			&\quad 
			+
			C_1\left( 
			1 
			+
			\varepsilon \alpha_{\delta}
			\| (a,u) \|_{ {\mathscr{E}_{\Omega,\varepsilon}}(t)}
			+
			\| (a,u) \|_{ {\mathscr{A}_{\Omega,\varepsilon;\alpha_\delta}^{q,r}}(t)}
			\right)\\
			&\quad \qquad 
			\times 
			\left(
			\| (a,u) \|_{ {\mathscr{A}_{\Omega,\varepsilon;\alpha_\delta}^{q,r}}(t)}^2
			+
			\| (a,u) \|_{ {\mathscr{A}_{\Omega,\varepsilon;\alpha_\delta}^{q,r}}(t)}^{\frac{r}{r-1}}
			\| (a,u) \|_{ {\mathscr{E}_{\Omega,\varepsilon}}(t)}^{\frac{r-2}{r-1}}
			\right)\\
			&\quad 
			+
			C_1
			\left(
			\varepsilon \alpha_{\delta}
			\| (a,u) \|_{ {\mathscr{E}_{\Omega,\varepsilon}}(t)}
			+
			\| (a,u) \|_{ {\mathscr{A}_{\Omega,\varepsilon;\alpha_\delta}^{q,r}}(t)}
			\right)^2\\
			&\quad+
			C_1
			\left(
			\varepsilon \alpha_{\delta}
			\| (a,u) \|_{ {\mathscr{E}_{\Omega,\varepsilon}}(t)}
			+
			\| (a,u) \|_{ {\mathscr{A}_{\Omega,\varepsilon;\alpha_\delta}^{q,r}}(t)}
			\right)
			\| (a,u) \|_{ {\mathscr{E}_{\Omega,\varepsilon}}(t)}\\
			&\quad
			+
			C_1
			\varepsilon
			\| (a,u) \|_{ {\mathscr{E}_{\Omega,\varepsilon}}(t)}^2
			+
			C_1
			\| (a,u) \|_{ {\mathscr{A}_{\Omega,\varepsilon;\alpha_\delta}^{q,r}}(t)}^{\frac{r}{2(r-1)}}
			\| (a,u) \|_{ {\mathscr{E}_{\Omega,\varepsilon}}(t)}^{\frac{r-2}{2(r-1)}+1}.
		\end{split}
	\end{align}
	for some positive constant $C_1=C_1(\mu,P,q,r)$.
	\vskip1pc
	\noindent
	{\it Step.2{:} A priori estimates for the auxiliary norm.}
	By Lemmas \ref{lemm:low-ene} and \ref{lemm:AE}, it holds
	\begin{align}\label{pf:main-A-1}
		\begin{split}
			&
			\| (a,u) \|_{\widetilde{L^{r}}( 0,t ; \dB_{q,\infty}^{\frac{3}{q}-3+\frac{4}{r}} )}^{\ell;|\Omega|\varepsilon}
			\leq 
			C\delta+\\
			&\quad 
			+
			C\left( 
			1 
			+
			\varepsilon \alpha_{\delta}
			\| (a,u) \|_{ {\mathscr{E}_{\Omega,\varepsilon}}(t)}
			+
			\| (a,u) \|_{ {\mathscr{A}_{\Omega,\varepsilon;\alpha_\delta}^{q,r}}(t)}
			\right)\\
			&\quad \qquad 
			\times 
			\left(
			\| (a,u) \|_{ {\mathscr{A}_{\Omega,\varepsilon;\alpha_\delta}^{q,r}}(t)}^2
			+
			\| (a,u) \|_{ {\mathscr{A}_{\Omega,\varepsilon;\alpha_\delta}^{q,r}}(t)}^{\frac{r}{r-1}}
			\| (a,u) \|_{ {\mathscr{E}_{\Omega,\varepsilon}}(t)}^{\frac{r-2}{r-1}}
			\right)\\
			&\quad
			+
			C
			\varepsilon
			\| (a,u) \|_{ {\mathscr{E}_{\Omega,\varepsilon}}(t)}
			\| (a,u) \|_{ {\mathscr{A}_{\Omega,\varepsilon;\alpha_\delta}^{q,r}}(t)}.
		\end{split}
	\end{align}
	From Lemmas \ref{lemm:middle-str-2}, \ref{lemm:E}, and \ref{lemm:AE}, it follows that
	\begin{align}\label{pf:main-A-2}
		\begin{split}
			&
			\| (a,u) \|_{\widetilde{L^{r}}( 0,t ; \dB_{q,\infty}^{\frac{3}{q}-3+\frac{4}{r}} )}^{m;|\Omega|\varepsilon,\frac{\beta_0}{\varepsilon}} \\
			& \leq {}
			C\delta
			+
			C|\Omega|^{-\frac{1}{r}}
			\left(
			1 + \| (a,u) \|_{ {\mathscr{E}_{\Omega,\varepsilon}}(t)}
			\right)\| (a,u) \|_{ {\mathscr{E}_{\Omega,\varepsilon}}(t)}^2\\
			&\quad 
			+
			C|\Omega|^{-\frac{1}{r}}
			\left(
			\| (a,u) \|_{ {\mathscr{A}_{\Omega,\varepsilon;\alpha_\delta}^{q,r}}(t)}^2
			+
			\| (a,u) \|_{ {\mathscr{A}_{\Omega,\varepsilon;\alpha_\delta}^{q,r}}(t)}^{\frac{r}{r-1}}
			\| (a,u) \|_{ {\mathscr{E}_{\Omega,\varepsilon}}(t)}^{\frac{r-2}{r-1}}
			\right).
		\end{split}
	\end{align}
	Lemmas \ref{lemm:middle-str-1} and \ref{lemm:AE} yield
	\begin{align}\label{pf:main-A-3}
			\notag		
			\| (a,u) \|_{L^r( 0,t ; \dB_{q,\infty}^{\frac{3}{q}-3+\frac{4}{r}} )}^{m;|\Omega|\varepsilon,\alpha_{\delta}} 
			& \leq {}
			C\delta
			+
			C\alpha_{\delta}^{\frac{1}{r}}|\Omega|^{-\frac{1}{r}}
			\left(
			\varepsilon \alpha_{\delta}
			\| (a,u) \|_{ {\mathscr{E}_{\Omega,\varepsilon}}(t)}
			+
			\| (a,u) \|_{ {\mathscr{A}_{\Omega,\varepsilon;\alpha_\delta}^{q,r}}(t)}
			\right)^2\\
			&\quad 
			+
			C\alpha_{\delta}^{\frac{1}{r}}|\Omega|^{-\frac{1}{r}}
			\| (a,u) \|_{ {\mathscr{A}_{\Omega,\varepsilon;\alpha_\delta}^{q,r}}(t)}^{\frac{r}{r-1}}
			\| (a,u) \|_{ {\mathscr{E}_{\Omega,\varepsilon}}(t)}^{\frac{r-2}{r-1}}.
	\end{align}
	It follows from Lemma \ref{lemm:middle-ene-1}, \ref{lemm:a-propri-high} {(}with $p=q${)}, and \ref{lemm:AE} that
	\begin{align}\label{pf:main-A-4}
		\begin{split}
			&
			\| (a,u) \|_{\widetilde{L^{\infty}}( 0,t ; \dB_{q,1}^{\frac{3}{q}-1} ) \cap L^1( 0,t ; \dB_{q,1}^{\frac{3}{q}+1} )}^{m;\alpha_{\delta},\frac{\beta_0}{\varepsilon}}\\
			&
			+
			\varepsilon \| a \|_{\widetilde{L^{\infty}}( 0,t ; \dB_{p,1}^{\frac{3}{q}})}^{h;\frac{\beta_0}{\varepsilon}}
			+
			\frac{1}{\varepsilon} \| a \|_{L^1( 0,t ; \dB_{q,1}^{\frac{3}{q}})}^{h;\frac{\beta_0}{\varepsilon}}
			+
			\| u \|_{\widetilde{L^{\infty}}( 0,t ; \dB_{q,1}^{\frac{3}{q}-1}) \cap L^1( 0,t ; \dB_{q,1}^{\frac{3}{q}+1})}^{h;\frac{\beta_0}{\varepsilon}}\\
			&\quad
			\leq
			C
			\delta
			+
			C
			\left(
			\varepsilon \alpha_{\delta}
			\| (a,u) \|_{ {\mathscr{E}_{\Omega,\varepsilon}}(t)}
			+
			\| (a,u) \|_{ {\mathscr{A}_{\Omega,\varepsilon;\alpha_\delta}^{q,r}}(t)}
			\right)^2\\
			&\qquad
			+
			C
			\| (a,u) \|_{ {\mathscr{A}_{\Omega,\varepsilon;\alpha_\delta}^{q,r}}(t)}^{\frac{r}{r-1}}
			\| (a,u) \|_{ {\mathscr{E}_{\Omega,\varepsilon}}(t)}^{\frac{r-2}{r-1}}.
		\end{split}
	\end{align}
	Thus, gathering \eqref{pf:main-A-1}, \eqref{pf:main-A-2}, \eqref{pf:main-A-3}, and \eqref{pf:main-A-4}, we obtain
	\begin{align}\label{A}
		\begin{split}
			&\| (a,u) \|_{ {\mathscr{A}_{\Omega,\varepsilon;\alpha_\delta}^{q,r}}(t)}
			\leq 
			C_2
			\delta+\\
			&\quad 
			+
			C_2
			\left( 
			1 
			+
			\varepsilon \alpha_{\delta}
			\| (a,u) \|_{ {\mathscr{E}_{\Omega,\varepsilon}}(t)}
			+
			\| (a,u) \|_{ {\mathscr{A}_{\Omega,\varepsilon;\alpha_\delta}^{q,r}}(t)}
			\right)\\
			&\quad \qquad 
			\times 
			\left(
			\| (a,u) \|_{ {\mathscr{A}_{\Omega,\varepsilon;\alpha_\delta}^{q,r}}(t)}^2
			+
			\| (a,u) \|_{ {\mathscr{A}_{\Omega,\varepsilon;\alpha_\delta}^{q,r}}(t)}^{\frac{r}{r-1}}
			\| (a,u) \|_{ {\mathscr{E}_{\Omega,\varepsilon}}(t)}^{\frac{r-2}{r-1}}
			\right)\\
			&\quad+
			C_2
			|\Omega|^{-\frac{1}{r}}
			\left(
			1 + \| (a,u) \|_{ {\mathscr{E}_{\Omega,\varepsilon}}(t)}
			\right)\| (a,u) \|_{ {\mathscr{E}_{\Omega,\varepsilon}}(t)}^2\\
			&\quad 
			+
			C_2
			|\Omega|^{-\frac{1}{r}}
			\left(
			\| (a,u) \|_{ {\mathscr{A}_{\Omega,\varepsilon;\alpha_\delta}^{q,r}}(t)}^2
			+
			\| (a,u) \|_{ {\mathscr{A}_{\Omega,\varepsilon;\alpha_\delta}^{q,r}}(t)}^{\frac{r}{r-1}}
			\| (a,u) \|_{ {\mathscr{E}_{\Omega,\varepsilon}}(t)}^{\frac{r-2}{r-1}}
			\right)\\
			&\quad 
			+
			C_2
			\left(
			\varepsilon \alpha_{\delta}
			\| (a,u) \|_{ {\mathscr{E}_{\Omega,\varepsilon}}(t)}
			+
			\| (a,u) \|_{ {\mathscr{A}_{\Omega,\varepsilon;\alpha_\delta}^{q,r}}(t)}
			\right)^2\\
			&\quad
			+
			C_2
			\| (a,u) \|_{ {\mathscr{A}_{\Omega,\varepsilon;\alpha_\delta}^{q,r}}(t)}^{\frac{r}{r-1}}
			\| (a,u) \|_{ {\mathscr{E}_{\Omega,\varepsilon}}(t)}^{\frac{r-2}{r-1}}
		\end{split}
	\end{align}
	for some positive constant $C_2=C_2(\mu,P,q,r)$.
	\vskip1pc
	\noindent
	{\it Step.3{:} The {continuity} argument.}
	Now, we prove $T_{\Omega,\varepsilon}^{\rm max} = \infty$ by the continuous argument via the a priori estimates for the $ {\mathscr{E}_{\Omega,\varepsilon}}(t)$ and $ {\mathscr{A}_{\Omega,\varepsilon;\alpha_\delta}^{q,r}}(t)$ norms established in Step.1 and Step.2.
	To this end, we consider the time
	\begin{align}
		T_{\Omega,\varepsilon}^*
		:=
		\sup
		\left\{
		t \in (0,T_{\Omega,\varepsilon}^{\rm max}{)}\ ;\ 
		\begin{aligned}
			&
			\| (a,u) \|_{ {\mathscr{E}_{\Omega,\varepsilon}}(t)}
			\leq
			2C_3
			\mathscr{D}_{\varepsilon}[a_0,u_0],\\
			&
			\| (a,u) \|_{ {\mathscr{A}_{\Omega,\varepsilon;\alpha_\delta}^{q,r}}(t)}
			\leq 
			2C_3\delta
		\end{aligned}
		\right\},
	\end{align}
	where $C_3:=\max\{C_1,C_2\}$.
	Suppose by contradiction that $T_{\Omega,\varepsilon}^*<\infty$.
	Then, by  {Proposition \ref{lemm:LWP}}, we see that $0< T_{\Omega,\varepsilon}^* < T_{\Omega,\varepsilon}^{\rm max}$.
	It follows from \eqref{E} and \eqref{A} that for $0<t < T_{\Omega,\varepsilon}^*$
	\begin{align}
		\| (a,u) \|_{ {\mathscr{E}_{\Omega,\varepsilon}}(t)}
		\leq{}& 
		C_3
		\mathscr{D}_{\varepsilon}[a_0,u_0]\\
		&
		+
		4C_3^2
		\left(
		1 + 2C_3\varepsilon\alpha_{\delta}\mathscr{D}_{\varepsilon}[a_0,u_0]
		\right)
		\left(
		\delta^2
		+
		\delta^{\frac{r}{r-1}}
		\mathscr{D}_{\varepsilon}[a_0,u_0]^{\frac{r-2}{r-1}}
		\right)\\
		& 
		+
		4C_3^3
		\left(
		\varepsilon
		\alpha_{\delta}
		\mathscr{D}_{\varepsilon}[a_0,u_0]
		+
		\delta  
		\right)^2\\
		& 
		+
		4C_3^3
		\left(
		\varepsilon
		\alpha_{\delta}
		\mathscr{D}_{\varepsilon}[a_0,u_0]
		+
		\delta  
		\right)
		\mathscr{D}_{\varepsilon}[a_0,u_0]\\
		&
		+
		4C_3^3\varepsilon 
		\mathscr{D}_{\varepsilon}[a_0,u_0]^2
		+
		4C_3^3
		\delta^{\frac{r}{2(r-1)}}
		\mathscr{D}_{\varepsilon}[a_0,u_0]^{\frac{r-2}{2(r-1)}+1}\\
		\leq{}&
		C_3
		\mathscr{D}_{\varepsilon}[a_0,u_0]\\
		&
		+
		4C_3^2
		\left(
		1 + 2C_3\delta
		\right)
		\left(
		\delta^2
		+
		\delta^{\frac{r}{r-1}}
		\mathscr{D}_{\varepsilon}[a_0,u_0]^{\frac{r-2}{r-1}}
		\right)\\
		& 
		+
		16C_3^3
		\delta^2
		+
		12C_3^3
		\delta  
		\mathscr{D}_{\varepsilon}[a_0,u_0]
		+
		4C_3^3
		\delta^{\frac{r}{2(r-1)}}
		\mathscr{D}_{\varepsilon}[a_0,u_0]^{\frac{r-2}{2(r-1)}+1}
	\end{align}
	and
	\begin{align}
		\| (a,u) \|_{ {\mathscr{A}_{\Omega,\varepsilon;\alpha_\delta}^{q,r}}(t)}
		\leq{}&
		C_3\delta
		+
		4C_3^3
		\left(
		1 + 2C_3\varepsilon\alpha_{\delta}\mathscr{D}_{\varepsilon}[a_0,u_0] + 2C_3 \delta
		\right)\\
		&\qquad \qquad \times
		\left(
		\delta^2
		+
		\delta^{\frac{r}{r-1}}
		\mathscr{D}_{\varepsilon}[a_0,u_0]^{\frac{r-2}{r-1}}
		\right)\\
		&
		+
		4C_3^3|\Omega|^{-\frac{1}{r}}
		\left( 1 + 2C_3 \mathscr{D}_{\varepsilon}[a_0,u_0] \right)\mathscr{D}_{\varepsilon}[a_0,u_0]^2\\
		&
		+
		4C_3^3|\Omega|^{-\frac{1}{r}}\left(
		\delta^2
		+
		\delta^{\frac{r}{r-1}}
		\mathscr{D}_{\varepsilon}[a_0,u_0]^{\frac{r-2}{r-1}}
		\right)\\
		&
		+
		4C_3^3
		\left(
		\varepsilon \alpha_{\delta} \mathscr{D}_{\varepsilon}[a_0,u_0]
		+
		\delta
		\right)^2\\
		&
		+
		4C_3^3
		\delta^{\frac{r}{r-1}}
		\mathscr{D}_{\varepsilon}[a_0,u_0]^{\frac{r-2}{r-1}}\\
		\leq{}&
		C_3\delta 
		+
		4C_3^3(1 + 4C_3\delta)
		\left(
		\delta^2
		+
		\delta^{\frac{r}{r-1}}
		\mathscr{D}_{\varepsilon}[a_0,u_0]^{\frac{r-2}{r-1}}
		\right)\\
		&
		+
		4C_3^3|\Omega|^{-\frac{1}{r}}
		\left( 1 + 2C_3 \mathscr{D}_{\varepsilon}[a_0,u_0] \right)\mathscr{D}_{\varepsilon}[a_0,u_0]^2\\
		&
		+
		4C_3^3|\Omega|^{-\frac{1}{r}}\left(
		\delta^2
		+
		\delta^{\frac{r}{r-1}}
		\mathscr{D}_{\varepsilon}[a_0,u_0]^{\frac{r-2}{r-1}}
		\right)\\
		&
		+
		16C_3^3
		\delta^2
		+
		4C_3^3
		\delta^{\frac{r}{r-1}}
		\mathscr{D}_{\varepsilon}[a_0,u_0]^{\frac{r-2}{r-1}}.
	\end{align}
	Here, there holds
	\begin{gather}
		\frac{1}{C_4}
		\| (a_0,u_0) \|_{\dB_{2,1}^{\frac{1}{2}}}
		\leq 
		\mathscr{D}_{\varepsilon}[a_0,u_0]
		\leq 
		C_4
		\mathscr{D}[a_0,u_0]
	\end{gather}
	for some positive constant $C_4=C_4(\mu)$,
     {
    where we have set 
    \begin{align}
        \mathscr{D}[a_0,u_0]
		:=
		\| a_0 \|_{\dB_{2,\infty}^{-\frac{3}{2}}\cap\dB_{2,1}^{\frac{3}{2}}}
		+
		\| u_0 \|_{\dB_{2,\infty}^{-\frac{3}{2}}\cap\dB_{2,1}^{\frac{1}{2}}}
		+
		\| a_0 \|_{\dB_{2,1}^{\frac{3}{2}}}
		\| a_0 \|_{\dB_{2,\infty}^{-\frac{3}{2}}}.
    \end{align}
    }
	Now, we choose the positive constant $\delta$ so that 
	$\delta \leq \| (a_0,u_0) \|_{\dB_{2,1}^{\frac{1}{2}}}$ and
	\begin{align}
		&
		4C_3^2C_4
		\left(
		1 + 2C_3\delta
		\right)
		\left(
		\delta
		+
		\delta^{\frac{1}{r-1}}
		(C_4\mathscr{D}[a_0,u_0])^{\frac{r-2}{r-1}}
		\right)\\
		&\quad
		+
		28C_3^3C_4
		\delta
		+
		4C_3^3
		\delta^{\frac{r}{2(r-1)}}
		(C_4\mathscr{D}[a_0,u_0])^{\frac{r-2}{2(r-1)}}
		\leq 
		\frac{1}{2},\\
		&
		4C_3^3(1 + 4C_3\delta)
		\left(
		\delta
		+
		\delta^{\frac{1}{r-1}}
		(C_4\mathscr{D}[a_0,u_0])^{\frac{r-2}{r-1}}
		\right)\\
		&\quad 
		+
		16C_3^3
		\delta
		+
		4C_3^3
		\delta^{\frac{1}{r-1}}
		(C_4\mathscr{D}[a_0,u_0])^{\frac{r-2}{r-1}}
		\leq \frac{1}{6}.
	\end{align}
	For this $\delta$, we choose a constant $\widetilde{\Omega_{\delta}} \geqslant \Omega_{\delta}$ so that
	\begin{align}
		&
		4C_3^3C_4|\Omega|^{-\frac{1}{r}}
		\left( 1 + 2C_3C_4 \mathscr{D}[a_0,u_0] \right)\mathscr{D}[a_0,u_0]^2
		\leq
		\frac{\delta}{6},\\
		&
		4C_3^3|\Omega|^{-\frac{1}{r}}
		\left(
		\delta
		+
		\delta^{\frac{1}{r-1}}
		(C_4\mathscr{D}[a_0,u_0])^{\frac{r-2}{r-1}}
		\right)
		\leq 
		\frac{1}{6}
	\end{align}
	for all $|\Omega|\geqslant \widetilde{\Omega_{\delta}}$.
	Then, we obtain
	\begin{align}
		\| (a,u) \|_{ {\mathscr{E}_{\Omega,\varepsilon}}(t)}
		\leq{}& 
		C_3
		\mathscr{D}_{\varepsilon}[a_0,u_0]\\
		&
		+
		4C_3^2
		\left(
		1 + 2C_3\delta
		\right)
		\left(
		\delta^2
		+
		\delta^{\frac{r}{r-1}}
		(C_4\mathscr{D}[a_0,u_0])^{\frac{r-2}{r-1}}
		\right)\\
		& 
		+
		16C_3^3
		\delta^2
		+
		12C_3^3C_4
		\delta  
		\mathscr{D}_{\varepsilon}[a_0,u_0]\\
		&
		+
		4C_3^3
		\delta^{\frac{r}{2(r-1)}}
		(C_4\mathscr{D}[a_0,u_0])^{\frac{r-2}{2(r-1)}}
		\mathscr{D}_{\varepsilon}[a_0,u_0]\\
		\leq{}& 
		C_3
		\mathscr{D}_{\varepsilon}[a_0,u_0]\\
		&
		+
		4C_3^2C_4
		\left(
		1 + 2C_3\delta
		\right)
		\left(
		\delta
		+
		\delta^{\frac{1}{r-1}}
		(C_4\mathscr{D}[a_0,u_0])^{\frac{r-2}{r-1}}
		\right)
		\mathscr{D}_{\varepsilon}[a_0,u_0]\\
		& 
		+
		28C_3^3C_4
		\delta  
		\mathscr{D}_{\varepsilon}[a_0,u_0]\\
		&
		+
		4C_3^3
		\delta^{\frac{r}{2(r-1)}}
		(C_4\mathscr{D}[a_0,u_0])^{\frac{r-2}{2(r-1)}}
		\mathscr{D}_{\varepsilon}[a_0,u_0]\\
		\leq{}& 
		\frac{3}{2}C_3
		\mathscr{D}_{\varepsilon}[a_0,u_0]
	\end{align}
	and
	\begin{align}
		\| (a,u) \|_{ {\mathscr{A}_{\Omega,\varepsilon;\alpha_\delta}^{q,r}}(t)}
		\leq{}&
		C_3\delta 
		+
		4C_3^3(1 + 4C_3\delta)
		\left(
		\delta
		+
		\delta^{\frac{1}{r-1}}
		(C_4\mathscr{D}[a_0,u_0])^{\frac{r-2}{r-1}}
		\right)\delta\\
		&
		+
		4C_3^3C_4|\Omega|^{-\frac{1}{r}}
		\left( 1 + 2C_3C_4 \mathscr{D}[a_0,u_0] \right)\mathscr{D}[a_0,u_0]^2\\
		&
		+
		4C_3^3|\Omega|^{-\frac{1}{r}}
		\left(
		\delta
		+
		\delta^{\frac{1}{r-1}}
		(C_4\mathscr{D}[a_0,u_0])^{\frac{r-2}{r-1}}
		\right)\delta\\
		&
		+
		16C_3^3
		\delta^2
		+
		4C_3^3
		\delta^{\frac{r}{r-1}}
		(C_4\mathscr{D}[a_0,u_0])^{\frac{r-2}{r-1}}\\
		\leq{}&
		\frac{3}{2}C_3\delta,
	\end{align}
	which contradict to the definition of $T_{\Omega,\varepsilon}^*$.
	Hence, we have $T_{\Omega,\varepsilon}^{\rm max} = T_{\Omega,\varepsilon}^* = \infty$ {and complete the proof.}
\end{proof}

\subsection*{Data availability}
Data sharing {is} not applicable to this article as no datasets were generated or analyzed during the current study.

\subsection*{Conflict of interest statement}
The author declares no conflicts of interest.

\subsection*{Acknowledgments}
The authors deeply appreciate Professor {Shuichi} Kawashima for valuable comments on Lemma \ref{lemm:ene-low}.
The first author was supported by Grant-in-Aid for Research Activity Start-up, Grant Number JP23K19011.
The second author was partly supported by JSPS KAKENHI, Grant Number JP21K13826.
    
\appendix    
{
\section{Key Lemmas}\label{sec:a}
We provide several estimates in the Chemin--Lerner spaces.
    \begin{lemm}\label{lemm:prod-1}
	Let $I \subset \mathbb{R}$ be a time interval.
	Let $s_1,s_2 \in \mathbb{R}$ and $1 \leq p_1,p_2,r,r_1,r_2 \leq \infty$ satisfy
	\begin{gather}
		{s_1 + s_2 \geqslant \max\Mp{0,3\sp{\frac{1}{p_1}+\frac{1}{p_2}-1}}}, \quad
		s_1 \leq \frac{3}{p_1}, \quad
		s_2 < \min\Mp{\frac{3}{p_1},\frac{3}{p_2}},
	\end{gather}
    and
    \begin{align}
        \frac{1}{r}=\frac{1}{r_1} + \frac{1}{r_2}.
    \end{align}
	Then, there exists a positive constant $C=C(s_1,s_2,p_1,p_2)$ such that
	\begin{align}
		\| fg \|_{\widetilde{L^r}(I; \dB_{p_2, \infty}^{s_1 + s_2 - \frac{3}{p_1}}) }
		\leq{}&
		C 
		\| f \|_{ \widetilde{L^{r_1}}(I; \dB_{p_1, 1}^{s_1} ) }
		\| g \|_{ \widetilde{L^{r_2}}(I; \dB_{p_2, \infty}^{s_2} ) }
	\end{align}
	for all $f \in \widetilde{L^{r_1}}(I; \dB_{p_1, 1}^{s_1}( \mathbb{R}^3 ) )$ and $g \in \widetilde{L^{r_2}}(I; \dB_{p_2, \infty}^{s_2}(\mathbb{R}^3))$.
\end{lemm}
\begin{proof}
	We first apply the para-product decomposition on $fg$ as 
	\begin{align}
		fg = T_fg + R(f,g) + T_gf,
	\end{align}
	where $T_fg$ and $R(f,g)$ are defined by
	\begin{align}
		T_fg := \sum_{j \in \mathbb{Z}} \sum_{k \leq j-3}\Delta_k f  \Delta_jg,
		\qquad
		R(f,g) := \sum_{j \in \mathbb{Z}} \sum_{|k-j|\leq 2} \Delta_{j}f \Delta_kg.
	\end{align}
	It suffices to consider the estimate of $R(f,g)$ in the case of $1/p_1 + 1/p_2 \geqslant 1$ since the other part are exactly same as in \cite{Xin-Xu-21}*{Proposition 5.1}. 
    We see from \cite{Ba-Ch-Da-11}*{Theorem 2.52} and $p_2 \leq p_1'$ that  
	\begin{align}
		\n{R(f,g)}_{\widetilde{L^r}(I; \dB_{p_2, \infty}^{s_1 + s_2 - \frac{3}{p_1}})}
		\leq{}&
		C 
		\n{R(f,g)}_{\widetilde{L^r}(I; \dB_{1, \infty}^{s_1 + s_2 + 3 (1 - \frac{1}{p_1} - \frac{1}{p_2})})}
		\\
		\leq {}&
		C 
		\n{f}_{\widetilde{L^{r_1}}(I;\dB_{p_1,1}^{s_1})}
		\n{g}_{\widetilde{L^{r_2}}(I;\dB_{p_1',\infty}^{s_2 + 3 (\frac{1}{p_1'} - \frac{1}{p_2})})}, 
		\\
		\leq {}&
		C 
		\n{f}_{\widetilde{L^{r_1}}(I;\dB_{p_1,1}^{s_1})}
		\n{g}_{\widetilde{L^{r_2}}(I;\dB_{p_2,\infty}^{s_2})}.
	\end{align}
	Thus, we complete the proof.
\end{proof}
 {The following lemma may be proved along the same argument as in \cite{Fu-24}*{Lemmas 2.6 and 2.9}.}
\begin{lemm}\label{lemm:prod-2}
	Let $I \subset \mathbb{R}$ be a time interval.
	Let
	$1 \leq q, \sigma, r, r_1,r_2,r_3,r_4 \leq \infty$
	and
	$s, s_1, s_2, s_3, s_4 \in \mathbb{R}$
	satisfy 
	\begin{gather}
		2 \leq q \leq 4,\quad
		\frac{1}{r} = \frac{1}{r_1} + \frac{1}{r_2} = \frac{1}{r_3} + \frac{1}{r_4},\\
		s_1, s_4 \leq 3 \left( \frac{2}{q} - \frac{1}{2} \right),\quad
		s = s_1 +s_2 = s_3 + s_4 > 0.
	\end{gather}
	Then, there exists a positive constant $C = C(q,s,s_1,s_2,s_3,s_4,\sigma)$
	such that
	\begin{align}
		\| fg \|_{\widetilde{L^r}(I;\dB_{2,\sigma}^{s - 3( \frac{2}{q} - \frac{1}{2})})}^{\ell,\beta}
		\leq
		C
		\left(
		\| f \|_{ \widetilde{L^{r_1}}(I; \dB_{q,1}^{s_1})}^{\ell,\beta}
		\| g \|_{ \widetilde{L^{r_2}}(I; \dB_{q,\sigma}^{s_2})}^{\ell,4\beta}
		+
		\| f \|_{ \widetilde{L^{r_3}}(I; \dB_{q,\sigma}^{s_3})}
		\| g \|_{ \widetilde{L^{r_4}}(I; \dB_{q,1}^{s_4})}
		\right)
	\end{align}
	for all $0 < \beta \leq \infty$, $f$ and $g$ provided that the right{-}hand side is finite.
	Moreover, if $s_1,s_4 < 3 (2\slash q - 1 \slash2)$,
	then it holds
	\begin{align}
		\| fg \|_{\widetilde{L^{r}}(I; \dB_{2,\sigma}^{s - 3( \frac{2}{q} - \frac{1}{2})})}^{\ell,\beta}
		\leq
		C
		\left(
		\| f \|_{\widetilde{L^{r_1}}(I; \dB_{q,\infty}^{s_1})}^{\ell,\beta}
		\| g \|_{\widetilde{L^{r_2}}(I; \dB_{q,\sigma}^{s_2})}^{\ell,4\beta}
		+
		\| f \|_{\widetilde{L^{r_3}}(I; \dB_{q,\sigma}^{s_3})}
		\| g \|_{\widetilde{L^{r_4}}(I; \dB_{q,\infty}^{s_4})}
		\right).
	\end{align}
\end{lemm}
\begin{lemm}\label{lemm:compo}
	Let $I \subset \mathbb{R}$ be a time interval.
	Let $F$ be a smooth function on some interval $J \subset \mathbb{R}$ with $0\in J$ and $F(0)=0$.
	Let 
	$ 1 \leq p,r,\sigma \leq \infty$
	and 
	$s>0$.
	Then, for each $R>0$ with $(-R,R) \subset J$, there exists a positive constant $C=C(s,p,\sigma,R)$ such that 
	\begin{align}
		\| F(a) \|_{\widetilde{L^r}(I; \dB_{p,\sigma}^s)} 
		\leq 
		C 
		\| a \|_{\widetilde{L^r}(I; \dB_{p,\sigma}^s)}
	\end{align}
	for all $a \in \widetilde{L^r}(I; \dB_{p,\sigma}^s(\mathbb{R}^3)) \cap {L^{\infty}}(I; L^{\infty}(\mathbb{R}^3))$ with $\| a \|_{{L^{\infty}}(I; L^{\infty})} \leq R$.
\end{lemm}
}


\begin{bibdiv}
\begin{biblist}
\bib{Ba-Ch-Da-11}{book}{
	author={Bahouri, Hajer},
	author={Chemin, Jean-Yves},
	author={Danchin, Rapha\"{e}l},
	title={Fourier analysis and nonlinear partial differential equations},
	publisher={Springer, Heidelberg},
	date={2011},
}
\bib{Ba-Ma-Ni-97}{article}{
	author={Babin, A.},
	author={Mahalov, A.},
	author={Nicolaenko, B.},
	title={Regularity and integrability of $3$D Euler and Navier--Stokes
		equations for rotating fluids},
	journal={Asymptot. Anal.},
	volume={15},
	date={1997},
	number={2},
	pages={103-150},
}
\bib{Ba-Ma-Ni-00}{article}{
	author={Babin, A.},
	author={Mahalov, A.},
	author={Nicolaenko, B.},
	title={Global regularity of 3D rotating Navier--Stokes equations for
		resonant domains},
	journal={Appl. Math. Lett.},
	volume={13},
	date={2000},
	number={4},
	pages={51-57},
}
\bib{Ba-Ma-Ni-01}{article}{
	author={Babin, A.},
	author={Mahalov, A.},
	author={Nicolaenko, B.},
	title={3D Navier--Stokes and Euler equations with initial data
		characterized by uniformly large vorticity},
	journal={Indiana Univ. Math. J.},
	volume={50},
	date={2001},
	number={Special Issue},
	pages={1-35},
}
\bib{Bo-Fa-Pr-22}{article}{
	author={Bocchi, Edoardo},
	author={Fanelli, Francesco},
	author={Prange, Christophe},
	title={Anisotropy and stratification effects in the dynamics of fast
		rotating compressible fluids},
	journal={Ann. Inst. H. Poincar\'{e} C Anal. Non Lin\'{e}aire},
	volume={39},
	date={2022},
	number={3},
	pages={647--704},
}
{
\bib{Ch-06}{article}{
   author={Charve, Fr\'ed\'eric},
   title={Global well-posedness and asymptotics for a geophysical fluid
   system},
   journal={Comm. Partial Differential Equations},
   volume={29},
   date={2004},
   number={11-12},
   pages={1919--1940},
}
}
\bib{Ch-Da-10}{article}{
   author={Charve, Fr\'{e}d\'{e}ric},
   author={Danchin, Rapha\"{e}l},
   title={A global existence result for the compressible Navier--Stokes
   equations in the critical $L^p$ framework},
   journal={Arch. Ration. Mech. Anal.},
   volume={198},
   date={2010},
   number={1},
   pages={233--271},
}
{
\bib{Ch-De-Ga-Gr-02}{article}{
	author={Chemin, J.-Y.},
	author={Desjardins, B.},
	author={Gallagher, I.},
	author={Grenier, E.},
	title={Anisotropy and dispersion in rotating fluids},
	book={
		series={Stud. Math. Appl.},
		volume={31},
		publisher={North-Holland, Amsterdam},
	},
	date={2002},
	pages={171-192},
}
\bib{Ch-De-Ga-Gr-06}{book}{
	author={Chemin, J.-Y.},
	author={Desjardins, B.},
	author={Gallagher, I.},
	author={Grenier, E.},
	title={Mathematical geophysics},
	series={Oxford Lecture Series in Mathematics and its Applications},
	volume={32},
	publisher={The Clarendon Press, Oxford University Press, Oxford},
	date={2006},
}
}
\bib{Ch-Mi-Zh-10-1}{article}{
   author={Chen, Qionglei},
   author={Miao, Changxing},
   author={Zhang, Zhifei},
   title={Global well-posedness for compressible Navier-Stokes equations
   with highly oscillating initial velocity},
   journal={Comm. Pure Appl. Math.},
   date={2010},
   pages={1173--1224},
}
\bib{Cu-Ro-Be-11}{book}{
	author={Cushman-Roisin, Benoit},
	author={Beckers, Jean-Marie},
	title={Introduction to geophysical fluid dynamics: physical and numerical aspects},
	publisher={Academic press},
	date={2011},
}
\bib{Da-00}{article}{
	author={Danchin, Rapha\"{e}l},
	title={Global existence in critical spaces for compressible Navier--Stokes
		equations},
	journal={Invent. Math.},
	volume={141},
	date={2000},
	number={3},
	pages={579--614},
}
\bib{Da-02-R}{article}{
   author={Danchin, Rapha\"{e}l},
   title={Zero Mach number limit in critical spaces for compressible
   Navier--Stokes equations},
   journal={Ann. Sci. \'{E}cole Norm. Sup. (4)},
   volume={35},
   date={2002},
   pages={27--75},
}
\bib{Da-14}{article}{
   author={Danchin, Rapha\"{e}l},
   title={A Lagrangian approach for the compressible Navier--Stokes
   equations},
   journal={Ann. Inst. Fourier (Grenoble)},
   volume={64},
   date={2014},
   pages={753--791},
}
\bib{Da-He-16}{article}{
   author={Danchin, Rapha\"{e}l},
   author={He, Lingbing},
   title={The incompressible limit in $L^p$ type critical spaces},
   journal={Math. Ann.},
   volume={366},
   date={2016},
   pages={1365--1402},
}
\bib{Dan-Xu-17}{article}{
   author={Danchin, Rapha\"{e}l},
   author={Xu, Jiang},
   title={Optimal time-decay estimates for the compressible Navier--Stokes
   equations in the critical $L^p$ framework},
   journal={Arch. Ration. Mech. Anal.},
   volume={224},
   date={2017},
   pages={53--90},
}
\bib{Fa-21}{article}{
	author={Fanelli, Francesco},
	title={Incompressible and fast rotation limit for barotropic
		Navier--Stokes equations at large Mach numbers},
	journal={Phys. D},
	volume={428},
	date={2021},
	pages={Paper No. 133049, 20},
}
\bib{Fe-Ga-Ge-No-12}{article}{
	author={Feireisl, Eduard},
	author={Gallagher, Isabelle},
	author={Gerard-Varet, David},
	author={Novotn\'{y}, Anton\'{\i}n},
	title={Multi-scale analysis of compressible viscous and rotating fluids},
	journal={Comm. Math. Phys.},
	volume={314},
	date={2012},
	number={3},
	pages={641-670},
}
\bib{Fe-Ga-No-12}{article}{
	author={Feireisl, Eduard},
	author={Gallagher, Isabelle},
	author={Novotn\'{y}, Anton\'{\i}n},
	title={A singular limit for compressible rotating fluids},
	journal={SIAM J. Math. Anal.},
	volume={44},
	date={2012},
	number={1},
	pages={192-205},
}
\bib{Fe-No-14}{article}{
	author={Feireisl, Eduard},
	author={Novotn\'{y}, Anton\'{\i}n},
	title={Multiple scales and singular limits for compressible rotating
		fluids with general initial data},
	journal={Comm. Partial Differential Equations},
	volume={39},
	date={2014},
	number={6},
	pages={1104-1127},
}
\bib{Fu-24}{article}{
	author={Fujii, Mikihiro},
	title={Low Mach number limit of the global solution to the compressible
		Navier--Stokes system for large data in the critical Besov space},
	journal={Math. Ann.},
	volume={388},
	date={2024},
	number={4},
	pages={4083--4134},
}
{
\bib{Fuj-Wat-25}{article}{
   author={Fujii, Mikihiro},
   author={Watanabe, Keiichi},
   title={Compressible Navier-Stokes-Coriolis system in critical Besov
   spaces},
   journal={J. Differential Equations},
   volume={428},
   date={2025},
   pages={747--795},
}}
\bib{Ha-11}{article}{
   author={Haspot, Boris},
   title={Existence of global strong solutions in critical spaces for
   barotropic viscous fluids},
   journal={Arch. Ration. Mech. Anal.},
   volume={202},
   date={2011},
   number={2},
   pages={427--460},
}
\bib{Iw-Ta-13}{article}{
	author={Iwabuchi, Tsukasa},
	author={Takada, Ryo},
	title={Global solutions for the Navier--Stokes equations in the rotational
		framework},
	journal={Math. Ann.},
	volume={357},
	date={2013},
	number={2},
	pages={727-741},
}
\bib{Iw-Ta-15}{article}{
	author={Iwabuchi, Tsukasa},
	author={Takada, Ryo},
	title={Dispersive effect of the Coriolis force and the local
		well-posedness for the Navier--Stokes equations in the rotational
		framework},
	journal={Funkcial. Ekvac.},
	volume={58},
	date={2015},
	number={3},
	pages={365-385},
}
{
\bib{Ka-95}{book}{
   author={Kato, Tosio},
   title={Perturbation theory for linear operators},
   series={Classics in Mathematics},
   publisher={Springer-Verlag, Berlin},
   date={1995},
   pages={xxii+619},
}}
\bib{Ko-Le-Ta-14-1}{article}{
	author={Koh, Youngwoo},
	author={Lee, Sanghyuk},
	author={Takada, Ryo},
	title={Dispersive estimates for the Navier--Stokes equations in the
		rotational framework},
	journal={Adv. Differential Equations},
	volume={19},
	date={2014},
	number={9-10},
	pages={857-878},
}
\bib{Ma-Ni-79}{article}{
   author={Matsumura, Akitaka},
   author={Nishida, Takaaki},
   title={The initial value problem for the equations of motion of viscous
   and heat-conductive gases},
   journal={J. Math. Kyoto Univ.},
   volume={20},
   date={1980},
   number={1},
   pages={67--104},
}
\bib{Sa-18}{book}{
   author={Sawano, Yoshihiro},
   title={Theory of Besov spaces},
   series={Developments in Mathematics},
   volume={56},
   publisher={Springer, Singapore},
   date={2018},
}
\bib{Xin-Xu-21}{article}{
   author={Xin, Zhouping},
   author={Xu, Jiang},
   title={Optimal decay for the compressible Navier-Stokes equations without
   additional smallness assumptions},
   journal={J. Differential Equations},
   date={2021},
   pages={543--575},
}

\end{biblist}
\end{bibdiv}
\end{document}